\RequirePackage{pdf14}% also disables `\pdfobjcompresslevel` --- needed for kdp, comment for arXiv
 
\documentclass[twoside]{book}
\usepackage{my-book}
\hypersetup{pdftitle={Extra pearls in graph theory}}
\makeindex

\newcommand{\arxiv}[2]{#1} %for arXiv
\arxiv{}{\pdfminorversion=3}

\newcommand{\spell}[2]{#1} %for book

\begin{document}
\spell{}{\pagestyle{empty}\renewcommand\includegraphics[2][{}]{}\def\emph{\textit}\renewcommand\footnote[1]{\ (#1)}\renewcommand\z{}\renewcommand\texorpdfstring[2]{}\renewcommand\section[1]{SECTION. {#1} SECTION.}}

\title{\tt Extra pearls in graph theory}
%\subtitle{\tt A minimalist introduction}
\author{\tt Anton Petrunin}
\date{}
\maketitle

I used these topics together with ``Pearls in graph theory'' by Nora Hartsfield and Gerhard Ringel \cite{hartsfield-ringel} to teach an undergraduate course in graph theory at the Pennsylvania State University.
I tried to keep clarity and simplicity on the same level.

Hope that someone will find it useful for something.

\medskip

I want to thank 
Semyon Alesker,
Alexander Gil,
Rostislav Matveyev,
Alexei Novikov,
Dmitri Panov,
and Lukeria Petrunina for help.

\begin{flushright}
Anton Petrunin
\end{flushright}

{\sloppy

\null\vfill\noindent{\begin{lpic}[t(-0mm),b(-0mm),r(0mm),l(0mm)]{pics/by-sa(0.5)}\end{lpic}
This work is licensed under the Creative Commons Attribution-ShareAlike 4.0 International License.
To view a copy of this license, visit \texttt{http://creativecommons.org/licenses/by-sa/4.0/}}

}

\thispagestyle{empty}
\newpage
\tableofcontents

\chapter{Introduction}

\section{Terminology}

\begin{wrapfigure}{r}{33 mm}
\begin{tikzpicture}[scale=1.4,
  thick,main node/.style={circle,draw,font=\sffamily\bfseries,minimum size=3mm}]

  \node[main node] (1) at (0,15/6) {$a$};
  \node[main node] (2) at (1,15/6){$b$};
  \node[main node] (11) at (1.5,10/6){$c$};
  \node[main node] (12) at (.5,10/6) {$d$};

  \path[every node/.style={font=\sffamily\small}]
   (11) edge [out=20,in=75,looseness=8] node[above] {} (11)
   (1) edge (2)
   (2) edge[bend left] (12)
   (2) edge[bend right] (12)
   (2) edge (11)
   (11) edge (12);
\end{tikzpicture}
\end{wrapfigure}

The diagram on the right may describe regular flights of an airline.
It has six flights,  which serve four airports labeled by $a$, $b$, $c$, and~$d$.

For this and similar types of data, mathematicians use the notion of a \index{pseudograph}\emph{pseudograph}.

Formally, a pseudograph is a finite nonempty set of \index{vertex}\emph{vertices}  (in our example a vertex is an airport) 
and a finite collection of \index{edge}\emph{edges}; each edge connects two vertices (in the example above, an edge is a regular flight).
A pair of vertices can be connected by multiple edges;
such edges are called \index{parallel edges}\emph{parallel}
(in our example, it might indicate that the airline operates multiple flights a day between these airports). 
Also, an edge can connect a vertex to itself; such an edge is called a \index{loop}\emph{loop} (we might think of it as a sightseeing flight).

Thus, from a mathematical point of view, the diagram above describes a pseudograph with vertices $a$, $b$, $c$, $d$, and 
six edges, including one loop at $c$ and a pair of parallel edges between $b$ and $d$.

\smallskip

The number of edges coming from one vertex is called its \index{degree}\emph{degree}; loops are counted twice.
In the example above,
the degrees of $a,b,c,$ and $d$ are $1,4,4,$ and $3$, respectively.

A vertex with zero degree is called \index{isolated vertex}\emph{isolated}, and a vertex of degree one is called an \index{end vertex}\emph{end vertex}.

\smallskip

A pseudograph without loops is also called a \index{multigraph}\emph{multigraph},
and a multigraph without parallel edges is also called a \index{graph}\emph{graph}.
Most of the time we will work with graphs.

If $x$ and $y$ are vertices of a pseudograph $G$, we say that $x$ is \index{adjacent}\emph{adjacent} to $y$ if there is an edge between $x$ and $y$.
We say that a vertex $x$ is \index{incident}\emph{incident} with an edge $e$ if $x$ is an end vertex of $e$.

\section{Wolf, goat, and cabbage}

Тhe vertices and edges of the graph might have very different nature.
As an example, let us consider the following classic puzzle.

\begin{thm}{Puzzle}
A farmer purchased a wolf, a goat, and a cabbage;
he needs to cross a river with them.
He has a boat, but he can carry only himself and a single one of his purchases: the wolf, the goat, or the cabbage.

If left unattended together, the wolf would eat the goat, and the goat would eat the cabbage.

The farmer has to carry himself and his purchases to the far bank of the river, leaving each purchase intact. How can he do it?
\end{thm}

\parit{Solution.}
Let us denote the farmer by $*$, the river by ${\parallel}$
the wolf by $w$, the goat by $g$, and the cabbage by $c$.
For example, $wc{\parallel}{*}g$ means that the wolf and cabbage are on the left bank of the river and the goat with the farmer are on the right bank.

{\sloppy

The starting position is $wgc{*}{\parallel}$; that is, everyone is on the left bank.
The following graph describes all possible positions that can be achieved;
each edge is labeled by the transported purchase.

}

\begin{center}
\begin{tikzpicture}[scale=1.8,
  thick,main node/.style={circle,draw,font=\sffamily\bfseries,minimum size=1mm}]

  \node[main node] (1) at (.5,5/6) {{\small${wgc{*}}{\parallel}{}$}};
  \node[main node] (2) at (1,0){{\small${wc}{\parallel}{{*}g}$}};
  \node[main node] (3) at (2,0){{\small${wc{*}}{\parallel}{g}$}};
  \node[main node] (4) at (2.5,5/6){{\small${w}{\parallel}{{*}gc}$}};
 \node[main node] (5) at (2.5,-5/6){{\small${c}{\parallel}{{*}gw}$}};
  \node[main node] (6) at (3.5,-5/6){{\small${cg{*}}{\parallel}{w}$}};
  \node[main node] (7) at (3.5,5/6){{\small${wg{*}}{\parallel}{c}$}};
  \node[main node] (8) at (4,0){{\small${g}{\parallel}{{*}wc}$}};
  \node[main node] (9) at (5,0){{\small${g{*}}{\parallel}{wc}$}};
  \node[main node] (10) at (5.5,-5/6){{\small${}{\parallel}{{*}wgc}$}};
  \path[every node/.style={font=\sffamily\small}]
   (1) edge node[auto]{$g$}(2)
   (2) edge node{}(3)
   (3) edge node[auto]{$c$}(4)
   (5) edge node[auto]{$w$}(3)
   (6) edge node[auto]{$g$}(5)
   (4) edge node[auto]{$g$}(7)
   (7) edge node[auto]{$w$}(8)
   (8) edge node[auto]{$c$}(6)
   (8) edge node[auto]{}(9)
   (9) edge node[auto]{$g$}(10);
\end{tikzpicture}
\end{center}

This graph shows that the farmer can achieve ${}{\parallel}{{*}wgc}$ by legal moves.
It solves the problem and also shows that there are exactly two different solutions,
assuming that the farmer does not want to repeat the same position twice. 
\qeds

Often, a graph comes with extra structure; for example, labeling of edges and/or vertices as in the example above.

Formally speaking, the described graph representation can be applied to nearly any puzzle.
However, it might be impossible to apply due to a huge number of states and transitions;
think of Rubik's cube or the 15 puzzle.

Here is a small variation of another classic puzzle.

\begin{thm}{Puzzle} Four missionaries and four cannibals must cross a river using a boat that can carry at most two people.
The boat cannot cross the river by itself with no people on board.
If there are missionaries on one bank, they cannot be outnumbered by cannibals;
otherwise the missionaries will be eaten.

Can the missionaries and cannibals safely cross the river?
\end{thm}

Let us introduce a notation to describe the states and transitions in this puzzle.
The river will be denoted by ${\parallel}$;
let $*$ denote the boat, we will write the number of cannibals on each side of ${\parallel}$, and the number of missionaries by subscript. 
For example, $4_2^*{\parallel}0_2$ means that on the left bank we have four cannibals, two missionaries, and the boat (these two missionaries will be eaten), and on the right bank there are no cannibals and two missionaries.

\begin{thm}{Exercise}\label{ex:cannibals}
Draw the corresponding graph.
Conclude that the puzzle has no solutions.
\end{thm}

%\begin{thm}{Exercise}\label{ex:ab+ac+bc}
%Draw the graph of all algebraic formulas that can be reached from $a\cdot b+a\cdot c+b\cdot c$ applying the commutative law of addition and the distributive law.
%\end{thm}

\chapter{Ramsey numbers}

In this chapter, we discuss Ramsey numbers.
The will be used as an illustration of the next chapter.

\section{Ramsey numbers}

Recall that the Ramsey number $r(m,n)$ is the least positive integer such that every blue-red coloring of edges in the complete graph $K_{r(m, n)}$ contains a blue $K_m$ or a red $K_n$.

Switching colors in the definition shows that $r(m,n)=r(n,m)$ for any $m$ and $n$.
Therefore, we may assume that $m\le n$.

Note that $r(1,n)=1$ for any positive integer $n$.
Indeed, the one-vertex graph $K_1$ has no edges;
therefore we can say that all its edges are blue (as well as \textit{red} and \textit{deep green-cyan turquoise} at the same time).

\begin{thm}{Exercise}\label{ex:r(2,n)}
Show that $r(2,n)=n$ for any positive integer $n$.
\end{thm}

The following table from \cite{radziszowski} 
includes all known values of $r(m,n)$ for $n\ge m\ge 3$:

\begin{table}[ht!]\label{ramsey-table}
\centering{%
    \begin{tabular}{|c|*{9}{c|}}
      \hline
      \diagbox[width=.8cm, height=.8cm]{$\!\!m$}{$n\!\!$}
       & 1 & 2 & 3 & 4  & 5  & 6  & 7  & 8  & 9\\
      \hline
      1& 1 & 1 & 1 & 1  & 1  & 1  & 1  & 1  & 1\\
      \hline
      2& 1 & 2 & 3 & 4  & 5  & 6  & 7  & 8  & 9\\
      \hline 
      3& 1 & 3 & 6 & 9  & 14 & 18 & 23 & 28 & 36\\
      \hline
      4& 1 & 4 & 9 & 18 & 25 & ?  & ?  & ?  & ?\\
      \hline
    \end{tabular}
  }%
\end{table}

For example, the table says that $r(4,4)=18$.
In order to show this, one has to prove two inequalities $r(4,4)\z\ge 18$ and $r(4,4)\z\le 18$.
The inequality $r(4,4)\z\ge 18$ means that there is a blue-red coloring of edges of $K_{17}$ that has no monochromatic $K_4$; this will be done in Exercise~\ref{ex:K8+K17}.
The inequality $r(4,4)\le 18$ means that in any blue-red coloring of $K_{18}$ there is a monochromatic~$K_4$.
The latter follows from $r(3,4)\le 9$ via the inequality \ref{eq:ramsey-inq} below, and
a proof of $r(3,4)\le 9$ is given in \cite[4.3]{hartsfield-ringel}.

\section{Binomial coefficients}

In this section, we review properties of binomial coefficients that will be needed further.

\index{binomial coefficient}\emph{Binomial coefficients} will be denoted by $\tbinom{n}{m}$;
they can be defined as unique numbers satisfying the identity
\[(a+b)^n=\tbinom{n}{0}\cdot a^0\cdot b^n+\tbinom{n}{1}\cdot a^1\cdot b^{n-1}+\dots +\tbinom{n}{n}\cdot a^n\cdot b^{0}\eqlbl{eq:binom-thm}
\]
for any real numbers $a,b$ and integer $n\ge 0$.
This identity is called \index{binomial expansion}\emph{binomial expansion}.
It can be used to derive some identities on binomial coefficients; for example, 
\[\tbinom{n}{0}+\tbinom{n}{1}+\dots +\tbinom{n}{n}=(1+1)^n=2^n.\eqlbl{eq:binom-2n}\]

The number $\tbinom{n}{m}$ plays an important role in combinatorics, as
it gives the number of ways that $m$ objects can be chosen from $n$ different objects.
This value can be found explicitly using the formula
\[\tbinom nm=\frac{n!}{m!\cdot (n-m)!}.\]

Note that all $\tbinom{n}{m}$ different ways to choose $m$ objects from $n$ different objects fall into two categories: (1) those which include the last object --- there are $\tbinom{n-1}{m-1}$ of them, and (2) those which do not include it --- there are $\tbinom{n-1}{m}$ of them.
It follows that 
\[\tbinom{n}{m}=\tbinom{n-1}{m-1}+\tbinom{n-1}{m}.\eqlbl{eq:binomial}\]
This identity will be used in the proof of Theorem~\ref{thm:ramsey-up}.

\section{Upper bound}

Recall that according to Theorem 4.3.2 in \cite{hartsfield-ringel}, the inequality
\[r(m,n) \le r(m-1, n) + r(m, n-1)\eqlbl{eq:ramsey-inq}\]
holds for all integers $m,n\ge 2$.

In other words, any value $r(m,n)$ in the table above cannot exceed the sum of values in the cells directly above and on the left from it.
The inequality \ref{eq:ramsey-inq} might be strict; for example,
\[r(3,4)=9<4+6=r(2,4)+r(3,3).\]

\begin{thm}{Theorem}\label{thm:ramsey-up}
For any positive integers $m,n$ we have that  
\[r(m,n)\le \tbinom{m+n-2}{m-1}.\]
\end{thm}

\parit{Proof.}
Set 
\[s(m,n)=\tbinom{m+n-2}{m-1}=\tfrac{(m+n-2)!}{(m-1)!\cdot(n-1)!},\]
so we need to prove the following inequality: 
\[r(m,n)\le s(m,n).\eqlbl{eq:r<s}\]
Note that from \ref{eq:binomial}, we get the identity
\[s(m,n)=s(m-1,n)+s(m,n-1)\eqlbl{eq:binomial-s}\]
which is similar to the inequality \ref{eq:ramsey-inq}.

Further note that $s(1,n)=s(n,1)=1$ for any positive integer $n$.
Indeed, $s(1,n)=\tbinom{n-1}{0}$, and there is only one choice of $0$ objects from the given $n-1$.
Similarly $s(n,1)=\tbinom{n-1}{n-1}$, and there is only one choice of $n-1$ objects from the given $n-1$.

The above observations make it possible to calculate the values of $s(m,n)$ recursively.
The following table provides some of its
\begin{table}[ht!]
\centering{%
    \begin{tabular}{|c|*{9}{c|}}
      \hline
      \diagbox[width=.8cm, height=.8cm]{$\!\!m$}{$n\!\!$}
       & 1 & 2 & 3 & 4  & 5  & 6  & 7  & 8  & 9\\
      \hline
      1& 1 & 1 & 1 & 1  & 1  & 1  & 1  & 1  & 1\\
      \hline
      2& 1 & 2 & 3 & 4  & 5  & 6  & 7  & 8  & 9\\
      \hline 
      3& 1 & 3 & 6 & 10 & 15 & 21 & 28 & 36 & 45\\
      \hline
      4& 1 & 4 & 10& 20 & 35 & 56 & 84  & 120  & 165\\
      \hline
    \end{tabular}
  }%
\end{table}
 values.
The inequality \ref{eq:r<s} means that any value in this table cannot exceed the corresponding value in the table for $r(m,n)$ on page~\pageref{ramsey-table}. 
The latter is nearly evident from \ref{eq:ramsey-inq} and \ref{eq:binomial-s};
let us show it formally.

Since
\[r(1,n)=r(n,1)=s(1,n)=s(n,1)=1,\]
the inequality \ref{eq:r<s} holds if $m=1$ or $n=1$.

Assume the inequality \ref{eq:r<s} does not hold for some $m$ and $n$.
Choose a \index{minimal criminal}\emph{minimal criminal} pair $(m,n)$;
that is, a pair with minimal value $m+n$ such that \ref{eq:r<s} does not hold.
From above we have that $m,n\ge2$.
Since $m+n$ is minimal, we have that
\[r(m-1,n)\le s(m-1,n)\quad \text{and}\quad r(m,n-1)\le s(m,n-1)\]
summing these two inequalities and applying \ref{eq:ramsey-inq} together with \ref{eq:binomial-s}
we get \ref{eq:r<s} --- a contradiction.
\qeds

\begin{thm}{Corollary}\label{cor:4^n}
The inequality
\[r(n,n)\le \tfrac14\cdot 4^n\] 
holds for any positive integer $n$.
\end{thm}

\parit{Proof.}
By \ref{eq:binom-2n}, we have that 
$\tbinom{k}{m}\le2^k$.
Applying Theorem~\ref{thm:ramsey-up}, we get that
\begin{align*}
r(n,n)&\le \tbinom{2\cdot n-2}{n-1}\le
\\
&\le2^{2\cdot n-2}=
\\
&=\tfrac14\cdot 4^n.
\end{align*}
\qedsf

\section{Lower bound}

In order to show that 
\[r(m,n)\ge s+1,\] 
it is sufficient to color the edges of $K_s$ in red and blue so that it has no red $K_m$ and no blue $K_n$.
Equivalently, it is sufficient to decompose $K_s$ into two subgraphs with no isomorphic copies of $K_m$ in the first one and no isomorphic copies of $K_n$ in the second one.

\begin{wrapfigure}{o}{38mm}
\centering
\begin{lpic}[t(-2 mm),b(0 mm),r(0 mm),l(0 mm)]{mppics/pic-21(1)}
\end{lpic}
\bigskip
\begin{lpic}[t(-0 mm),b(0 mm),r(0 mm),l(0 mm)]{mppics/pic-22(1)}
\end{lpic}
\end{wrapfigure}

For example, the subgraphs in the decomposition of $K_5$ in the diagram have no monochromatic triangles;
this implies that $r(3,3)\ge 6$.
We already showed that for any decomposition of $K_6$ into two subgraphs,
one of the subgraphs has a triangle;
that is, $r(3,3)=6$.

Similarly, to show that $r(3,4)\ge 9$, we need to construct a decomposition of $K_{8}$ into two subgraphs $G$ and $H$ such that $G$ contains no triangle $K_3$ and $H$ contains no  $K_4$.
In fact, in any decomposition of $K_9$ into two subgraphs,
either the first subgraph contains a triangle or the second contains a $K_4$.
That is, $r(3,4)=9$ \cite[see][p. 82--83]{hartsfield-ringel}.

Applying \ref{eq:ramsey-inq}, we get 
\[r(4,4)\le 2\cdot r(3,4)=18.\]
Therefore, to prove that $r(4,4)= 18$, we need to show that  $r(4,4)\ge 18$;
that is, we need to construct a decomposition of $K_{17}$ into two subgraphs with no $K_4$.

The required decomposition is given on the  diagram.
The constructed decomposition is rationally symmetric; the first subgraph contains the chords of angle lengths 1, 2, 4, and 8 and the second contains all the chords of angle lengths 3, 5, 6, and 7.

\begin{figure}[ht!]
\centering
\begin{lpic}[t(-0 mm),b(0 mm),r(0 mm),l(0 mm)]{mppics/pic-23}
\end{lpic}
\end{figure}

\begin{thm}{Exercise}\label{ex:K8+K17}
Show that 

\begin{enumerate}[(a)]
\item In the decomposition of $K_8$ above, the left graph contains no triangle, and the right graph contains no $K_4$.
\item In the decomposition of $K_{17}$, neither graph contains any $K_4$.
\end{enumerate}
\end{thm} 

For larger values $m$ and $n$, the problem of finding the exact lower bound for $r(m,n)$ quickly becomes too hard.
Even getting a reasonable estimate is challenging.
In the next chapter we will show how to obtain such an estimate by using probability.

\chapter{Probabilistic method}

We assume that the reader is familiar with the notion of a random variable and expected value,
at least on an intuitive level.
An introductory part in any textbook on probability should be sufficient;
see for example \cite{feller, lawler, williams}.

\section{Markov's inequality}

A \index{random variable}\emph{random variable} is a real number that depends on a random event.
We will consider random variables that take only finitely many values. 

Two random variables are called \index{independent random variables}\emph{independent} if the occurrence of one does not affect the other.

The \index{expected value}\emph{expected value} is the average of a large number of independently selected outcomes of the random variable.
The expected value of a random variable $X$ will be denoted by $\EE[X]$.
Suppose a random variable $X$ takes only values $x_1,\dots,x_n$ with probabilities $P_1,\dots,P_n$, respectively,
so $P_1+\dots+P_n=1$.
Then 
\[\EE[X]=P_1\cdot x_1+\dots+P_n\cdot x_n.
\eqlbl{eq:E(X)}
\]
For example, if $X$ is the result of rolling a die, then it takes each value $1,2,\dots,6$ with probability $\tfrac16$;
therefore 
\[\EE[X]=\tfrac16\cdot 1+\dots+\tfrac16\cdot6=3.5.\]

\begin{thm}{Claim}\label{clm:E}
For any two random variables $X$ and $Y$ and real constant $c$ we have
\[\EE[X+Y]=\EE[X]+\EE[Y]
\qquad\text{and}\qquad
\EE[c\cdot X]=c\cdot \EE[X].\]
Moreover, if $X$ and $Y$ are independent, then 
\[\EE[X\cdot Y]=\EE[X]\cdot\EE[Y].\]
(We assume that all expected values in the formulas are well defined.)
\end{thm}

\begin{thm}{Markov's inequality}\index{Markov's inequality}
Suppose $Y$ is a nonnegative random variable and $c> 0$.
Denote by $P$ the probability of the event $Y\ge c$.
Then 
\[P\cdot c\le \EE[Y].
\eqlbl{eq:cebyshov}\]
\end{thm}

\parit{Proof.}
Consider another random variable $\bar Y$ such that $\bar Y=c$ if $Y\ge c$ and $\bar Y=0$ otherwise.

Note that $\bar Y\z\le Y$ and therefore
\[\EE[\bar Y]\le \EE[Y].\]
The random variable $\bar Y$ takes the value $c$ with probability $P$ and $0$ with probability $1-P$.
By \ref{eq:E(X)}, 
\[\EE[\bar Y]=P\cdot c;\] whence \ref{eq:cebyshov} follows.
\qeds

\section{Probabilistic method}

The probabilistic method allows us to prove the existence of graphs with certain properties without explicitly constructing them.
The idea is to show that if one randomly chooses a graph or its coloring from a specified class, then the probability that the result has the needed property is more than zero.
This implies that a graph with required the property exists.

Despite using probability, the final conclusion is determined for certain, without any possible error.

\medskip

\begin{thm}{Theorem}\label{thm:ramsey-lower}
Assume that the inequality 
\[\tbinom R n < 2^{{\binom n 2} - 1}\]
holds for a pair of positive integers $R$ and $n$.
Then $r(n,n)>R$.
\end{thm}

\parit{Proof.} 
We need to show that the complete graph $K_R$
admits a coloring of edges in red and blue such that it has no monochromatic subgraph isomorphic to $K_n$.

Let us color the edges randomly;
that is, color each edge independently in red or blue with equal chances.

Fix a set $S$ of $n$ vertices. 
Define the random variable $X(S)$ to be $1$ if every edge between the vertices in $S$ has the same color; otherwise set $X(S)=0$.

The random variable $X(S)$ may take the values $0$ or $1$.
The expected value of $X(S)$ is the probability that $X(S)=1$;
that is, all of the $\tbinom n 2=\tfrac{n\cdot(n-1)}{2}$
edges in $S$ have the same color. 
The probability that all the edges with the ends in $S$ are blue is ${2^{-\binom n 2}}$, and with the same probability all the edges are red.
Since these two possibilities exclude each other, 
\[\EE[X(S)]={2}\cdot {2^{-\binom n 2}}.\]
This holds for any $n$-vertex subset $S$ of the vertices of $K_R$.

Note that the number $Y$ of monochromatic $n$-subgraphs in $K_R$ is the sum of $X(S)$ over all possible $n$-vertex subsets $S$. 
The total number of such subsets is $\tbinom R n$.
By \ref{clm:E}, 
\[\EE[Y]=2\cdot \tbinom R n\cdot 2^{-\binom n 2}.\]

Denote by $P$ the probability that a random coloring of $K_R$ has at least one monochromatic $K_n$'s.
By Markov's inequality, 
\[P\le \EE[Y].\]

Since the number of monochromatic $K_n$'s is an integer, it follows that 
a random coloring of $K_R$ has no monochromatic $K_n$'s with probability at least $1-P$.
If $\EE[Y]<1$, it implies that this probability is positive;
in particular, at least one edge-coloring of $K_R$ has no monochromatic $K_n$. 
That is, if
$\tbinom R n < 2^{\binom n 2 - 1},$
then there is a coloring $K_R$ with no monochromatic $n$-subgraph.
\qeds

The following corollary implies that the function $n\mapsto r(n,n)$ grows at least exponentially. 

\begin{thm}{Corollary}\label{cor:2^n/2}
$r(n, n)> \tfrac1{8}\cdot 2^{\frac{n}{2}}$.
\end{thm}

\parit{Proof.}
Set $R=\lfloor\tfrac1{8}\cdot 2^{\frac{n}{2}}\rfloor$;
that is, $R$ is the largest integer $\le\tfrac1{8}\cdot 2^{\frac{n}{2}}$.

Note that 
\[2^{\binom n 2 - 1}> (2^{\frac{n-3}2})^n\ge R^n.\]
and
\[\tbinom R n=\frac{R\cdot(R-1)\cdots (R-n+1)}{n!}<  R^n.\]

Therefore,  
\[\tbinom R n<2^{\binom n 2 - 1}.\]
By Theorem~\ref{thm:ramsey-lower}, we get that $r(n,n)> R$.
\qeds

\begin{thm}{Exercise}\label{ex:number(ham-cycles)}
Assume the edges of the complete graph $K_{100}$ are colored randomly,
so each edge is colored independently in red or blue with equal chances. 
Show that the expected number of monochromatic Hamiltonian cycles in $K_{100}$ is larger than $10^{125}$.
\end{thm} 

\parbf{Remark.}
One might think that the exercise alone is sufficient to conclude that \textit{most} of the colorings of $K_{100}$ have a monochromatic Hamiltonian cycle.
Let us show that it is not that easy.
(It is still true that the probability of the existence of a monochromatic coloring is close to 1, but the proof requires more work.)

The total number of colorings of $K_{100}$ is $2^{\binom{100}2}>10^{1400}$.
Therefore, in principle, it might happen that $99.99\%$ of the colorings have no monochromatic Hamiltonian cycles and $0.01\%$ of the colorings contain all the monochromatic Hamiltonian cycles.
To keep the expected value above $10^{125}$,
this $0.01\%$ of the colorings should have fewer than $10^{130}$ of monochromatic cycles on average;
the latter does not seem impossible since the total number of Hamiltonian cycles in $K_{100}$ is $99!/2>10^{155}$.

\section{Counting proof}

In this section, we rewrite the proof of Theorem~\ref{thm:ramsey-lower} without the use of probability.
We do this to affirm that the probabilistic method provides a real proof, without any possible error.

In principle,  any probabilistic proof admits such a translation,
but in most cases, the translation is less intuitive. 

\parit{Counting proof of \ref{thm:ramsey-lower}.}
The graph $K_R$ has $\tbinom{R}{2}$ edges.
Each edge can be colored in blue or red;
therefore the total number of different colorings is \[\Omega=2^{\binom{R}{2}}.\]

Fix a subgraph isomorphic to $K_n$ in $K_R$.
Note that this graph is red in $\Omega/2^{\binom n2}$ different colorings
and blue in $\Omega/2^{\binom n2}$ colorings.

There are $\tbinom Rn$ different subgraphs isomorphic to $K_n$ in $K_R$.
Therefore, the total number of monochromatic $K_n$'s in all the colorings 
is 
\[M=\tbinom Rn\cdot\Omega\cdot  2/2^{\binom n2}.\]

If $M<\Omega$, then by the pigeonhole principle,
there is a coloring with no monochromatic $K_n$.
Hence the result.
\qeds

\section[\texorpdfstring{Graph of $n$-cube}{Graph of n-cube}]{Graph of $\bm{n}$-cube}

In this section, we give another classic application of the probabilistic method.

\begin{wrapfigure}{o}{25mm}
\vskip-0mm
\centering
\includegraphics{mppics/pic-24}
\vskip-0mm
\end{wrapfigure}

Let $Q_n$ denote the graph of the $n$-dimensional cube;
$Q_n$ has $2^n$ vertices, each labeled with a sequence of length $n$ consisting  of zeros and ones.
Тwo vertices are adjacent if their labels differ in only one digit.

The graph $Q_4$ is shown on the diagram.
Note that each vertex of $Q_n$ has degree $n$.

Recall that the \index{distance}\emph{distance} between two vertices in a graph $G$ is the length of a shortest path connecting the vertices.

Note that the distance between two vertices in $Q_n$ is the number of different digits in their sequences.
For example, two sequences $01101$ and $11011$ have different digits at positions $1$, $3$, and $4$,
and the corresponding vertices at distance $3$ from each other.

\begin{thm}{Problem}\label{prob:Qn}
Suppose $\ell(n)$ denotes the maximal number of vertices in $Q_n$ at a distance more than $n/3$ from each other.
Then $\ell(n)$ grows exponentially in $n$;
moreover, $\ell(n)\ge 1.05^n$. 
\end{thm}

To solve the problem, one has to construct a set with at least $1.05^n$ vertices in $Q_n$ that are far from each other.
It is challenging to construct such a set explicitly.
Instead, we will show that 
if one randomly selects that many vertices, they are far from each other with positive probability.
To choose a random vertex in $Q_n$, one can toss a coin $n$ times, writing 1 for a head and 0 for a tail, and then select the vertex labeled by the obtained sequence.

The following exercise guides you to a solution of the problem.
The same argument shows that for any coefficient $k<\tfrac12$, the maximal number of vertices in $Q_n$ at a distance larger than $k\cdot n$ from each other grows exponentially in $n$.
According to Exercise~\ref{ex:lin-Qn}, the case $k= \tfrac12$ is very different.

\begin{thm}{Exercise}\label{ex:Qn-dist}
Let $P_n$ denote the probability that two randomly chosen vertices in $Q_n$ at a distance $\le\tfrac n3$ between them.
\begin{enumerate}[(a)]

\item\label{Pn} Use Claim~\ref{clm:coin} to show that 
\[P_n<.95^n.\]

\item\label{kPn} Assume $k$ vertices  $v_1,\dots ,v_k$ in $Q_n$ are fixed and $v$ is a random vertex.
Show that $k\cdot P_n$ is the expected number of vertices $v_i$ that lies at a distance $\le\tfrac n3$ from $v$.
Use Markov's inequality to show that $v$ lies at a distance larger than $\tfrac n3$ from each of $v_i$ with probability at least $1-k\cdot P_n$.

\item\label{ex:Qn-dist:end} Apply \ref{Pn} and \ref{kPn} to show that there are at least $1.05^n$ vertices in $Q_n$ at a distance larger than $\tfrac n3$ from each other.
\end{enumerate}
\end{thm}

\begin{thm}{Claim}\label{clm:coin}
The probability $P_n$ of obtaining fewer than one-third heads after $n$ fair tosses of a coin decays exponentially in $n$;
in fact $P_n<.95^n$ for any $n$.
\end{thm}

\parit{Proof.}
Let us introduce independent $n$ random variables $X_1,\z\dots X_n$;
each $X_i$ returns the number of heads after the $i$-th toss of the coin;
in particular, each $X_i$ takes values $0$ or $1$ with the probability $\tfrac12$ each.
We need to show that the probability $P_n$ of the event $X_1+\dots+X_n\le\tfrac n3$ is less than 
$.95^n$.

Consider the random variable 
\[Y=2^{-X_1-\dots-X_n}.\]

Note that $P_n$ is the probability of the event that $Y\ge 2^{-\frac n3}$.
Also, $Y>0$. 
By Markov's inequality, we get that
\[P_n\cdot 2^{-\frac n3}\le \EE[Y].\]

The random variable $2^{-X_i}$ takes the two values $1$ and $\tfrac12=2^{-1}$ with the probability of $\tfrac12$ each.
Therefore, 
\[\EE[2^{-X_i}]=\tfrac12\cdot 1+\tfrac12\cdot\tfrac12=\tfrac 34.\]

Note that 
\[Y=2^{-X_1}\cdots 2^{-X_n}.\]
Since the random variables $X_i$ are independent, \ref{clm:E} implies that
\[\EE[Y]=\left(\tfrac34\right)^n.\]

It follows that 
\[P_n\le \left(\tfrac34\cdot 2^{\frac13}\right)^n< .95^n.\]
\qedsf

\begin{thm}{Advanced exercise}\label{ex:lin-Qn}
\begin{enumerate}[(a)]
\item \label{ex:lin-Qn:n+1} Show that $Q_n$ contains at most $n+1$ vertices at a distance larger than $\tfrac n2$ from each other.
\item \label{ex:lin-Qn:2n} Show that $Q_n$ contains at most $2\cdot n$ vertices at a distance at least $\tfrac n2$ from each other.
\end{enumerate}
\end{thm}

\section{Remarks}

The existence of Ramsey number $r(m,n)$ for any $m$ and $n$ is the first result in the so-called \textit{Ramsey theory}. 
A typical theorem in this theory states that any large object of a certain type contains a very ordered piece of a given size.
We recommend the book by Matthew Katz and Jan Reimann \cite{katz-reimann} on the subject. 

Corollaries \ref{cor:4^n} and \ref{cor:2^n/2} imply that 
\[\tfrac18\cdot 2^{\frac12\cdot n}\le r(n,n)\le \tfrac14\cdot 2^{2\cdot n}.\]

It is unknown if these inequalities can be essentially improved.%
\footnote{This question might look insignificant at first glance, but it is considered one of the major problems in combinatorics \cite{gowers}.}
More precisely, it is unknown whether there are constants $c>0$ and $\alpha>\tfrac12$ such that the inequality
\[r(n,n)\ge c\cdot 2^{\alpha\cdot n}\]
holds for any $n$.
Similarly, it is unknown whether there are constants $c$ and $\alpha<2$ such that the inequality
\[r(n,n)\le c\cdot 2^{\alpha\cdot n}\]
holds for any $n$.

The probabilistic method was introduced by Paul Erd\H os.
It finds applications in many areas of mathematics, not only in graph theory.

Note that the probabilistic method is nonconstructive;
if it is used to prove the existence of a certain object, then it is still uncontrollably hard to describe a concrete example.

More involved examples of proofs based on the probabilistic method deal with \textit{typical properties} of random graphs.

To describe the concept, let us consider the following \textit{random process} that generates a graph $G_n$ with $n$ vertices.

Fix a positive integer $n$. 
Consider a graph $G_n$ with the vertices labeled by $1,\dots,n$,
where the existence of an edge between every pair of vertices is decided independently by flipping a coin.

Note that the described process depends only on $n$, and as a result, we can get a graph isomorphic to any given graph with $n$ vertices.

\begin{thm}{Exercise}\label{ex:prob(isom)}
Let $H$ be a graph with $n$ vertices.
Denote by $\alpha$ the probability that $G_n$ is isomorphic to $H$.
Show that 
\[1/2^{\binom n2}\le \alpha\le n!/2^{\binom n2}.\]

Suppose $n=10$;
describe graphs $H_1$ and $H_2$ such that these inequalities become equalities.
\end{thm}

Fix a property of a graph (for example, connectedness)
and denote by $\alpha_n$ the probability that $G_n$ has this property.
We say that the property is \index{typical property}\emph{typical} if $\alpha_n\to 1$ as $n\to \infty$.

\begin{thm}{Exercise}\label{ex:diam=2}
Show that random graphs typically have a diameter of 2.
That is, the probability that $G_n$ has a diameter of 2 converges to~1 as $n\to \infty$.
\end{thm}

Note that from the exercise above, it follows that in the described random process, \textit{the random graphs are typically connected}.

\begin{thm}{Exercise}\label{ex:typ(K100)}
Show that random graphs typically have a subgraph isomorphic to $K_{100}$.
That is, the probability that $G_n$ has a subgraph isomorphic to $K_{100}$ converges to~1 as $n\to \infty$.
\end{thm}

The following theorem gives a deeper illustration of the probabilistic method with the use of typical properties;
a proof can be found in \cite[Chapter 44]{aigner-ziegler}.

\begin{thm}{Theorem}
Given positive integers $g$ and $k$, there is a graph $G$ with girth at least $g$ and a chromatic number at least $k$. %??? defs girth and chromatic number
\end{thm}

\chapter{Deletion and contraction}
\label{Deletion and contraction}

\section{Definitions}

\begin{wrapfigure}{r}{40 mm}
\vskip-12mm
\centering
\includegraphics{mppics/pic-31}
\vskip-0mm
\end{wrapfigure}

Let $G$ be a pseudograph with a marked edge~$e$.
Denote by $G- e$ the pseudograph obtained from $G$ by deleting $e$,
and by $G/e$ the pseudograph obtained from $G$ by contracting the edge $e$ to a point; see diagram.

Assume $G$ is a graph; that is, $G$ has no loops and no parallel edges.
In this case, $G- e$ is also a graph.
However, $G/e$ might have parallel edges, but no loops; that is, $G/e$ is a multigraph.

If $G$ is a multigraph, then so is $G- e$.
If the edge $e$ is parallel to $f$ in $G$, then $f$ in $G/e$ becomes a loop; that is, $G/e$ is a pseudograph in general.

\section{Number of spanning trees}

Recall that $s(G)$ denotes the number of spanning trees in the pseudograph $G$.

{
\begin{wrapfigure}{o}{42 mm}
\vskip-4mm
\centering
\includegraphics{mppics/pic-32}
\vskip-0mm
\end{wrapfigure}

An edge $e$ in a connected graph $G$ is called the \index{bridge}\emph{bridge}, if deletion of $e$ makes the graph disconnected;
in this case, the remaining graph has two connected components which are called \index{bank}\emph{banks}.

\begin{thm}{Exercise}\label{ex:bridge}
Assume that the graph $G$ contains a bridge between banks $H_1$ and $H_2$.
Show that
\[s(G)=s(H_1)\cdot s(H_2).\]
\end{thm}

}

\begin{thm}{Deletion-plus-contraction formula}\index{deletion-plus-contraction formula}
\label{thm:deletion-plus-contraction}
Let $e$ be an edge in the pseudograph $G$.
Assume $e$ is not a loop, then the following identity holds
\[s(G)=s(G- e)+s(G/e).\eqlbl{eq:deletion-plus-contraction}\]

\end{thm}

{

\begin{wrapfigure}{r}{43 mm}
\vskip-12mm
\centering
\includegraphics{mppics/pic-33}
\vskip-0mm
\end{wrapfigure}

It is convenient to write the identity \ref{eq:deletion-plus-contraction} using a diagram as on the picture; the arrows point from one multigraph to multigraphs with the same total number of spanning trees;
the edge $e$ is marked in~$G$.

}

\parit{Proof.} Note that the spanning trees of $G$ can be subdivided into two groups:
(1)~those which contain the edge $e$ and (2)~those which do not.
For the trees in the first group, the contraction of $e$ to a point  gives a spanning tree in~$G/e$, while the trees in the second group are also spanning trees in~$G- e$.

Moreover, both of the described correspondences are one-to-one.
Hence the formula follows.
\qeds

A spanning tree cannot have loops, so removing all loops does not change the number of spanning trees.
Let us state it.

\begin{thm}{Claim}
If $e$ is a loop in a pseudograph $G$, then 
\[s(G)=s(G- e).\]

\end{thm}

The proof of the following claim uses the deletion-plus-contraction formula.

\begin{thm}{Claim}
If one removes an end vertex $w$ from a pseudograph $G$, then in the obtained graph $G- w$
the number of spanning trees remains unchanged;
that is,
\[s(G)=s(G- w).\eqlbl{eq:deletion-plus-contraction-w}\]

\end{thm}

\parit{Proof.}
Denote by $e$ the only edge incident to $w$. 
Note that the graph $G- e$ is not connected, since the vertex $w$ is isolated.
In particular,
$s(G- e)=0$.
On the other hand, $G/e=G- w$.
Therefore, the deletion-plus-contraction formula \ref{eq:deletion-plus-contraction} implies~\ref{eq:deletion-plus-contraction-w}.
\qeds

On the diagrams, we may use a two-sided arrow ``$\leftrightarrow$'' for the graphs with equal numbers of the spanning trees.
For example, using the deletion-plus-contraction formula together with the claims, we can draw the following diagram, which in particular implies the following identity:
\[s(G)=2\cdot s(H).\]

\begin{wrapfigure}{o}{58 mm}
\vskip-0mm
\centering
\includegraphics{mppics/pic-34}
\vskip-0mm
\end{wrapfigure}

Note that the deletion-plus-contraction formula gives an algorithm to calculate  the value $s(G)$ for a given pseudograph $G$.
Indeed, for any edge $e$, both graphs $G- e$ and $G/e$ have a smaller number of edges.
That is, the deletion-plus-contraction formula reduces the problem of finding the number of the trees to simpler graphs;
applying this formula a few times we can reduce the question to a collection of graphs with an evident answer for each.
In the next section we will show how it works.

\section{Fans and their relatives}

Recall that \index{Fibonacci numbers}\emph{Fibonacci numbers} $f_n$ are defined using the recursive identity 
$f_{n+1}=f_n+f_{n-1}$
with $f_1=f_2=1$.
The sequence of Fibonacci numbers starts as
\[1,1,2,3,5,8,13,\dots\]

The graphs of the following type are called \label{page:fan}\index{fan}\emph{fans}; 
\begin{figure}[ht!]
\centering
\includegraphics{mppics/pic-35}
\end{figure}
a fan with $n+1$ vertex will be denoted by $F_n$. 

\begin{thm}{Theorem}\label{thm:fans}
$s(F_n)=f_{2\cdot n}$.
\end{thm}

\parit{Proof.}
Applying the deletion-plus-contraction formula, we can draw the following infinite diagram.
\begin{figure}[ht!]
\centering
\includegraphics{mppics/pic-36}
\end{figure}
(We ignore loops and end vertices since they do not change the number of spanning trees.)
In addition to the fans $F_n$ we use its variations $F_n'$, which differ from $F_n$ by an extra parallel edge.

Set $a_n=s(F_n)$ and $a'_n=s(F'_n)$.
From the diagram, we get the following two recursive relations:
\begin{align*}
a_{n+1}&=a'_n+a_n,
\\
a'_n&=a_n+a'_{n-1}.
\end{align*}
That is, in the sequence 
\[a_1,a_1',a_2,a_2',a_3\dots\]
every number starting from $a_2$ is the sum of the previous two.

Further note that $F_1$ has two vertices connected by a unique edge,
and  $F'_1$ has two vertices connected by a pair of parallel edges.
Hence $a_1=1=f_2$ and $a_1'=2=f_3$ and therefore 
\[a_n=f_{2\cdot n}\]
for any $n$.\qeds

\parbf{Comments.}
We can deduce a recursive relation for $a_n$, without using~$a_n'$:
\begin{align*}
a_{n+1}&=a_n'+a_n=
\\
&=2\cdot a_n+a'_{n-1}=
\\
&=3\cdot a_n-a_{n-1}.
\end{align*}
This is a special case of the so-called \index{constant-recursive sequences}\emph{constant-recursive sequences}.
The general term of constant-recursive sequences can be expressed by a closed formula ---
read \cite{jordan} if you wonder how.
In our case, it is
\[a_n=\tfrac1{\sqrt{5}}\cdot
\left(
(\tfrac{3+\sqrt{5}}2)^n-(\tfrac{3-\sqrt{5}}2)^n
\right).\]
Since $a_n$ is an integer and $0<\tfrac1{\sqrt{5}}\cdot(\tfrac{3-\sqrt{5}}2)^n<1$ for any $n\ge 1$, a shorter formula can be written
\[a_n
=
\left\lfloor\tfrac1{\sqrt{5}}\cdot(\tfrac{3+\sqrt{5}}2)^n\right\rfloor,\]
where $\lfloor x\rfloor$ denotes floor of $x$; 
that is, $\lfloor x\rfloor$ is the maximal integer that does not exceed $x$.

\begin{thm}{Exercise}\label{ex:zig-zag}
Consider the sequence of zig-zag graphs $Z_n$ of the following type:
\begin{figure}[ht!]
\centering
\includegraphics{mppics/pic-37}
\end{figure}

Show that $s(Z_n)=f_{2\cdot n}$ for any $n$. 
\end{thm}

\begin{thm}{Exercise}\label{ex:ladder}
Let us denote by $b_n$ the number of spanning trees in the \label{page:ladder}\index{ladder}\emph{$n$-step ladder} $L_n$; that is, in the graph of the following type:

\begin{figure}[ht!]
\centering
\includegraphics{mppics/pic-38}
\end{figure}

Apply the method we used for the fans $F_n$ to show that the sequence $b_n$ satisfies the following linear recursive relation:
\[b_{n+1}=4\cdot b_n-b_{n-1}.\]

\end{thm}

Note that $b_1=1$ and $b_2=4$; applying the exercise, we could calculate the first numbers of the sequence $(b_n)$:
\[1,4,15,56,209,780,2911,\dots \]

The following exercise is analogous, but more complicated.

\begin{thm}{Advanced exercise}\label{ex:wheel}
Recall that a wheel $W_n$ is the graph of following type:

\begin{figure}[ht!]
\centering
\includegraphics{mppics/pic-30}
\end{figure}

Show that the sequence $c_n=s(W_n)$ satisfies the following recursive relation:
\[c_{n+1}=4\cdot c_n-4\cdot c_{n-1}+c_{n-2}.\]

\end{thm}

Using the exercise above and applying induction, one can show that 
\[c_n=f_{2\cdot n+1}+f_{2\cdot n-1}-2=l_{2\cdot  n}-2\]
for any $n$.
The numbers $l_n=f_{n+1}+f_{n-1}$ are called \index{Lucas numbers}\emph{Lucas numbers};
they pop up in combinatorics as often as Fibonacci numbers.  

\section{Remarks}

The \textit{deletion-plus-contraction formula} together with Kirchhoff's rules was used in the solution of the so-called \index{squaring the square problem}\emph{squaring the square problem}.
The history of this problem and its solution are discussed in a book by Martin Gardiner \cite[Chapter 17]{gardiner}.

The proof of the recurrence relation above is given by Mohammad Hassan Shirdareh Haghighi and Khodakhast Bibak \cite{shirdareh-haghighi-bibak};
this problem is also discussed in the book by Ronald Graham, Donald Knuth, and Oren Patashnik~\cite{graham-knuth-patashnik}.

\chapter{Matrix theorem}

\section{Adjacency matrix}

Let us describe a way to encode a given multigraph $G$ with $p$ vertices by an $p{\times}p$ matrix.
First, enumerate the vertices of the multigraph by numbers from $1$ to $p$;
such a multigraph will be called \index{labeled graph}\emph{labeled}. 
Consider the matrix $A=A_G$ with the component $a_{i,j}$ equal to the number of edges from the $i$-th vertex to the $j$-th vertex of $G$.

This matrix $A$ is called the \index{adjacency matrix}\emph{adjacency matrix} of $G$.
Note that $A$ is \index{symmetric matrix}\emph{symmetric}; that is, $a_{i,j}=a_{j,i}$ for any pair $i,j$.
Also, the diagonal components of $A$ vanish; that is, $a_{i,i}=0$ for any $i$.

{

\begin{wrapfigure}{o}{28 mm}
\begin{tikzpicture}[scale=1.4,
  thick,main node/.style={circle,draw,font=\sffamily\bfseries,minimum size=3mm}]

  \node[main node] (1) at (0,15/6) {$1$};
  \node[main node] (2) at (1,15/6){$2$};
  \node[main node] (11) at (1.5,10/6){$3$};
  \node[main node] (12) at (.5,10/6) {$4$};

  \path[every node/.style={font=\sffamily\small}]
  
   (1) edge node{}(2)
   (1) edge node{}(12)
   (2) edge[bend left] node{}(12)
   (2) edge[bend right] node{}(12)
   (2) edge node{}(11);
\end{tikzpicture}
\end{wrapfigure}

For example, for the labeled multigraph $G$ shown on the diagram, we get the following adjacency matrix:
\[A=\left(
\begin{matrix}
0&1&0&1
\\
1&0&1&2
\\
0&1&0&0
\\
1&2&0&0
\end{matrix}
\right).\]

}

\begin{thm}{Exercise}\label{ex:n(walks)}
Let $A$ be the adjacency matrix of a labeled multigraph.
Show that the components $b_{i,j}$ of the $n$-th power $A^n$ is the number of walks of length $n$ in the graph from vertex $i$ to vertex $j$. 
\end{thm}

\section{Kirchhoff minor}

In this section, we construct a special matrix called the \textit{Kirchhoff minor}, associated with a pseudograph,
and discuss its basic properties.
This matrix will be used in the next section in a formula for the number of spanning trees in a pseudograph~$G$.
Since loops do not change the number of spanning trees, we can remove all of them.
In other words, we can (and will) always assume that $G$ is a multigraph. 

Fix a multigraph $G$ and consider its adjacency matrix $A=A_G$;
it is a $p{\times}p$ symmetric matrix with zeros on the diagonal.

\begin{enumerate}
\item Revert the signs of the components of $A$ and exchange the zeros on the diagonal for the degrees of the corresponding vertices. 
The obtained matrix $A'$ is called the \index{Kirchhoff matrix}\emph{Kirchhoff matrix}, {}\emph{Laplacian matrix} or {}\emph{admittance matrix} of the graph $G$.

\item Delete from $A'$ the last column and the last row;
the obtained matrix $M=M_G$ will be called the \index{Kirchhoff minor}\emph{Kirchhoff minor} of the labeled pseudograph $G$.
\end{enumerate}

{

\begin{wrapfigure}{r}{28 mm}
\vskip-4mm
\begin{tikzpicture}[scale=1.4,
  thick,main node/.style={circle,draw,font=\sffamily\bfseries,minimum size=3mm}]

  \node[main node] (1) at (0,15/6) {$1$};
  \node[main node] (2) at (1,15/6){$2$};
  \node[main node] (11) at (1.5,10/6){$3$};
  \node[main node] (12) at (.5,10/6) {$4$};

  \path[every node/.style={font=\sffamily\small}]
  
   (1) edge node{}(2)
   (1) edge node{}(12)
   (2) edge[bend left] node{}(12)
   (2) edge[bend right] node{}(12)
   (2) edge node{}(11);
\end{tikzpicture}
\end{wrapfigure}

For example, the labeled multigraph $G$ on the diagram has the following Kirchhoff matrix and Kirchhoff minor:
\[A'=\left(
\begin{matrix}
2&-1&0&-1
\\
-1&4&-1&-2
\\
0&-1&1&0
\\
-1&-2&0&3
\end{matrix}
\right),
\quad 
M=\left(
\begin{matrix}
2&-1&0
\\
-1&4&-1
\\
0&-1&1
\end{matrix}
\right).\]

}

\begin{thm}{Exercise}\label{ex:Kirchhoff-row}
Show that in any Kirchhoff matrix $A'$ the sum of the components in each row or column vanishes.
Conclude that 
\[\det A'=0.\]

\end{thm}

\begin{thm}{Exercise}\label{ex:minor>graph}
Draw a labeled pseudograph with the following Kirchhoff minor:
\[\left(
\begin{matrix}
4&-1&-1&-1&0
\\
-1&4&-1&0&-1
\\
-1&-1&4&-1&-1
\\
-1&0&-1&4&-1
\\
0&-1&-1&-1&4
\end{matrix}
\right).\]
\end{thm}

\begin{thm}{Exercise}\label{ex:sum-kirchhoff}
Show that the sum of all components in every column of the Kirchhoff minor is nonnegative.

Moreover, the sum of all components in the $i$-th column vanishes if and only if the $i$-th vertex is not adjacent to the last vertex.
\end{thm}

\parbf{Relabeling.}
Let us understand what happens with Kirchhoff minor and its determinant as we swap two labels distinct from the last one.

\begin{wrapfigure}[6]{o}{28 mm}
\begin{tikzpicture}[scale=1.4,
  thick,main node/.style={circle,draw,font=\sffamily\bfseries,minimum size=3mm}]

  \node[main node] (1) at (0,15/6) {$1$};
  \node[main node] (2) at (1,15/6){$3$};
  \node[main node] (11) at (1.5,10/6){$2$};
  \node[main node] (12) at (.5,10/6) {$4$};

  \path[every node/.style={font=\sffamily\small}]

   (1) edge node{}(2)
   (1) edge node{}(12)
   (2) edge[bend left] node{}(12)
   (2) edge[bend right] node{}(12)
   (2) edge node{}(11);
\end{tikzpicture}
\end{wrapfigure}

For example, if we swap the labels $2$ and $3$ in the graph above,
we get another labeling shown on the diagram.
Then the corresponding Kirchhoff minor will be 
\[
M'=\left(
\begin{matrix}
2&0&-1
\\
0&1&-1
\\
-1&-1&4
\end{matrix}
\right),
\]
which is obtained from $M$ by swapping columns $2$ and $3$ followed by swapping rows $2$ and $3$.

Note that swapping a pair of columns or rows changes the sign of the determinant.
Therefore, swapping one pair of rows and one pair of columns does not change the determinant.
Summarizing we get the following:

\begin{thm}{Observation}\label{observaiton:swap}
Assume $G$ is a labeled graph with $p$ vertices, and $M_G$ is its Kirchhoff minor.
If we swap two labels $i,j<p$, then the corresponding Kirchhoff minor $M'_G$ can be obtained from $M_G$ by swapping columns $i$ and $j$ and then swapping rows $i$ and $j$.
In particular, 
\[\det M'_G=\det M_G.\]

\end{thm}

\parbf{Deletion and contraction.}
Let us understand what happens with the Kirchhoff minor if we delete or contract an edge in the labeled multigraph.
(If after the contraction of an edge we get loops, we remove it; this way we obtain a multigraph.)

Assume an edge $e$ connects the first and the last vertex of a labeled multigraph $G$ as shown in the following example:
\begin{center}
\begin{tikzpicture}[scale=1.4,
  thick,main node/.style={circle,draw,font=\sffamily\bfseries,minimum size=3mm}]

  \node[main node] (1) at (0,15/6) {$1$};
  \node[main node] (2) at (1,15/6){$2$};
  \node[main node] (11) at (1.5,10/6){$3$};
  \node[main node] (12) at (.5,10/6) {$4$};

  \path[every node/.style={font=\sffamily\small}]
  
   (1) edge node{}(2)
   (12) edge node[auto]{$e$}(1)
   (2) edge[bend left] node{}(12)
   (2) edge[bend right] node{}(12)
   (2) edge node{}(11);
\node[align=center, yshift=2em] (title) 
    at (current bounding box.north)
    {$G$};
\end{tikzpicture}
\hskip5mm
\begin{tikzpicture}[scale=1.4,
  thick,main node/.style={circle,draw,font=\sffamily\bfseries,minimum size=3mm}]

  \node[main node] (1) at (0,15/6) {$1$};
  \node[main node] (2) at (1,15/6){$2$};
  \node[main node] (11) at (1.5,10/6){$3$};
  \node[main node] (12) at (.5,10/6) {$4$};

  \path[every node/.style={font=\sffamily\small}]
  
   (1) edge node{}(2)
   (2) edge[bend left] node{}(12)
   (2) edge[bend right] node{}(12)
   (2) edge node{}(11);
\node[align=center, yshift=2em] (title) 
    at (current bounding box.north)
    {$G-e$};
\end{tikzpicture}
\hskip10mm
\begin{tikzpicture}[scale=1.4,
  thick,main node/.style={circle,draw,font=\sffamily\bfseries,minimum size=3mm}]

  \node[main node] (2) at (1,15/6){$2$};
  \node[main node] (11) at (1.5,10/6){$3$};
  \node[main node] (12) at (.5,10/6) {$4$};

  \path[every node/.style={font=\sffamily\small}]
  
  (2) edge node{}(12)
   (2) edge[bend left] node{}(12)
   (2) edge[bend right] node{}(12)
   (2) edge node{}(11);
   \node[align=center, yshift=2em] (title) 
    at (current bounding box.north)
    {$G/e$};
\end{tikzpicture}
\end{center}

Note that deleting $e$ only reduces the corner component of $M_G$ by one,
while contracting it removes the first row and column.
That is, since 
\[M_G=
\left(
\begin{matrix}
2&-1&0
\\
-1&4&-1
\\
0&-1&1
\end{matrix}
\right),\]
we have
\[M_{G- e}=\left(
\begin{matrix}
1&-1&0
\\
-1&4&-1
\\
0&-1&1
\end{matrix}
\right)
\quad\text{and}\quad
M_{G/e}=\left(
\begin{matrix}
4&-1
\\
-1&1
\end{matrix}
\right).\]

Summarizing the above discussion, we get the following:

\begin{thm}{Observation}\label{observaiton:dpc}
Assume $e$ is an edge of a labeled multigraph $G$ between the first and last vertices and $M_G$ is the Kirchhoff minor of~$G$.
Then 
\begin{enumerate}[(a)]
\item the Kirchhoff minor  $M_{G- e}$ of $G- e$ can be obtained from  $M_G$ by subtracting 1 from the corner element with index (1,1);
\item the Kirchhoff minor $M_{G/e}$ of $G/e$ can be obtained  from  $M_G$ by removing the first row and the first column in $M_G$.
\end{enumerate}

In particular, applying the cofactor expansion of a determinant, we get that
\[\det M_G=\det M_{G- e}+\det M_{G/ e}.\]

\end{thm}

Note that the last formula resembles the deletion-plus-contraction formula.
This observation will be a key to the proof of the matrix theorem; see next section.

\section{Matrix theorem}

\begin{thm}{Matrix theorem}\label{thm:matrix}
Let $M$ be the Kirchhoff minor of a labeled multigraph $G$ with at least two vertices.
Then
\[s(G)=\det M,\eqlbl{eq:matrix-formula}\]
where $s(G)$ denotes the number of spanning trees in $G$.
\end{thm}

\parit{Proof.}
Denote by $d$ the degree of the last vertex in $G$.

Assume $d=0$.
Then $G$ is not connected and therefore $s(G)=0$.
On the other hand, the sum in each row of $M_G$ vanishes (compare to Exercise~\ref{ex:sum-kirchhoff}).
Hence the sum of all columns in $M_G$ vanish;
in particular, the columns in $M_G$ are linearly dependent and hence $\det M_G=0$.
Hence the equality \ref{eq:matrix-formula} holds if $d=0$.

As usual, we denote by $p$ and $q$ the number of vertices and edges in $G$; by the assumption we have that $p\ge 2$.

\begin{wrapfigure}{o}{25 mm}
\vskip0mm
\centering
\includegraphics{mppics/pic-41}
\end{wrapfigure}

Assume $p=2$; that is, $G$ has two vertices and $q$ parallel edges connecting them.
Clearly, $s(G)=q$.
Further note that $M_G=(q)$; that is,  the Kirchhoff minor $M_G$ is a $1{\times}1$ matrix with a single component $q$.
In particular, $\det M_G=q$ and therefore the equality \ref{eq:matrix-formula} holds.

Assume the equality \ref{eq:matrix-formula} does not hold in general;
choose a \index{minimal criminal}\emph{minimal criminal} graph $G$;
that is, a graph that minimize the value $p+q$ among the graphs violating \ref{eq:matrix-formula}.

From above we have that $p>2$ and $d>0$.
Note that we may assume that the first and last vertices of $G$ are adjacent;
otherwise permute pair of labels $1$ and some $j<p$  and apply Observation~\ref{observaiton:swap}.
Denote by $e$ the edge between the first and last vertex.

Note that the total number of vertices and edges in the pseudographs $G- e$ and $G/e$ are smaller than $p+q$.
Therefore, we have that
\begin{align*}
s(G- e)&=\det M_{G- e},
&
s(G/e)&=\det M_{G/e}.
\end{align*}
Applying these two identities together with the deletion-plus-contraction formula 
and Observation~\ref{observaiton:dpc}, we get that
\begin{align*}
s(G)&=s(G- e)+s(G/e)=
\\
&=\det M_{G- e}+\det M_{G/e}=
\\
&=\det M_G;
\end{align*}
that is, the identity \ref{eq:matrix-formula} holds for $G$ --- a contradiction.
\qeds

\begin{thm}{Exercise}\label{ex:K33W6Q3}
Fix a labeling for each of the following graphs, 
find its Kirchhoff minor and use the matrix theorem to find the number of spanning trees.
\begin{enumerate}[(a)]
\item $s(K_{3,3})$;
\item $s(W_6)$;
\item $s(Q_3)$.
\end{enumerate}
(Use \href{https://matrix.reshish.com/determinant.php}{\texttt{https://matrix.reshish.com/determinant.php}}, or any other matrix calculator.)
\end{thm}

\section{Calculation of determinants}

In this section, we recall key properties of the determinant which will be used in the next section.

Let $M$ be an $n{\times}n$-matrix; that is, a table $n{\times}n$, filled with numbers called {}\emph{components} of the matrix.
The determinant $\det M$ is a polynomial of the $n^2$ components of $M$, satisfying the following conditions:
\begin{enumerate}
 \item\label{1} The unit matrix has determinant 1; that is,
\[
\det\left(
\begin{matrix}
1&0&\cdots&0
\\
0&1&\ddots&\vdots
\\
\vdots&\ddots&\ddots&0
\\
0&\cdots&0&1
\end{matrix}
\right)=1.
\]
\item\label{2} If we multiply each component of one of the rows of  by a number~$\lambda$, then for the obtained matrix $M'$, we have
\[\det M'=\lambda\cdot \det M.\]

\item\label{3} If one adds (or subtracts) term-by-term components
in one row to another row, then the obtained matrix $M'$ has the same determinant; that is
\[\det M'= \det M.\]
\end{enumerate}
These three conditions define the determinant in a uniquely. 
We will not give a proof of the statement; it is not evident and not complicated (sooner or later you will have to learn it if it is not done already).

\begin{thm}{Exercise}\label{ex:det}
Show that the following property follows from the properties above.
\end{thm}

\begin{enumerate}[resume]
 \item Interchanging any pair of rows of a matrix multiplies its determinant by $-1$; that is, if a matrix $M'$ is obtained from a matrix $M$ by permuting two of its rows, then 
\[\det M'=-\det M.\]
\end{enumerate}

The determinant of $n\times n$-matrix can be written explicitly as a sum of $n!$ terms. 
For example, 
\[
a_1{\cdot} b_2{\cdot} c_3+a_2{\cdot} b_3{\cdot} c_1+a_3{\cdot} b_1{\cdot} c_2-a_3{\cdot} b_2{\cdot} c_1-a_2{\cdot} b_1{\cdot} c_3-a_1{\cdot} b_3{\cdot} c_2\]
is the determinant of the matrix
\[M=\left(
\begin{matrix}
a_1&b_1&c_1
\\
a_2&b_2&c_2
\\
a_3&b_3&c_3
\end{matrix}
\right).\]
However, the properties described above give a more convenient and faster way to calculate the determinant, especially for larger~$n$.

Let us show it in one example which will be needed in the next section.
\begin{align*}
\det\left(
\begin{matrix}
4&-1&-1&-1
\\
-1&4&-1&-1
\\
-1&-1&4&-1
\\
-1&-1&-1&4
\end{matrix}
\right)
=
&\det\left(
\begin{matrix}
1&1&1&1
\\
-1&4&-1&-1
\\
-1&-1&4&-1
\\
-1&-1&-1&4
\end{matrix}
\right) 
=
\\
=
&\det\left(
\begin{matrix}
1&1&1&1
\\
0&5&0&0
\\
0&0&5&0
\\
0&0&0&5
\end{matrix}
\right)
=
\\
=
5^3\cdot
&\det\left(
\begin{matrix}
1&1&1&1
\\
0&1&0&0
\\
0&0&1&0
\\
0&0&0&1
\end{matrix}
\right)=
\\
=
5^3\cdot&\det\left(
\begin{matrix}
1&0&0&0
\\
0&1&0&0
\\
0&0&1&0
\\
0&0&0&1
\end{matrix}
\right)=
\\
=5^3&.
\end{align*}
Let us describe what we used on each line above:
\begin{enumerate}
\item property \ref{3} three times --- we add to the first row each of the remaining rows;
\item property \ref{3} three times --- we add the first row to each of the remaining three rows;
\item property \ref{2} three times;
\item property \ref{3} three times --- we subtract from the first row the remaining three rows;
\item property \ref{1}.
\end{enumerate}

\section{Cayley formula}

Recall that a \index{complete graph}\emph{complete graph} is a graph where each pair of vertices is connected by an edge;
a complete graph with $p$ vertices is denoted by~$K_p$.

Note that every vertex of $K_p$ has degree $p-1$.
Therefore, the Kirchhoff minor $M=M_{K_p}$ in the matrix formula \ref{eq:matrix-formula} for $K_p$ is the following $(p-1)\times (p-1)$-matrix:

\[
M=\left(
\begin{matrix}
p{-}1&-1&\cdots&-1
\\
-1&p{-}1&\ddots&\vdots
\\
\vdots&\ddots&\ddots&-1
\\
-1&\cdots&-1&p{-}1
\end{matrix}
\right).
\]

Following the steps at the end of the previous section, we get
\begin{align*}
\det\left(
\begin{matrix}
p{-}1&-1&\cdots&-1
\\
-1&p{-}1&\ddots&\vdots
\\
\vdots&\ddots&\ddots&-1
\\
-1&\cdots&-1&p{-}1
\end{matrix}
\right)
&=\det\left(
\begin{matrix}
1&1&\cdots&1
\\
-1&p{-}1&\ddots&\vdots
\\
\vdots&\ddots&\ddots&-1
\\
-1&\cdots&-1&p{-}1
\end{matrix}
\right)
=
\\
&=
\det\left(
\begin{matrix}
1&1&\cdots&1
\\
0&p&\ddots&\vdots
\\
\vdots&\ddots&\ddots&0
\\
0&\cdots&0&p
\end{matrix}
\right)=
\\
&=  p^{p-2}.
\end{align*}
That is,
\[\det M=p^{\,p-2}.\]
Therefore, applying the matrix theorem, we get the following:

\begin{thm}{Cayley formula}\label{thm:cayley}
 \[s(K_p)=p^{p-2};\]
that is, the number of spanning trees in the complete graph $K_p$ is $p^{p-2}$.
\end{thm}

\begin{thm}{Exercise}\label{ex:s(Kp-e)}
Show $s(K_p-e)=(p-2)\cdot p^{n-3}$ for any edge $e$ in~$K_p$. 
\end{thm}

\begin{thm}{Exercise}\label{ex:s(Kmn)}
Use the matrix theorem to show that 
\[s(K_{m,n})=m^{n-1}\cdot n^{m-1}.\]
\end{thm}

\section{Remarks}

There is a strong connection between counting spanning trees of a given graph,
calculations of currents in an electric chain 
and random walks; a good survey is given in the book by Peter Doyle and Laurie Snell \cite{doyle-snell}.
Let us give some examples.

Assume that the graph $G$ describes an electric chain;
each edge has resistance one Ohm,  and a battery is connected to the vertices $a$ and $b$.
Assume that the total current between these vertices is one Ampere.
The following procedure calculates the current $I_e$ along a given edge $e$.

Fix a direction of $e$.
Note that any spanning tree $T$ of $G$ has exactly one the following three properties:
(1) the edge $e$ appears on the (necessarily unique) path from $a$ to $b$ in $T$ with a positive orientation,
(2) the edge $e$ appears on the path from $a$ to $b$ in $T$ with a negative orientation,
(3) the edge $e$ does not appear on the path from $a$ to $b$ in $T$.
Denote by $s_+$, $s_-$, and $s_0$ the number of the trees in each group.
Clearly,
\[s(G)=s_++s_-+s_0.\]
Then  the current $I_e$ can be calculated using the following formula:
\[I_e=\frac{s_+-s_-}{s(G)}\cdot I.\]
This statement can be proved by checking Kirchhoff's rules for the currents calculated by this formula.

There are many other applications of Kirchhoff's rules to graph theory.
For example, they can be used to prove the Euler's formula
\[p-q+r=2,\]
where $p$, $q$, and $r$ denote the number of vertices, edges, and regions in a plane drawing of graphs \cite{levi}.

\chapter{Graph-polynomials}\label{chap:Graph-polynomials}

Counting problems often lead to an organized collection of numbers.
Sometimes it is convenient to consider a polynomial with these numbers as coefficients.
If it is done in a smart way, then the algebraic structure of the obtained polynomial reflects the original combinatorial structure of the graph.

\section{Chromatic polynomial}

Denote by $P_G(x)$ the number of different colorings of the graph $G$ in $x$ colors such that the ends of each edge get different colors.

\begin{thm}{Exercise}\label{ex:PG=PHPH}
Assume that a graph $G$ has exactly two connected components $H_1$ and $H_2$.
Show that 
\[P_G(x)=P_{H_1}(x)\cdot P_{H_2}(x)\]
for any $x$.
\end{thm}

\begin{thm}{Exercise}\label{ex:PWn}
Show that for any integer $n\ge 3$,
\[P_{W_n}(x+1)=(x+1)\cdot P_{C_n}(x),\]
where $W_n$ denotes a wheel with $n$ spokes and $C_n$ is a cycle of length~$n$.
\end{thm}

\begin{thm}{Deletion-minus-contraction formula}\label{thm:deletion-minus-contraction}\index{deletion-minus-contraction formula}
Let $e$ be an edge in a pseudograph $G$.
Then
\[P_G(x)=P_{G- e}(x)-P_{G/e}(x).
\eqlbl{eq:deletion-minus-contraction}\]

\end{thm}

\parit{Proof.}
The valid colorings of $G- e$ can be divided into two groups: 
(1) those where the ends of the edge $e$ get different colors --- these remain valid colorings of $G$ and (2) those where the ends of $e$ get the same color --- each of such colorings corresponds to a unique coloring of $G/e$.
Hence
\[P_{G- e}(x)=P_G(x)+P_{G/e}(x),\]
which is equivalent to the deletion-minus-contraction formula \ref{eq:deletion-minus-contraction}.
\qeds

Note that if the pseudograph $G$ has loops, then $P_G(x)=0$ for any $x$.
Indeed, in a valid coloring, the ends of a loop should get different colors, which is impossible.

The latter can also be proved using the deletion-minus-contraction formula.
Indeed, if $e$ is a loop in $G$, then $G/e=G- e$.
Therefore, $P_{G- e}(x)\z=P_{G/e}(x)$ and
\[P_G(x)=P_{G- e}(x)-P_{G/e}(x) =0.\]

Similarly, removing a parallel edge from a pseudograph $G$ does not change the value $P_G(x)$ for any $x$.
Indeed, if $e$ has a parallel edge $f$, then in $G/e$ the edge $f$ becomes a loop.
Therefore, $P_{G/e}(x)=0$ for any $x$ and by the deletion-minus-contraction formula we get that
\[P_G(x)=P_{G- e}(x).\]
The same identity can be seen directly --- any admissible coloring of $G- e$ is also admissible in $G$ --- since the ends of $f$ get different colors, so does $e$. 

Summarizing the discussion above:
the problem of finding $P_G(x)$ for a pseudograph $G$ can be reduced to the case when $G$ is a graph; that is, $G$ has no loops and no parallel edges.
Indeed, if $G$ has a loop, then $P_G(x)=0$ for all $x$.
Further, removing one of the parallel edges from $G$ does not change $P_G(x)$.

Recall that polynomial $P$ of $x$ is an expression of the following type
\[P(x)=a_0+a_1\cdot x+\dots+a_n\cdot x^n,\]
with constants $a_0,\dots, a_n$, which are called {}\emph{coefficients} of the polynomial.
The coefficient $a_0$ is called the \index{free term}\emph{free term} of the polynomial.
If $a_n\ne 0$, it is called the \index{leading coefficient}\emph{leading coefficient} of $P$;
in this case $n$ is the degree of $P$.
If the leading coefficient is 1, then the polynomial is called \index{monic polynomial}\emph{monic}.

\begin{thm}{Theorem}\label{thm:chromatic-polynomial}
Let $G$ be a pseudograph with $p$ vertices.
Then $P_G(x)$ is a polynomial with integer coefficients and a vanishing free term.

Moreover, if $G$ has a loop, then $P_G(x)\equiv 0$;
otherwise, $P_G(x)$ is monic and has degree $p$.
\end{thm}

Based on this result we can call $P_G(x)$ the \index{chromatic polynomial}\emph{chromatic polynomial} of the graph~$G$.
The deletion-minus-contraction formula will play the central role in the proof.

\parit{Proof.}
As usual, denote by $p$ and $q$ the number of vertices and edges in $G$.
To prove the first part, we will use induction on $q$.

As the base case, consider the \index{null graph}\emph{null graph} \index{$N_p$ (null graph)}$N_p$;
that is, the graph with $p$ vertices and no edges.
Since $N_p$ has no edges, any coloring of $N_p$ is admissible.
We have $x$ choices for each of $n$ vertices therefore
\[P_{N_p}(x)=x^p.\]
In particular, the function $x\mapsto P_{N_p}(x)$ is given by a monic polynomial of degree $p$ with integer coefficients and a vanishing free term.

Assume that the first statement holds for all pseudographs with at most $q-1$ edges.
Fix a pseudograph $G$ with $q$ edges. 
Applying the deletion-minus-contraction formula for some edge $e$ in $G$, we get that
\[P_G(x)=P_{G- e}(x)-P_{G/e}(x).\eqlbl{eq:deleting-contracting}\]
Note that the pseudographs $G- e$ and $G/e$ have $q-1$ edges.
By the induction hypothesis, $P_{G- e}(x)$ and $P_{G/e}(x)$ are polynomials with integer coefficients and vanishing free terms.
Hence \ref{eq:deleting-contracting} implies the same for $P_G(x)$.

If $G$ has a loop, then $P_G(x)=0$, as $G$ has no valid colorings.
It remains to show that if $G$ has no loops, then $P_G(x)$ is a monic polynomial of degree $p$.

Assume that the statement holds for any multigraph $G$ with at most $q-1$ edges and at most $p$ vertices.

Fix a multigraph $G$ with $p$ vertices and $q$ edges.
Note that $G- e$ is a multigraph with $p$ vertices and $q-1$ edges.
By assumption, its chromatic polynomial $P_{G- e}$ is monic of degree $p$.

Further the pseudograph $G/e$ has $p-1$ vertices,
and its chromatic polynomial $P_{G/e}$ either vanishes or has degree $p-1$.
In both cases the difference $P_{G- e}-P_{G/e}$ is a monic polynomial of degree~$p$.
It remains to apply \ref{eq:deleting-contracting}.
\qeds

\begin{thm}{Advanced exercise}\label{ex:PGpqn}
Let $G$ be a graph with $p$ vertices, $q$ edges, and $n$ connected components.
Show that 
\[P_G(x)=x^p-a_{p-1}\cdot x^{p-1}+a_{p-2}\cdot x^{p-2}+\dots+(-1)^{p-n}a_n\cdot x^n,\]
where $a_n,\dots,a_{p-1}$ are positive integers and $a_{p-1}=q$.
\end{thm}

\begin{thm}{Exercise}\label{ex:PTCpFnLn}
Use induction and the deletion-minus-contraction formula to show that 
\begin{enumerate}[(a)]
\item\label{ex:PTCpFnLn:tree} $P_{T}(x)=x\cdot(x-1)^q$ for any tree $T$ with $q$ edges;
\item $P_{C_p}(x)=(x-1)^p+(-1)^p\cdot(x-1)$ for the cycle $C_p$ of length $p$.
\item if $G$ is a graph with $p$ vertices, then $P_G(x)\ge 0$
for any $x\ge p-1$.
\item $P_{F_n}(x)=x\cdot(x-1)\cdot(x-2)^{n-1}$, where $F_n$ denotes the $n$-spine fan, defined on page \pageref{page:fan}.
\item $P_{L_n}(x)=x\cdot(x-1)\cdot(x^2-3\cdot x+3)^{n-1}$, where $L_n$ denotes the $n$-step ladder, defined on page \pageref{page:ladder}.
\end{enumerate}
\end{thm}

\begin{thm}{Exercise}\label{ex:P(tree)}
Show that a graph $G$ is a tree if and only if
\[P_G(x)= x\cdot(x-1)^{p-1}\]
for some positive integer $p$.
\end{thm}

\begin{thm}{Exercise}\label{ex:chrom(K_p)}
Show that 
\[P_{K_p}(x)=x\cdot(x-1)\cdots(x-p+1).\]

\end{thm}

\parbf{Remark.}
Note that for any graph $G$ with $p$ vertices, we have
\[P_{K_p}(x)\le P_G(x)\le P_{N_p}(x)\]
for any $x$.
Since both polynomials
\begin{align*}
P_{K_p}(x)&=x\cdot(x-1)\cdots(x-p+1)&&
\text{and}
&
P_{N_p}(x)&=x^p,
\end{align*}
are monic of degree $p$,
it follows that so is $P_G$.

Hence Exercise~\ref{ex:chrom(K_p)} leads to an alternative way to prove the second statement in Theorem~\ref{thm:chromatic-polynomial}.

\begin{thm}{Exercise}\label{ex:P=nonisom}
Construct a pair of nonisomorphic graphs with equal chromatic polynomials.
\end{thm}

\section{Matching polynomial}

Recall that a \index{matching}\emph{matching} in a graph is a set of edges without common vertices.

Given an integer $n\ge0$,
denote by $m_n=m_n(G)$ the number of matchings with $n$ edges in the graph $G$.

Note that for a graph $G$ with $p$ vertices and $q$ edges, we have 
$m_0(G)=1$, 
$m_1(G)=q$, 
and if $2\cdot n>p$, then $m_n(G)=0$.
The maximal integer $k$ such that $m_k(G)\ne0$ is called the \index{matching number}\emph{matching number} of $G$.
The expression 
\[M_G(x)=m_0+m_1\cdot x+\dots +m_k\cdot x^k\]
is called the \index{matching polynomial}\emph{matching polynomial} of $G$.

The matching polynomial $M_G(x)$ conveniently organizes the numbers $m_n(G)$ so we can work with all of them simultaneously.
For example, the degree of $M_G(x)$ is the matching number of $G$ and
the total number of matchings in $G$ is its value at $1$:  
\[M_G(1)=m_0+m_1+\dots +m_k.\]

\begin{thm}{Exercise}\label{ex:MG=MHMH}
Assume that a graph $G$ has exactly two connected components $H_1$ and $H_2$.
Show that 
\[M_G(x)=M_{H_1}(x)\cdot M_{H_2}(x)\]
for any $x$.
\end{thm}

\begin{thm}{Exercise}\label{ex:matchings}
Show that the values
\[\tfrac12\cdot[M_G(1)+ M_G(-1)]\quad\text{and}\quad\tfrac12\cdot[M_G(1)- M_G(-1)]\]
equal the number of matchings with even and odd numbers of edges, respectively.
\end{thm}

{

\begin{wrapfigure}{r}{33 mm}
\vskip-8mm
\centering
\includegraphics{mppics/pic-51}
\vskip-0mm
\end{wrapfigure}

Assume $e$ is an edge in a graph $G$.
Recall that the graph $G- e$ is obtained by deleting $e$ from $G$.
Let us denote by $G- [e]$ the graph obtained by deleting the vertices of $e$ with all their edges from $G$;
that is, if $e$ connects two vertices $v$ and $w$, then 
\[G- [e]=G- \{v,w\}.\]

}

The following exercise is analogous to the deletion-contraction formulas \ref{thm:deletion-minus-contraction} and \ref{thm:deletion-plus-contraction}.

\begin{thm}{Exercise}\label{ex:deletion-deletion-total}
Let $G$ be a graph.
\begin{enumerate}[(a)]
\item\label{ex:deletion-deletion} Show that
\[M_G(x)=M_{G-e}(x)+x\cdot M_{G- [e]}(x)\]
for any edge $e$ in $G$.

\item\label{ex:deletion-deletion-K} Use part \ref{ex:deletion-deletion} to show that the matching polynomials of complete graphs satisfy the following recursive relation:
\[M_{K_{n+1}}(x)=M_{K_{n}}(x)+n\cdot x\cdot M_{K_{n-1}}(x).\]

\item\label{ex:deletion-deletion-KK} Use \ref{ex:deletion-deletion-K} to calculate $M_{K_n}(x)$ for $1\le n\le 6$.
\end{enumerate}

\end{thm}

\section{Spanning-tree polynomial}

Consider a connected graph $G$ with $p$ vertices;
assume $p\ge 2$.

Let us prepare independent variables $x_1,\dots,x_p$, one for each vertex of $G$.
For each spanning tree $T$ in $G$ consider the monomial 
\[x_1^{d_1-1}\cdots x_p^{d_p-1},\]
where $d_i$ denotes the degree of the $i$-th vertex in $T$.

The tree $T$ has $p-1$ edges and therefore 
$d_1+\dots+d_p=2\cdot(p-1)$.
It follows that the total degree of the monomial is $p-2$.

The sum of these monomials is a polynomial of degree $p-2$ of $p$ variables $x_1,\dots, x_p$.
This polynomial will be called the \index{spanning-tree polynomial}\emph{spanning-tree polynomial} of $G$ and it will be denoted by 
$S_G(x_1,\dots,x_p)$.

\begin{wrapfigure}{o}{25 mm}
\vskip-0mm
\centering
\includegraphics{mppics/pic-52}
\vskip-0mm
\end{wrapfigure}

For example, the graph $G$ shown on the diagram has three spanning trees; each is obtained by deleting one of the edges in the cycle $xyz$.
Abusing notation slightly, let us use the same label for a vertex in $G$ and for the corresponding variable.
The monomial for the tree obtained by deleting the edge $xy$ is $x^2\cdot y\cdot z$.
Indeed, in this tree the vertex $x$ has degree 3, vertices $y$ and $z$ have degree $2$ and the remaining vertices $u$, $v$, and $w$ have degree~1.
The monomials for the other two trees are $x^2\cdot z^2$ and $x^3\cdot z$.
Therefore, 
\[S_G(x,y,z,u,v,w)= x^2\cdot y\cdot z+x^2\cdot z^2+x^3\cdot z.\]
Note that $S_G(x,y,z,u,v,w)$ does not depend on $u$, $v$, and $w$ since the corresponding vertices have degree 1 in any spanning tree of $G$.

Note that $s(G)=S_G(1,\dots,1)$ for any graph $G$; that is, $S_G(1,\dots,1)$ is the total number of spanning trees in $G$.
Indeed, each spanning tree in $G$ contributes one monomial to $S_G$, and each monomial contributes $1$ to the value $S_G(1,\dots,1)$.
The following exercise shows that the polynomial $S_G$ keeps a lot more information about spanning trees in $G$.   

\begin{thm}{Exercise}\label{ex:SG}
Let $S_G(x_1,\dots,x_p)$ be the spanning-tree polynomial of a graph $G$.
Show the following:
\begin{enumerate}[(a)]
\item\label{ex:SG:a} $S_G(0,1,\dots,1)$ is the number of spanning trees with a leaf at the first vertex.

\item\label{ex:SG:b} The coefficient of $S_G$ in front of $x_1\cdots x_{p-2}$ equals the number of paths of length $p-1$ connecting $(p-1)$-th and $p$-th vertices.

\item\label{ex:SG:c} The partial derivative
\[\frac{\partial S_G}{\partial x_1}(0,1,\dots,1)\]
is the number of spanning trees in $G$ with degree 2 at the first vertex.

\item\label{ex:SG:d} The two values 
\[\tfrac12\cdot\left[S_G(1,1,\dots,1)\pm S_G(-1,1,\dots,1)\right]\]
are the numbers of spanning trees in $G$ with odd or even degrees at the first vertex, respectively.
\end{enumerate}
\end{thm}

\begin{thm}{Theorem}\label{thm:spanning-tree-polynomial}
\[S_{K_p}(x_1,\dots,x_p)=(x_1+\dots +x_p)^{p-2},\]
where $K_p$ is the complete graph with $p\ge 2$ vertices.
\end{thm}

Note that the theorem above generalizes the Cayley formula (\ref{thm:cayley}).
Indeed, 
\[s(K_p)=S_{K_p}(1,\dots,1)=
(1+\dots+1)^{p-2}=p^{p-2}.\]

In the proof we will use the following algebraic lemma; its proof is left to the reader.

\begin{thm}{Lemma}\label{lem:polyx}
Let $P(x_1,\dots,x_n)$ be a polynomial.
Assume 
\[P(x_1,\dots,x_{n-1},0)=0.\]
Then $P$ is divisible by $x_n$;
that is, there is a polynomial $Q(x_1,\dots,x_n)$ such that
\[P(x_1,\dots,x_n)=x_n\cdot Q(x_1,\dots,x_n).\]
\end{thm}

\parit{Proof of \ref{thm:spanning-tree-polynomial}.}
Let us apply induction on $p$;
the base case $p=2$ is evident.

Assume that the statement holds for $p-1$; that is,
\[S_{K_{p-1}}(x_1,\dots,x_{p-1})-(x_1+\dots+x_{p-1})^{p-3}=0.
\eqlbl{eq:S-K(n-1)}\]
We need to show that 
\[S_{K_p}(x_1,\dots,x_{p-1},x_p)-(x_1+\dots+x_{p-1}+x_p)^{p-2}=0.
\eqlbl{eq:S-K(n)}\]

First, let us show that the equality holds if $x_p=0$; that is,
\[S_{K_p}(x_1,\dots,x_{p-1},0)-(x_1+\dots+x_{p-1})^{p-2}=0.
\eqlbl{clm:sum-up}\]
Indeed, $S_{K_p}(x_1,\dots,x_{p-1},0)$ is the sum of all monomials in $S_{K_p}$ without~$x_p$.
Each of these monomials corresponds to a spanning tree $T$ in $K_p$ with $d_p=1$; in other words, $T$ has a leaf at $x_p$.%
\footnote{We use $x_i$ as a label of a vertex in $K_p$ and as the corresponding variable.}
Note that the tree $T$ is obtained from another tree $T'$ with the vertices $x_1,\dots,x_{p-1}$ 
by adding an edge from $x_p$ to $x_i$ for some $i<p$.

Note that the monomial in $S_{K_p}$ that corresponds to $T$ equals 
to the product of $x_i$ times the monomial in $S_{K_{p-1}}$
that corresponds to $T'$.
To get the sum of all monomials in $S_{K_p}$ without $x_p$, we need to sum up these products for all $i<p$ and all the monomials in $S_{K_{p-1}}$; this way we get 
\[S_{K_{p-1}}(x_1,\dots,x_{p-1})\cdot(x_1+\dots+x_{p-1}).\]
By \ref{eq:S-K(n-1)}, the latter equals
\[(x_1+\dots+x_{p-1})^{p-2}\]
which implies \ref{clm:sum-up}.

Now assume \ref{eq:S-K(n)} does not hold.
Denote by $P(x_1,\dots,x_p)$ the left-hand side in \ref{eq:S-K(n)}.
Observe that $P$ is a \index{symmetric polynomial}\emph{symmetric} \index{homogeneous polynomial}\emph{homogeneous} polynomial of degree $p\z-2$.
That is, any permutation of values $x_1,\dots, x_p$ does not change $P(x_1,\dots,x_p)$ and each monomial in $P$ has total degree $p-2$.

By Lemma~\ref{lem:polyx} and \ref{clm:sum-up}, we get that $P$ is divisible by $x_p$.
Since $P$ is symmetric it is divisible by each $x_i$;
that is, there is a polynomial $Q(x_1,\dots,x_p)$ such that 
\[P(x_1,\dots,x_p)=x_1\cdots x_p\cdot Q(x_1,\dots,x_p).\]
Since $P\ne 0$, the total degree of $P$ is at least $p$.
But $P$ has degree $p-2$, a contradiction.
\qeds

\begin{thm}{Exercise}\label{ex:S(Kp)}
Show that the number of spanning trees in $K_p$ with degree $k$ at the first vertex equals $\tbinom{p-2}{k-1}\cdot (p-1)^{p-k-1}$.
\end{thm}

\begin{thm}{Exercise}\label{ex:S(Kmn)}
Assume that the vertices of the left part of $K_{m,n}$ have corresponding variables $x_1,\dots,x_m$ and the vertices in the right part have corresponding variables $y_1,\dots,y_n$. 
Show that
\[S_{K_{m,n}}(x_1,\dots,x_m,y_1,\dots,y_n)=(x_1+\dots +x_m)^{n-1}\cdot(y_1+\dots +y_n)^{m-1}.\]
Conclude that $s(K_{m,n})=m^{n-1}\cdot n^{m-1}$.

\end{thm}

\section{Remarks}

Very good expository papers on chromatic polynomials are written by
Ronald Read \cite{read} and Alexandr Evnin \cite{evnin-chnom}. 
Matching polynomials are discussed in a paper by Christopher Godsil and Ivan Gutman \cite{godsil-gutman}.

Our discussion of spanning-tree polynomials is based on a modification by Fedor Petrov \cite{petrov} of the original proof by Arthur Cayley~\cite{cayley}.    

\chapter{Generating functions}\label{Generating functions}

In this chapter, we discuss one connection between graph theory and power series;
we prove the exponential formula (Theorem~\ref{thm:exp-formula})
and discuss its applications.
This material is similar to Chapter~\ref{chap:Graph-polynomials}, but more challenging.

\section{Exponential generating functions}

{\sloppy

The power series 
\[A(x)=a_0+a_1\cdot x+\tfrac12\cdot a_2\cdot x^2+\dots+\tfrac1{n!}\cdot a_n\cdot x^n+\dots\]
is called the \index{exponential generating function}\emph{exponential generating function} of the sequence $a_0,a_1,\dots$

}

If the series $A(x)$ converges in some neighborhood of zero, then it defines a function which remembers all information of the sequence~$a_n$.
The latter follows since 
\[A^{(n)}(0)=a_n;\eqlbl{eq:An=an}\]
that is, the $n$-th derivative of $A(x)$ at $0$ equals $a_n$.

However, without assuming the convergence, we can treat $A(x)$ as a formal power series.
We are about to describe how to add, multiply, take the derivative, and do other operations with formal power series.

\parbf{Sum and product.}
Consider two exponential generating functions
\begin{align*}
A(x)&=a_0+a_1\cdot x+\tfrac12\cdot a_2\cdot x^2+\tfrac16\cdot a_3\cdot x^3+\dots
\\
B(x)&=b_0+b_1\cdot x+\tfrac12\cdot b_2\cdot x^2+\tfrac16\cdot b_3\cdot x^3+\dots
\end{align*}
We will write 
\[S(x)=A(x)+B(x),\quad 
P(x)=A(x)\cdot B(x)\]
if the power series $S(x)$ and $P(x)$ are obtained from $A(x)$ and $B(x)$ by opening the parentheses of these formulas and combining like terms.

It is straightforward to check that $S(x)$ is the exponential generating function for the sequence  
\begin{align*}
s_0&=a_0+b_0,
\\
s_1&=a_1+b_1,
\\
&\dots
\\
s_n&=a_n+b_n.
\end{align*}
The product $P(x)$ is also exponential generating function for the sequence
\[
\begin{aligned}
p_0&=a_0\cdot b_0,
\\
p_1&=a_0\cdot b_1+a_1\cdot b_0,
\\
p_2&=a_0\cdot b_2+2\cdot a_1\cdot b_1+a_2\cdot b_0,
\\
p_3&=a_0\cdot b_3+3\cdot a_1\cdot b_2+3\cdot a_2\cdot b_1+a_3\cdot b_0,
\\
&\dots
\\
p_n&=\sum_{i=0}^n\tbinom ni\cdot a_i\cdot b_{n-i}.
\end{aligned}
\eqlbl{eq:multiplication}
\]

\begin{thm}{Exercise}\label{ex:B=xA}
Assume $A(x)$ is the exponential generating function of the sequence $a_0,a_1,\dots$
Show that $B(x)=x\cdot A(x)$ corresponds to the sequence $b_n=n\cdot a_{n-1}$.
\end{thm}

\parbf{Composition.}
Once we define addition and multiplication of power series we can also substitute one power series into another.
For example, if $a_0=0$, then the expression 
\[E(x)=e^{A(x)}\] is another power series obtained by substituting $A(x)$ for $x$ in
\[e^x=1+x+\tfrac12\cdot x^2+\tfrac16\cdot x^3+\dots\eqlbl{eq:ex}\]
It is harder to express the sequence $e_n$ corresponding to $E(x)$ in terms of $a_n$, but it is easy to find the first few terms.
Since we assume $a_0=0$, we have
\begin{align*}
e_0&=1,
\\
e_1&=a_1,
\\
e_2&=a_2+2\cdot a_1^2,
\\
e_3&=a_3+6\cdot a_1\cdot a_2,
\\
&\dots
\end{align*}

\parbf{Derivative.}
The derivative of 
\[A(x)=a_0+a_1\cdot x+\tfrac12\cdot a_2\cdot x^2+\dots+\tfrac1{n!}\cdot a_n\cdot x^n+\dots\]
is defined as the following formal power series 
\[A'(x)=a_1+a_2\cdot x+\tfrac12\cdot a_3\cdot x^2+\dots+\tfrac1{n!}\cdot a_{n+1}\cdot x^n+\dots\]
Note that $A'(x)$ coincides with the ordinary derivative of $A(x)$ if the latter converges.

Also, $A'(x)$ is the exponential generating function of the sequence $a_1,a_2,a_3,\dots$
obtained from the original sequence $a_0,a_1,a_2,\dots$
by deleting the zero-term and shifting the indexes by 1.

\begin{thm}{Exercise}\label{ex:B=xA'}
Let $A(x)$ be the exponential generating function of the sequence $a_0,a_1,a_2\dots$
Describe the sequence $b_n$ with the exponential generating function
\[B(x)=x\cdot A'(x).\]
\end{thm}

\parbf{Calculus.}
If $A(x)$ converges and
\[E(x)=e^{A(x)},\] 
then we have 
\[\ln E(x)=A(x).\]
Also by taking the derivative of $E(x)=e^{A(x)}$ we get that
\begin{align*}
E'(x)
&=e^{A(x)}\cdot A'(x)=
\\
&= E(x)\cdot A'(x).
\end{align*}

These identities have a perfect meaning in terms of formal power series
and they still hold without convergence.
We will not prove it, but the proof is not difficult.

\section{Fibonacci numbers}

Recall that Fibonacci numbers, denoted as $f_n$, are defined using the recursive identity 
$f_{n+1}=f_n+f_{n-1}$
with initial conditions $f_0=0$ and $f_1=1$.

\begin{thm}{Exercise}\label{ex:exp(Fn)}
Let $F(x)$ be the exponential generating function of Fibonacci numbers $f_n$.
\begin{enumerate}[(a)]
\item\label{ex:exp(Fn):F''} Show that it satisfies the following differential equation
\[F''(x)=F(x)+F'(x).\]
\item\label{ex:exp(Fn):F(x)} Conclude that 
\[F(x)=\frac{1}{\sqrt5}\cdot\left(e^{\frac{1+\sqrt{5}}{2}\cdot x}- e^{\frac{1-\sqrt{5}}{2}\cdot x}\right).\]
\item\label{ex:exp(Fn):Binet} Use the identity \ref{eq:An=an} to derive the \index{Binet's formula}\emph{Binet's formula}: 
\[f_n=\tfrac{1}{\sqrt5}\cdot\left((\tfrac{1+\sqrt{5}}{2})^n-(\tfrac{1+\sqrt{5}}{2})^n\right).\]
\end{enumerate}

\end{thm}

\section{Exponential formula}

Fix a set of graphs $\mathcal{S}$.
Denote by $c_n=c_n(\mathcal{S})$ the number of spanning subgraphs of $K_n$ isomorphic to one of the graphs in $\mathcal{S}$. 
Let \[C(x)=C_{\mathcal{S}}(x)\] be the exponential generating function of the sequence~$c_n$.

\begin{thm}{Theorem}\label{thm:exp-formula}
Let $\mathcal{S}$ be a set of connected graphs. 

\begin{enumerate}[(a)]
\item\label{thm:exp-formula:Wk}
Fix a positive integer $k$.
Denote by $w_n$ the number of spanning subgraphs of $K_n$ with exactly $k$ connected components such that each component is  isomorphic to one of the graphs in~$\mathcal{S}$.
Then
\[W_k(x)=\tfrac1{k!}\cdot C_{\mathcal{S}}(x)^k,\]
where $W_k(x)$ is the exponential generating function of the sequence $w_n$.

\item\label{thm:exp-formula:all}
Let $a_n$ be the number of \textit{all} spanning subgraphs of $K_n$ such that each connected component is isomorphic to one of the graphs in~$\mathcal{S}$.
Let $A(x)$ be the exponential generating function of the sequence $a_n$.
Then
\[1+A(x)=e^{C_{\mathcal{S}}(x)}.\]
\end{enumerate}

\end{thm}

Taking the logarithm and derivative of the formula in \textit{\ref{thm:exp-formula:all}},
we get the following:

\begin{thm}{Corollary}\label{cor:exp-formula}
Assume $A(x)$ and $C(x)$ as in Theorem~\ref{thm:exp-formula}\ref{thm:exp-formula:all}.
Then
\[\ln [1+A(x)]=C(x)\quad\text{and}\quad A'(x)=[1+A(x)]\cdot C'(x).\]
\end{thm}

The second formula in this corollary provides a recursive formula for the corresponding sequences which will be important later.

\parit{Proof of \ref{thm:exp-formula}; \ref{thm:exp-formula:Wk}.}
Denote by $v_n$ the number of spanning subgraphs of $K_n$ which have $k$ ordered connected components and each connected component is isomorphic to one of the graphs in $\mathcal{S}$.
Let $V_k(x)$ be the corresponding generating function.

Note that for each graph described above
there are $k!$ ways to order its $k$ components.
Therefore, $w_n=\tfrac1{k!}\cdot v_n$ for any $n$ and
\[W_k(x)=\tfrac1{k!}\cdot V_k(x).\]
Hence it is sufficient to show that 
\[V_k(x)=C(x)^k.
\eqlbl{eq:V=Ck}\]

To prove the latter identity, we apply induction on $k$ and the multiplication formula \ref{eq:multiplication} for exponential generating functions.
The base case $k=1$ is evident.

Assuming that the identity \ref{eq:V=Ck} holds for $k$;
we need to show that 
\[V_{k+1}=V_k(x)\cdot C(x).\eqlbl{eq:VkC}\]

Assume that a spanning graph with $k+1$ ordered connected components of $K_n$ is given.
Denote by $m$ the number of vertices in the first $k$ components.
There are $\tbinom nm$ ways to choose these $m$ vertices among $n$ vertices of $K_n$. 
For each choice, we have
$v_m$ ways to choose a spanning subgraph with $k$ components in it.
The last component has $m-n$ vertices and we have $c_{n-m}$ ways to choose a subgraph from $\mathcal{S}$.
All together, we get
\[\tbinom nm\cdot v_m\cdot c_n\]
spanning graphs with $k+1$ ordered connected components and $m$ vertices in the first $k$ components.
Summing it up for all $m$, we get the multiplication formula \ref{eq:multiplication} for exponential generating functions; hence \ref{eq:VkC} follows.

\parit{\ref{thm:exp-formula:all}.}
To find $a_n$ we need to add the numbers of spanning graphs with $1, 2,\dots$ components.
That is,
\begin{align*}A(x)&=W_1(x)+W_2(x)+\dots=
\tag{by \ref{thm:exp-formula:Wk}}
\\
&=C(x)+\tfrac12\cdot C(x)^2+\tfrac16\cdot C(x)^3+\dots=
\\
&=e^{C(x)}-1.
\end{align*}
The last equality follows from the Taylor expansion for $e^x$; see \ref{eq:ex}.
\qeds

\section{Perfect matchings}

Recall that a \index{perfect matching}\emph{perfect matching} is a 1-factor of the graph. 
In other words, it is a set of isolated edges that covers all the vertices.
Note that if a graph admits a perfect matching, then the number of its vertices is even.

The \index{double factorial}\emph{double factorial} is defined as the product of all the integers from $1$ up to some non-negative integer $n$ that have the same parity (odd or even) as $n$;
the double factorial of $n$ is denoted by $n!!$.
For example, 
\[9!! = 9\cdot 7 \cdot 5 \cdot 3 \cdot  1 = 945
\quad\text{and}\quad
10!!=10\cdot8\cdot6\cdot4\cdot2=3840.
\]

\begin{thm}{Exercise}\label{ex:perfect-matching}
Let $a_n$ denote the number of perfect matchings in $K_n$.
Show that 
\begin{enumerate}[(a)]
 \item $a_2=1$;
 \item $a_n=0$ for any odd $n$;
 \item\label{ex:perfect-matching:recursion} $a_{n+1}=n\cdot a_{n-1}$ for any integer $n\ge 2$.
 \item\label{ex:perfect-matching:n!!} 
 Conclude that $a_n=0$ and $a_{n+1}=n!!$ for any odd $n$.
\end{enumerate}

\end{thm}

Let us present a proof of the second part of Exercise~\ref{ex:perfect-matching}\ref{ex:perfect-matching:n!!} based on Theorem~\ref{thm:exp-formula}.

\parit{Solution.} 
Denote by $a_n$ the number of perfect matchings in $K_{n}$ and let $A(x)$ be the corresponding exponential generating function.

Note that a perfect matching can be defined as a spanning subgraph such that each connected component is isomorphic to $K_2$.
So we can apply the formula in Theorem~\ref{thm:exp-formula} for the set $\mathcal{S}$ consisting of only one graph $K_2$.

Note that if $K_n$ contains a spanning subgraph isomorphic to $K_2$,
then $n=2$.
It follows that $c_2(\mathcal{S})=1$ and $c_n(\mathcal{S})=0$ for $n\ne 2$.
Therefore, 
\[C(x)=C_{\mathcal{S}}(x)=\tfrac12\cdot x^2.\]

By Theorem~\ref{thm:exp-formula}\ref{thm:exp-formula:all},
\begin{align*}
1+A(x)&=e^{C(x)}=
\\
&=e^{\frac12\cdot x^2}=
\\
&=1+\tfrac12\cdot x^2+\tfrac1{2\cdot 4}\cdot x^4+\tfrac1{6\cdot 8}\cdot x^6+\dots+\tfrac1{n!\cdot 2^n}\cdot x^{2\cdot n}+\dots
\end{align*}
That is,
\[
\tfrac1{(2\cdot n-1)!}\cdot a_{2\cdot n-1}=0
\quad
\text{and}
\quad
\tfrac1{(2\cdot n)!}\cdot a_{2\cdot n}=\tfrac1{n!\cdot 2^n}\]
for any positive integer $n$.
In particular, \begin{align*}
a_{2\cdot n}&=\frac{(2\cdot n)!}{n!\cdot 2^n}=
\\
&=\frac{1\cdot 2\cdots (2\cdot n)}{2\cdot4 \cdots (2\cdot n)}=
\\
&=1\cdot 3\cdots (2\cdot n-1)=
\\
&=(2\cdot n-1)!!
\end{align*}

That is, $a_n=0$ for any odd $n$ and $a_n=(n-1)!!$ for even $n$.
\qeds

\parbf{Remark.}
Note that by Corollary~\ref{cor:exp-formula}, we also have
\[A'(x)=[1+A(x)]\cdot x,\]
which is equivalent to the recursive identity
\[a_{n+1}=n \cdot a_{n-1}\]
in Exercise~\ref{ex:perfect-matching}\ref{ex:perfect-matching:recursion}.

\section{All matchings}

Now let $\mathcal{S}$ be the set of two graphs $K_1$ and $K_2$.
Evidently, $c_1(\mathcal{S})\z=c_2(\mathcal{S})=1$.
Further, we have that $c_n(\mathcal{S})=0$ for all $n\ge 3$ since $K_n$ contains no spanning subgraphs isomorphic to $K_1$ or $K_2$.

Therefore, the exponential generating function of the sequence $c_n(\mathcal{S})$ is a polynomial of degree 2
\[C(x)=x+\tfrac12\cdot x^2.\]

Note that a matching in a graph $G$ can be identified with a spanning subgraph with all connected components isomorphic to  $K_1$ or $K_2$.
If we denote by $a_n$ the number of all matchings and by $A(x)$ the corresponding exponential generating function, then by  Theorem~\ref{thm:exp-formula}\ref{thm:exp-formula:all}, we get that
\[A(x)=e^{x+\frac12\cdot x^2}-1.\]
Applying Corollary~\ref{cor:exp-formula}, we also have
\[A'(x)=[1+A(x)]\cdot (1+x).\]
The latter is equivalent to the following recursive formula for $a_n$:
\[a_{n+1}=a_n+n\cdot a_{n-1}.\eqlbl{an+nan-1}\]
Since $a_1=1$ and $a_2=2$, we can easily find first the few terms of this sequence:
\[1,2,4,10,26,\dots\]

\begin{thm}{Exercise}\label{ex:an+nan-1}
Prove formula \ref{an+nan-1} directly --- without using generating functions.
Compare to Exercise~\ref{ex:deletion-deletion-total}\ref{ex:deletion-deletion-K}.
\end{thm}

\section{Two-factors}

Let $\mathcal{S}$ be the set of all cycles.

Note that a \index{$2$-factor}\emph{$2$-factor} of a graph can be defined as a spanning subgraph with components isomorphic to cycles.
Denote by $a_n$ and $c_n$ the number of $2$-factors and \index{spanning cycle}\emph{spanning cycles}%
\footnote{Also known as \index{Hamiltonian cycle}\emph{Hamiltonian cycles}.} in $K_n$ respectively.
Let $A(x)$ and $C(x)$ be the corresponding exponential generating functions.

\begin{thm}{Exercise}\label{ex:ex:2-factor}
\begin{enumerate}[(a)]
\item\label{ex:2-factor:cn} Show that $c_1=c_2=0$ and 
\[c_n=(n-1)!/2\]
for $n\ge 3$.
In particular,
\[c_1=0, c_2=0, c_3=1, c_4=3, c_5=12, c_6=60.\]
\item\label{ex:2-factor:recursive} Use part \ref{ex:2-factor:cn} and the following identity from \ref{cor:exp-formula}
\[A'(x)=[1+A(x)]\cdot C'(x)\]
to find $a_1,\dots, a_6$.
\item Count the number of 2-factors in $K_1,\dots,K_6$ by hand, and compare it to the result in part \ref{ex:2-factor:recursive}.
\item\label{ex:2-factor:C} Use part \ref{ex:2-factor:cn} to conclude that 
\[C(x)=-\tfrac12\cdot\ln(1-x)-\tfrac12\cdot x-\tfrac14\cdot x^2.\]
\item Use part \ref{ex:2-factor:C} and Theorem~\ref{thm:exp-formula}\ref{thm:exp-formula:all}
to show that
\[A(x)=\frac{1}{ e^{\frac x2+\frac{x^2}4}\cdot\sqrt{1-x}}-1.\]
\end{enumerate}

\end{thm}

\section{Counting spanning forests}

Recall that a \index{forest}\emph{forest} is a graph without cycles.
Assume we want to count the number of spanning forests in $K_n$;
denote its number by $a_n$ and by $c_n$ the number of connected spanning forests. 
That is, $c_n$ is the number of spanning trees in $K_n$.
For example, $K_3$ has the following 7 spanning forests; therefore $a_3\z=7$.
\begin{figure}[H]%{r}{20 mm}
\vskip-0mm
\centering
\includegraphics{mppics/pic-81}
\vskip-0mm
\end{figure}

By Corollary~\ref{cor:exp-formula}, the following identity
\[A'(x)=[1+A(x)]\cdot C'(x)\]
holds for the corresponding exponential generating functions.

According to the Cayley theorem, $c_n=n^{n-2}$;
in particular,
\[c_1=1, c_2=1, c_3=3, c_4=16,\dots\]
Applying the product formula \ref{eq:multiplication}, we can use $c_n$ to calculate $a_n$ recurrently:
\begin{align*}
a_1&=c_1=1,
\\
a_2&=c_2+a_1\cdot c_1=
\\&=1+1\cdot 1=2,
\\
a_3&=c_3+2\cdot a_1\cdot c_2+ a_2\cdot c_1=
\\
&=3+2\cdot1\cdot 1+2\cdot 1=7,
\\
a_4&=c_4+3\cdot a_1\cdot c_3+3\cdot a_2\cdot c_2+a_3\cdot c_1=
\\&=16+3\cdot 1\cdot 3+3\cdot 2\cdot 1+7\cdot 1=38
\\
&\dots
\end{align*}

It is instructive to check by hand that there are exactly $38$ spanning forests in $K_4$.

For the general term of $a_n$, no simple formula is known;
however, the recursive formula above provides a sufficiently fast way to calculate its terms.

\section{Counting connected subgraphs}

Let $a_n$ be the number of all spanning subgraphs of $K_n$ and $c_n$ be the number of connected spanning subgraphs of $K_n$.
All 4 connected spanning subgraphs of $K_3$ are shown on the diagram; therefore $c_3=4$.

\begin{figure}[ht!]
\vskip-0mm
\centering
\includegraphics{mppics/pic-82}
\end{figure}

Assume $A(x)$ and $C(x)$ are the corresponding exponential generating functions.
These two series diverge for all $x\ne 0$;
nevertheless, the formula for formal power series in Theorem~\ref{thm:exp-formula}\ref{thm:exp-formula:all}
still holds, and by Corollary~\ref{cor:exp-formula} we can write
\[A'(x)=[1+A(x)]\cdot C'(x).\]

Note that $a_n=2^{\binom n2}$;
indeed, to describe a subgraph of $K_n$ we can choose any subset of $\tbinom n2$ edges of $K_n$, and $a_n$ is the total number of $\binom n2$ these independent choices.
In particular, the first few terms of $a_n$ are
\[a_1=1,\quad a_2=2,\quad a_3=8,\quad a_4=64,\quad\dots\]

Applying the product formula \ref{eq:multiplication}, we can calculate the first few terms of~$c_n$:
\begin{align*}
c_1&=a_1=1
\\
c_2&=a_2-a_1\cdot c_1=
\\
&=2-1\cdot 1=1,
\\
c_3&=a_3-2\cdot a_1\cdot c_2- a_2\cdot c_1=
\\
&=8-2\cdot1\cdot 1-2\cdot 1=4,
\\
c_4&=a_4-3\cdot a_1\cdot c_3-3\cdot a_2\cdot c_2-1\cdot a_3\cdot c_1=
\\
&=64-3\cdot 1\cdot 4-3\cdot 2\cdot 1-1\cdot 8\cdot 1=38,
\\
&\dots
\end{align*}

Note that in the previous section we found $a_n$ from $c_n$, and now we go in the opposite direction.
No closed formula is known for $c_n$,
but the recursive formula gives a sufficiently good way to calculate it.

\section{Remarks}

The method of generating functions was introduced and widely used by Leonard Euler;
the term \textit{generating function} was coined later by Pierre Laplace.
For more on the subject, we recommend the classic book by Frank Harary and Edgar Palmer \cite{harary-palmer}.

\chapter{Minimum spanning trees}

\section{Optimization problems}

We say that a graph $G$ has \index{weight}\emph{weighted edges} if each edge in $G$ is labeled by a real number, called its {}\emph{weight}.
In this case the {}\emph{weight of a subgraph} $H$ (briefly $\weight(H)$) is defined as the sum of the weights of its edges.

An {}\emph{optimization problem} typically asks to minimize the weight of a subgraph of a certain type.
A classic example is the so-called \index{traveling salesman problem}\emph{traveling salesman problem}.
It asks to find a Hamiltonian cycle with minimal weight (if it exists in the graph).

No fast algorithm is known to solve this problem and it is expected that no fast algorithm exists.
We say \textit{fast} if the required time depends polynomially on the number of vertices in the graph.
The best-known algorithms require exponential time, which is just a bit better than brute force checking of all possible Hamiltonian cycles.

In fact, the traveling salesman problem is a classic example of the so-called \textit{NP-hard problems}; it means that if you find a fast algorithm to solve it, then you solve the \index{P=NP problem}\emph{P=NP problem} --- the most important question in modern mathematics.

The so-called \index{nearest neighbor algorithm}\emph{nearest neighbor algorithm} gives a heuristic way to find a Hamiltonian cycle with small (but not necessarily the smallest) weight.
In the following description, we assume that the graph is complete.
\begin{enumerate}[(i)]
\item Start at an arbitrary vertex.
\item\label{NNA-main} Find the edge with minimal weight connecting the current vertex $v$ and an unvisited vertex $w$ and move to $w$.
\item Repeat step \ref{NNA-main} until all vertices are visited and then return to the original vertex.
This way we walked along a Hamiltonian cycle.
\end{enumerate}

\begin{wrapfigure}{o}{25 mm}
\vskip-0mm
\centering
\includegraphics{mppics/pic-90}
\end{wrapfigure}

The nearest neighbor algorithm is an example of \index{greedy algorithm}\emph{greedy algorithms};
in other words, it chooses the cheapest step at each stage.
The cheap choices at the beginning might force it to make expansive choices at the end,
so it may not find an optimal solution. 
In fact a greedy algorithm might give the worst solution.
For example, if we start from the vertex $a$ in the shown weighted graph, then the nearest neighbor algorithm produces the cycle $abcd$, which is the worst --- its total weight is $11=1+1+5+4$, and the other two  Hamiltonian cycles have weights $9$ and $10$.

{\sloppy

Let us list a few other classic examples of optimization problems.
Unlike the traveling salesman problem, there are fast algorithms which solve them.
\begin{itemize}
\item The \index{shortest path problem}\emph{shortest path problem} asks for a path with minimal weight, connecting two given vertices.
It can be thaut of as optimal driving directions between two locations.
\item The \label{assignment problem}\index{assignment problem}\emph{assignment problem} asks for a matching of a given size with minimal weight, in a bipartite graph.
It can be used to minimize the total cost of work, assuming that available workers charge different prices for each task.
This problem is closely related to the subject of the next chapter.
\item The \index{minimum spanning tree proble}\emph{minimum spanning tree problem} asks for a spanning tree with minimal weight.
For example, we may think of minimizing the cost to build a computer network between given locations.
This problem will be the main subject of the remaining sections of this chapter.
\end{itemize}

}

\section{Neighbors of a spanning tree}

Let $G$ be a connected graph.
Suppose $T$ and $T'$ are spanning trees in~$G$.
If 
\[T'=T-e+e'\]
for some edges $e$ and $e'$,
then we say that $T'$ is a \index{neighbor}\emph{neighbor} of $T$.

\begin{thm}{Exercise}\label{ex:neighbor-trees}
Let $G$ be a connected graph that is not a tree.
Suppose that any two distinct spanning trees in $G$ are neighbors.
Show that $G$ has exactly one cycle.
\end{thm}

\begin{thm}{Theorem}\label{thm:mst-iff}
Let $G$ be a graph with weighted edges.
A spanning tree $T$ in $G$ has minimal weight if and only if
\[\weight(T)\le \weight(T')\]
for any neighbor $T'$ of $T$.
\end{thm}

\begin{thm}{Exchange lemma}
Let $S$ and $T$ be spanning trees in a graph $G$.
Then for any edge $s$ in $S$ that is not in $T$, there is an edge $t$ in $T$, but not in $S$ such that the subgraphs
\[S'=S-s+t\quad\text{and}\quad T'=T+s-t\]
are spanning trees in $G$.
\end{thm}

\parit{Proof.}
Note that $S-s$ has two components, denote them by $S_1$ and $S_2$.

Let $P$ be a path in $T$ that connects the vertices of $s$.
Note that the path $P$ starts in $S_1$ and ends in $S_2$.
Therefore, one of the edges of $P$, say $t$, connects $S_1$ to $S_2$.

\begin{figure}[ht!]
\vskip-0mm
\centering
\includegraphics{mppics/pic-91}
\end{figure}

The subgraph $T+s$ has a cycle created by $P$ and~$s$.
Removing the edge $t$ from this cycle leaves a connected subgraph $T'=T-t+s$.
Evidently, $T'$ and $T$ have the same number of vertices and edges;
since $T$ is a tree, so is $T'$.

Further $S'=S-s+t=S_1+S_2+t$.
Since $t$ connects two connected subgraphs $S_1$ and $S_2$, the resulting subgraph $S'$ is connected.
Again $S'$ and $S$ have the same number of vertices and edges;
since $S$ is a tree, so is $S'$.
\qeds

\parit{Proof of \ref{thm:mst-iff}.} The ``only-if'' part is trivial --- if $T$ has minimal weight among all spanning trees, then, in particular, it has to be minimal among its neighbors. 
It remains to show the ``if'' part.

Suppose that a spanning tree $T$ has the minimal weight among its neighbors.
Consider a minimum weight spanning tree $S$ in $G$;
if there is more than one, then
we assume in addition that $S$ shares with $T$ the maximal number of edges. 

Arguing by contradiction, assume $T$ is not a minimum weight spanning tree;
in this case $T\ne S$.
Then there is an edge $s$ in $S$, but not in~$T$.
Let $t$ be an edge in $T$ provided by the exchange lemma; in particular, 
$S'=S-s+t$ and $T'=T+s-t$ are spanning trees.

Since $T'$ is a neighbor of $T$, we have 
\begin{align*}
\weight(T)&\le \weight(T')=
\\
&=\weight(T)+\weight(s)-\weight(t).
\intertext{Therefore, $\weight(s)\ge \weight(t)$ and}
\weight(S')&=\weight(S)-\weight(s)+\weight(t)\le
\\
&\le\weight(S).
\end{align*}
Since $S$ has minimal weight, we have an equality in the last inequality.
That is, $S'$ is another minimum weight spanning tree in $G$.
By construction, $S'$ has one more common edge with $T$;
the latter contradicts the choice of $S$.
\qeds

\begin{thm}{Exercise}\label{ex:w>2w}
Let $G$ be a graph with weighted edges and  $T$ a minimal weight spanning tree in $G$.

Suppose that for each edge $e$, we change its weight from $w$ to $2^w$.
Show that $T$ remains a minimal weight spanning tree in $G$ with the new weights.
\end{thm}

\section{Kruskal’s algorithm}

Kruskal’s algorithm finds a minimum weight spanning tree in a given connected graph.
If the graph is not connected, the algorithm finds a spanning tree in each component of the graph.
Let us describe it informally.
\begin{enumerate}[(i)]
\item Suppose $G$ is a graph with weighted edges.
Start with the spanning forest $F$ in $G$ formed by all vertices and no edges.
\item\label{Kruskal:main} Remove from $G$ an edge with minimal weight.
If this edge connects different trees in $F$, then add it to $F$.
\item Repeat the procedure while $G$ has any remaining edges.
\end{enumerate}

Note that Kruskal’s algorithm is \index{greedy algorithm}\emph{greedy}; it chooses the cheapest step at each stage.

\begin{thm}{Theorem}\label{thm:kruskal}
Kruskal’s algorithm finds a minimum weight spanning tree in any connected pseudograph with weighted edges.
\end{thm}

\parit{Proof.}
Evidently, the subgraph $F$ remains a forest at all stages.
If at the end, $F$ has more than one component, then the original graph $G$ could not have edges connecting these components;
otherwise, the algorithm would choose it at some stage.
Therefore, the resulting graph is a tree.

Suppose the spanning tree $T$ produced by Kruskal’s algorithm does not have minimal total weight.
Then $T$ does not satisfy Theorem~\ref{thm:mst-iff};
that is, there is a neighbor $T'=T+e'-e$ of $T$ such that 
\begin{align*}
\weight(T)&>\weight(T')=
\\
&=\weight(T)+\weight(e')-\weight(e).
\end{align*}
Therefore, $\weight(e')<\weight(e)$.

Denote by $F$ the forest at the stage right before adding $e$ to it.
Note that $F$ is a subgraph of $T-e$.
Since $e'$ connects different trees in $T-e$,
it connects different trees in $F$.
But $\weight(e)>\weight(e')$, therefore the algorithm should reject $e$ at this stage --- a contradiction.
\qeds

\section{Other algorithms}

Let us describe two more algorithms that find  a minimal weight spanning tree in a given connected graph $G$.
For simplicity, we assume that $G$ has distinct weights.

\parbf{Prim's algorithm.}
\begin{enumerate}[(i)]
\item Start with a tree $T$ formed by a single vertex $v$ in $G$.
\item\label{main:prim} Choose an edge $e$ with the minimal weight that connects $T$ to a vertex not yet in the tree.
Add $e$ to $T$.
\item Repeat step \ref{main:prim} until all vertices are in the tree $T$.
\end{enumerate}

\parbf{Borůvka's algorithm.}
\begin{enumerate}[(i)]
\item Start with the forest $F$ formed by all vertices of $G$ and no edges.
\item\label{main:boruvka} For each component $T$ of $F$, find an edge $e$ with minimal weight that connects $T$ to another component of $F$; add $e$ to $F$.
\item Repeat step \ref{main:boruvka} until $F$ has one component. 
\end{enumerate}

\begin{thm}{Exercise}\label{ex:PB}
Show that both (a) Prim's algorithm and (b) Borůvka's algorithm produce a minimum weight spanning tree for a given connected graph with all different weights.
\end{thm}

\begin{thm}{Exercise}\label{ex:KPB}
Compare Kruskal’s algorithm, Prim's algorithm, and Borůvka's algorithm.
Which one should work faster?
When and why?
\end{thm}

\begin{thm}{Exercise}\label{ex:deleting-algorithm}
Assume $G$ is a connected graph with weighted edges and all weights are different.
Let $T$ be a minimum weight spanning tree in $G$.
\begin{enumerate}[(a)]
\item\label{ex:deleting-algorithm:a} Suppose $v$ is an end vertex in $G$ and $e$ is its adjacent edge.
Show that $e$ belongs to $T$.
\item\label{ex:deleting-algorithm:b} Suppose $C$ is a cycle in $T$ and $f$ is an edge in $C$ with maximal weight.
Show that $f$ does not belong to $T$.
\item\label{ex:deleting-algorithm:c}
Come up with a minimum-spanning-tree algorithm based on observations \ref{ex:deleting-algorithm:a} and \ref{ex:deleting-algorithm:b}.
\end{enumerate}
\end{thm}

\section{Remarks}

The NP-hardness of the traveling salesman problem follows from a result by Richard Karp \cite{karp}.
The discussed algorithms are named after Joseph Kruskal, Robert Prim and Otakar Borůvka.
Joseph Kruskal also described the so-called \index{reverse-delete algorithm}\emph{reverse-delete algorithm} \cite{kruskal}, which is closely relevant to Exercise~\ref{ex:deleting-algorithm}.
Prim's algorithm was first discovered by Vojtěch Jarník,
rediscovered  independently two times: by Robert Prim and by Edsger Dijkstra.

For more on optimization problems, read the book by Dieter Jungnickel~\cite{jungnickel}.
A reader-friendly introduction to the P=NP problem is given in the book by 
Thomas Cormen,
Charles Leiserson,
Ronald Rivest, and
Clifford Stein \cite[Chapter 34]{cormen-leiserson-rivest-stein}.

\chapter[Marriage theorem and its relatives]{Marriage theorem\\ and its relatives}
\chaptermark{Marriage theorem}

\section{Alternating and augmenting paths}

Recall that a \index{matching}\emph{matching} in a graph is a set of edges without common vertices. 

Let $M$ be a matching in a graph $G$.
Suppose a path $P$ in $G$ alternates between edges from $M$ and edges not from $M$,
then $P$ is called an \index{alternating path}\emph{$M$-alternating path}.

If an alternating path connects two unmatched vertices of $G$, then it is called \index{augmenting path}\emph{$M$-augmenting}.
An $M$-augmenting path $P$ can be used to improve the matching $M$;
namely by deleting all the edges of $P$ in $M$
and adding the remaining edges of $P$, we obtain a new matching $M'$ with more edges.
This construction implies the following:

\begin{thm}{Observation}\label{obs:augmenting}
Assume $G$ is a graph and $M$ is a maximal matching in $G$.
Then $G$ has no $M$-augmenting paths.
\end{thm} 

\begin{wrapfigure}[5]{r}{28 mm}
\vskip-4mm
\centering
\includegraphics{mppics/pic-60}
\vskip-0mm
\end{wrapfigure}

On the diagrams we denote the edges in $M$ by solid lines and the remaining edges by dashed lines.

\begin{thm}{Exercise}\label{ex:augmenting-path}
Connect two unmatched vertices on the diagram by an augmenting path.
Use it to construct a larger matching.
\end{thm}

Recall that a \index{bigraph}\emph{bigraph} stands for a \index{bipartite graph}\emph{bipartite graph}.

\begin{thm}{Exercise}\label{ex:bigraph-matching}
Let $M$ be a matching in a bigraph $G$.
Suppose $P$ is an $M$-augmenting path.
Show that $P$ has its ends in the opposite parts of $G$.
\end{thm}

The following theorem states the converse to the observation.

\begin{thm}{Theorem}
Assume $M$ is a matching in a graph $G$.
If $M$ is not maximal, then $G$ contains an $M$-augmenting path.
\end{thm}

This theorem implies that the Hungarian algorithm \cite[Section 7.2]{hartsfield-ringel} produces a maximal matching.
Indeed, the Hungarian algorithm checks all $M$-alternating paths in the graph and it will find an $M$-augmenting path if it exists.
Once it is found it could be used to improve the matching $M$.
Therefore, by the theorem, theHungarian algorithm improves any nonmaximal matching;
in other words, it can stop only at a maximal matching.

\parit{Proof.}
Suppose $M$ is a nonmaximal matching; that is, there is a matching $M'$ with a larger number of edges.
Consider the subgraph $H$ of $G$ formed by all edges in $M$ and $M'$.

Note that each component of $H$ is either an $M$-alternating path or a cycle.
Each cycle in $H$ has the same number of edges from $M$ and $M'$.
Since $|M'|>|M|$, the one component of $H$ has to be a path that starts and ends with an edge in $M'$. 
Evidently, such a path is $M$-augmenting.
\qeds

\section{Marriage theorem}

Assume that $G$ is a graph and $S$ is a subset of its vertex set.
We say that a matching $M$ of $G$ {}\emph{covers} $S$ if any vertex in $S$ is incident to an edge in $M$.

Given a set of vertices $W$ in a graph $G$, the set $W'$ of all vertices adjacent to at least one of the vertices in $W$ will be called the \index{set of neighbors}\emph{set of neighbors} of $W$.
Note that if $G$ is a bigraph and $W$ lies in the left part, then $W'$ lies in the right part. 

\begin{thm}{Marriage theorem}\label{thm:marriage}
Let $L$ and $R$ be the left and right parts of a bigraph $G$.
There is a matching that covers $L$ if and only if for any $W\subset L$, we have 
\[|W'|\ge |W|\]
where $W'\subset R$ is the set of all neighbors of $W$.

\end{thm}

\parit{Proof.}
Assume that a matching $M$ is covering $L$.
Note that for any set  $W\subset L$, the set $W'$ of its neighbors includes the vertices matched with~$W$.
In particular,
\[|W'|\ge |W|;\]
it proves the ``only if'' part.

Consider a maximal matching $M$ of $G$.
To prove the ``if'' part, it is sufficient to show that $M$ covers $L$.
Arguing by contradiction, assume that there is a vertex $w$ in $L$ that is not incident to any edge in $M$.

Consider the maximal set $S$ of vertices in $G$ which are reachable from $w$ by 
$M$-alternating paths.
Denote by $W$ and $W'$ the set of left and right vertices in $S$ respectively.

Since $S$ is maximal, $W'$ is the set of neighbors of $W$. 
According to Observation \ref{obs:augmenting}, the matching $M$ provides a bijection between $W-w$ and $W'$.
In particular, 
\[|W|=|W'|+1;\] 
the latter contradicts the assumption.
\qeds

\begin{thm}{Exercise}\label{ex:1-factor}
Assume $G$ is an $r$-regular bigraph; $r\ge 1$.
Show that 
\begin{enumerate}[(a)]
\item\label{ex:1-factor:a} $G$ admits a 1-factor;
\item\label{ex:1-factor:b} the edge chromatic number of $G$ is $r$; in other words, $G$ can be decomposed into $1$-factors.
\end{enumerate}

\end{thm}

{

\begin{wrapfigure}{r}{25 mm}
\vskip-6mm
\centering
\includegraphics{mppics/pic-61}
\vskip-0mm
\end{wrapfigure}

In Exercise~\ref{ex:1-factor} one has to assume that $G$ is bipartite.
Indeed, the following exercise states that the graph shown on the diagram has no 1-factor,
while it is 3-regular.

\begin{thm}{Exercise}\label{ex:no-1-factor}
Prove that the graph shown on the diagram has no 1-factor.
\end{thm}

}

\begin{thm}{Exercise}\label{ex:kids}
Children from 25 countries, 10 kids from each, decided to stand in a rectangular formation with 25 rows of 10 children in each row.
Show that you can always choose one child from each row so that all 25 of them will be from different countries.
\end{thm}

\begin{thm}{Exercise}\label{ex:sons(king)}
The sons of the king divided the kingdom between each other into 23 parts of equal area --- one for each son.
Later a new son was born. 
The king proposed a new subdivision into 24 equal parts and gave one of the parts to the newborn son.

Show that each of the 23 older sons can choose a part of land in the new subdivision which overlaps with his old part.
\end{thm}

\begin{thm}{Exercise}\label{ex:nxn-table}
A table $n{\times}n$ is filled with nonnegative numbers.
Assume that the sum in each column and each row is 1.
Show that one can choose $n$ cells with positive numbers that do not share columns and rows. 
\end{thm}

\begin{thm}{Advanced exercise}\label{ex:camel17}
In a group of people, for some fixed $s$ and any $k$,
any $k$ girls like at least $k-s$ boys in total.
Show that then all but $s$ girls may get married to the boys they like.
\end{thm}

\section{Vertex covers}

A set $S$ of vertices in a graph is called a \index{vertex cover}\emph{vertex cover} if any edge is incident to at least one of the vertices in $S$.

\begin{thm}{Theorem}\label{thm:vertex-cover}
In any bigraph, the number of edges in a maximal matching equals the number of vertices in a minimal vertex cover.
\end{thm}

It is instructive to do the following exercise before reading the proof.

\begin{thm}{Exercise}\label{ex:two-paths}
Let $M$ be a maximal matching in a bigraph $G$.
Assume two unmatched vertices $l$ and $r$ lie on the opposite parts of $G$.
Show that no pair of $M$-alternating paths starting from $l$ and $r$ can have a common vertex.
\end{thm}

On the following diagram, a maximal matching is marked by solid lines;
the remaining edges of the graph are marked by dashed lines.
The vertices of the cover are marked in black and the remaining vertices in white;
the unmatched vertices are marked by a cross.

\begin{wrapfigure}{o}{30 mm}
\vskip-5mm
\centering
\includegraphics{mppics/pic-62}
\vskip1mm
\end{wrapfigure}

\parit{Proof.}
Fix a bigraph $G$;
denote by $L$ and $R$ its left and right parts.
Let $M$ be a matching, and $S$ be a vertex cover in $G$.

By the definition of vertex cover, any edge $m$ in $M$ is incident to at least one vertex in $S$.
Therefore, 
\[|S|\ge |M|.\] 
That is, the number of vertices in any vertex cover $S$ is at least as large as the number of edges in any matching $M$.

Now assume that $M$ is a maximal matching.
Let us construct a vertex cover $S$ such that $|S|\z=|M|$.

Denote by $U_L$ and $U_R$ the set of left and right unmatched vertices (these are marked by a cross on the diagram).
Denote by $Q_L$ and $Q_R$ the set of vertices in $G$ which can be reached by $M$-alternating paths starting from $U_L$ and from $U_R$, respectively.

Note that $Q_L$ and $Q_R$ do not overlap.
Otherwise, there would be an $M$-augmenting path from $U_L$ to $U_R$ (compare to Exercise~\ref{ex:two-paths}).
Therefore, $M$ is not maximal --- a contradiction.

Further note that if $m$ is an edge in $M$, then both of its end vertices lie either in $Q_L$ or $Q_R$, or neither.

Let us construct the set $S$ by taking one incident vertex (left or right) of each edge $m$ in $M$ by the following rule:
\emph{if $m$ connects vertices in $Q_L$, then include its right vertex in $S$;
otherwise include its left vertex}.
Since $S$ has exactly one vertex incident to each edge of $M$, we have
\[|S|=|M|.\]

It remains to prove that $S$ is a vertex cover;
that is, at least one vertex of any edge $e$ in $G$ is in $S$.

Note that if the left vertex of $e$ lies in $Q_L$, then $e$ is an edge on an $M$-alternating path starting from $U_L$. 
Therefore, the right vertex of $e$ also lies in $Q_L$.

Therefore, it is sufficient to consider only the following three cases:
\begin{itemize}
\item The edge $e$ has
its right vertex in $Q_L$ and its left vertex outside of $Q_L$.
In this case, both vertices of $e$ lie in~$S$.
\item The edge $e$ connects vertices in $Q_L$.
In this case, the right vertex of $e$ is in $S$.
\item The edge $e$ connects vertices outside of $Q_L$. 
In this case, the left vertex of $e$ is in~$S$.
\qeds
\end{itemize}

\begin{thm}{Exercise}\label{ex:rooks}
A few squares on a chessboard are marked.
Show that the minimal number of rows and columns that cover all marked squares is the same as the maximal number of rooks on the marked squares that do not threaten each other.
\end{thm}

\section{Edge cover}

A collection of edges $N$ in a graph is called an \index{edge cover}\emph{edge cover} if every vertex is incident with at least one of the edges in $N$.

\begin{wrapfigure}{r}{25 mm}
\vskip-2mm
\centering
\includegraphics{mppics/pic-63}
\vskip0mm
\end{wrapfigure}

On the diagram, two edge covers of the same graph are marked in solid lines.

\begin{thm}{Theorem}
Let $G$ be a connected graph with $p$ vertices and $p>1$.
Assume that a minimal edge cover $N$ of $G$ contains $n$ edges, and a maximal matching $M$ of $G$ contains $m$ edges.
Then
\[m+n=p.\]
\end{thm}

The following exercise will guide you thru the proof.

\begin{thm}{Exercise}\label{ex:monotree}
Suppose $G$, $N$, $M$, $p$, $n$, and $m$ are as in the theorem.
Denote by $m'$ the number of components in $N$.
\begin{enumerate}[(a)]
\item\label{ex:monotree:m'+n=p} Show that $N$ contains no paths of length 3 and no triangle.
Conclude that each component of $N$ is a \index{star}\emph{star};
that is, it is isomorphic to $K_{1,k}$ for some $k$.
Use it to show that 
\[m'+n=p.\]
\item\label{ex:monotree:m+n>=p} Choose $m'$ edges, one in each component of $N$.
Observe that it is a matching.
Conclude that $m'\le m$.
Use \ref{ex:monotree:m'+n=p} to show that 
\[m+n\ge p.\]
\item\label{ex:monotree:m+n=<p} Consider a subgraph $H$ of $G$ formed by all edges in the matching $M$ and one edge incident to each unmatched vertex in $G$.
Observe that $H$ is a vertex cover with $p-m$ edges. 
Conclude that 
\[m+n\le p.\]
\item Observe that \ref{ex:monotree:m+n>=p} and \ref{ex:monotree:m+n=<p} imply the theorem.
\end{enumerate}

\end{thm}

\section{Minimal cut}

Recall that a \index{directed graph}\emph{directed graph} (or briefly a \index{digraph}\emph{digraph})
is a graph where the edges have a direction (also known as orientation) associated with them;
that is, an edge in a digraph is defined as an \textit{ordered} pair of vertices.
A choice of directions on every edge of a graph is called its \index{orientation}\emph{orientation}.

\begin{thm}{Min-cut theorem}\label{thm:mincut}
Let $G$ be a digraph. 
Fix vertices $s$ and $t$ in~$G$.
Then the maximal number of directed paths from $s$ to $t$ which do not have common edges equals the minimal number of edges one can remove from $G$ so that there will be no directed path from $s$ to $t$.
\end{thm}

\parit{Proof.}
Denote by $m$ the maximal number of directed paths from $s$ to $t$ which do not have common edges.
Denote by $n$ the minimal number of edges one can remove from $G$ to disconnect $t$ from $s$; more precisely, after removing $n$ edges from $G$, there will be no directed path from $s$ to~$t$.

Let $P_1, \dots, P_m$ be a maximal collection of directed paths from $s$ to $t$ which have no common edges.
Note that in order to disconnect $t$ from $s$, we have to cut at least one edge in each path $P_1, \dots, P_m$.
In particular, $n\ge m$.

Consider the new orientation on $G$ where each path $P_i$ is oriented backward --- from $t$ to $s$.

Consider the set $S$ of the vertices which are reachable from $s$ by directed paths for this new orientation.

Assume $S$ contains $t$; that is, there is a path $Q$ from $s$ to $t$, which
can move along $P_i$ only backward.

(Further, the path $Q$ will be used the same way as the augmenting path in the proof of the marriage theorem.
In a sequence of moves, we will improve the collection $Q,P_1,\dots,P_m$ so that there will be no overlaps.
On the diagram, a case with $m=1$ that requires two moves is shown;
$P_1$ is marked by a solid line and $Q$ is marked by a dashed line.)

Since $P_1,\dots,P_m$ is a maximal collection, $Q$ overlaps with some of the paths $P_1,\dots,P_m$.
Without loss of generality, we can assume that $Q$ first overlaps with $P_1$ --- assume it meets $P_1$ at the vertex $v$ and leaves it at the vertex $w$.
Let us modify the paths $Q$ and $P_1$ the following way:
Instead of the path $P_1$ consider the path $P_1'$ that goes along $Q$ from $s$ to $v$ and after that goes along $P_1$ to $t$.
Instead of the path $Q$, consider the trail $Q'$ which goes along $P_1$ from $s$ to $w$ and after that goes along $Q$ to $t$.

\begin{wrapfigure}{o}{30 mm}
\vskip-0mm
\centering
\includegraphics{mppics/pic-64}
\vskip0mm
\end{wrapfigure}
If the constructed trail $Q'$ is not a path (that is, if $Q'$ visits some vertices several times), then we can discard some circuits from $Q'$ to obtain a genuine path, 
which we will still denote by $Q'$.

Note that the obtained collection of paths $Q', P_1',P_2\dots,P_m$ satisfies the same conditions as the original collection.
Further, since we discard the part of $P_1$ from $w$ to $v$, the total number of edges in $Q', P_1',P_2\dots,P_m$
is smaller than in the original collection $Q, P_1,P_2\z\dots,P_m$.
Therefore, by repeating the described procedure several times, we get $m+1$ paths without overlaps --- a contradiction.

{

\begin{wrapfigure}{r}{40 mm}
\vskip-2mm
\centering
\includegraphics{mppics/pic-65}
\vskip0mm
\end{wrapfigure}

It follows that $S\not\ni t$.

Note that all edges which connect $S$ to the remaining vertices of $G$ are directed toward $S$.
That is, every such edge which comes out of $S$ in the original orientation belongs to one of the paths $P_1,\dots,P_m$.

Moreover, for each path $P_i$ there is only one such edge.
In other words, if a path $P_i$ leaves $S$, then it cannot come back. 
Otherwise, $S$ could be made larger by moving backward along $P_i$.
Therefore, cutting one such edge in each path $P_1,\dots,P_m$ makes it impossible to leave $S$.
In particular, we can disconnect $t$ from $s$ by cutting $m$ edges from $G$; that is, $n\le m$.
\qeds

}

\parbf{Remark.}
The described process has the following physical interpretation.
Think of each path $P_1,\dots,P_m$, and $Q$ like water pipelines from $s$ to $t$.
\begin{figure}[ht!]%{r}{30 mm}
\vskip0mm
\centering
\includegraphics{mppics/pic-66}
\vskip0mm
\end{figure}
At each overlap of $Q$ with another path $P_i$, the water in $P_i$ and $Q$ runs in the opposite directions.
So we can cut the overlapping edges and reconnect the open ends of the pipes to each other while keeping the water flow from $s$ to $t$ unchanged.
As a result, we get $m+1$ pipes from $s$ to $t$ with no common edges and possibly some cycles which we can discard.
An example of this procedure for two paths $P_1$ and $Q$ is shown on the diagram;
as above, $P_1$ is marked by a solid line and $Q$ is marked by a dashed line.

{

\begin{wrapfigure}{r}{35 mm}
\vskip-0mm
\centering
\includegraphics{mppics/pic-67}
\end{wrapfigure}

\begin{thm}{Advanced exercise}\label{ex:min-cut-marriage}
Assume $G$ is a bigraph.
Let us add two vertices, $s$ and $t$, to $G$ so that $s$ is connected to each vertex in the left part of $G$, and $t$ is connected to each vertex in the right part of $G$.
Orient the graph from left to right.
Denote the obtained digraph by~$\hat G$.

Give another proof of the marriage theorem for a bigraph $G$, applying the min-cut theorem to the digraph $\hat G$. 
\end{thm}

}

\section{Remarks}

The marriage theorem was proved by Philip Hall in \cite{hall};
it has many applications in all branches of mathematics.
The theorem on vertex cover was discovered by D\'enes K\H{o}nig \cite{konig}, and rediscovered by Jen\H{o} Egerv\'ary \cite{egervary}.
The theorem on min-cut was proved by Peter Elias, Amiel Feinstein, and Claude Shannon \cite{elias-feinstein-shannon};
it was rediscovered by Lester Ford and Delbert Fulkerson \cite{ford-fulkerson}.

These theorems form a base for the Hungarian algorithm that solves the \index{assignment problem}\emph{assignment problem}; see page \pageref{assignment problem}.
This algorithm was found by Harold Kuhn, but, as it was discovered by François Ollivier, essentially the same algorithm was found much earlier by Carl Gustav Jacobi~\cite{ollivier}.

An extensive overview of the marriage theorem and its relatives is given by Alexandr Evnin in \cite{evnin}.

\chapter{Toroidal graphs}

Recall that a graph is called \index{planar graph}\emph{planar} if it can be drawn on a plane with no crossings;
the latter is equivalent to the existence of its drawing on a sphere.
Here and further, we say \index{drawing}\emph{drawing} for \textit{drawing with no crossings}.  

\begin{wrapfigure}{o}{40 mm}
\vskip-2mm
\centering
\includegraphics{mppics/pic-101}
\end{wrapfigure}

In this chapter, we will discuss drawing of graphs on a torus --- another surface that can be obtained by revolving a circle about an axis (this is the surface of a donut).

A graph is called \index{toroidal graph}\emph{toroidal}, if it can be drawn on a torus with no crossing.

If one cuts a torus along a parallel and a meridian as shown on the diagram,
then the obtained surface can be developed into a square.%
\footnote{Formally speaking, this means that there is a continuous map $f\:\square\to T$ from a square to torus  such that $f(x)=f(y)$ if and only if $x=y$ or $x$ and $y$ are corresponding points on the opposite sides of the square.}
It gives a convenient way to describe drawings of graphs on the torus which will be called \index{square diagram}\emph{square diagram}.
One only has to remember that the corresponding points on opposite sides of the square are identified in the torus;
in particular, all four vertices of the square correspond to one point in the torus.

{

\begin{wrapfigure}{r}{30 mm}
\vskip-4mm
\centering
\includegraphics{mppics/pic-102}
\end{wrapfigure}

For example, on the given square diagram you see a drawing of the  complete graph $K_5$;
in particular it shows that $K_5$ is toroidal.
Note that the edge $xy$,
after coming to the right side of the square, reappears at the corresponding point of the left side and goes further to~$y$.
Similarly, the edge $vw$ comes to the top and reappears at the bottom.
At these points, the edges cross the parallel and the meridian.

}

The following square diagrams show that $K_{4,4}$ and $K_7$ are toroidal graphs as well.

\begin{figure}[ht!]
\begin{minipage}{.48\textwidth}
\centering
\includegraphics{mppics/pic-104}
\end{minipage}\hfill
\begin{minipage}{.48\textwidth}
\centering
\includegraphics{mppics/pic-106}\label{K5-toroidal}
\end{minipage}

\medskip

\begin{minipage}{.45\textwidth}
\centering
\caption*{$K_{4,4}$}
\end{minipage}\hfill
\begin{minipage}{.45\textwidth}
\centering
\caption*{$K_7$}
\end{minipage}
\vskip-4mm
\end{figure}

\begin{thm}{Exercise}\label{ex:crossing1}
Show that any graph with crossing number 1 is toroidal.
That is, if a graph admits a drawing on a sphere with one crossing, then it admits a drawing on a torus with no crossings.
\end{thm}

\begin{thm}{Exercise}\label{ex:toroidal-graphs}
Show that each of the following graphs is toroidal; construct a corresponding square diagram in each case.

(a) \includegraphics{mppics/pic-107}
(b) \includegraphics{mppics/pic-108}
(c) \includegraphics{mppics/pic-117}
\end{thm}

\section{Simple regions}

Choose a drawing $G$ of a pseudograph on a torus or sphere;
it subdivides the surface into regions.
A region $R$ is called \index{simple region}\emph{simple} if its interior can be parameterized by an open plane disc.%
\footnote{Formally speaking, this means that there is a continuous bijection from the interior of $R$ to an open disc in the plane such that its inverse is also continuous.}

\begin{wrapfigure}{r}{30 mm}
\vskip-6mm
\centering
\includegraphics{mppics/pic-109}
\end{wrapfigure}

The annulus shown on the diagram is an example of a nonsimple region.
On a sphere it may appear only in drawings of nonconnected graphs.
That is, if a drawing of a pseudograph $G$ on a sphere has a nonsimple region, then $G$ is not connected.
This statement should be intuitively obvious, but the proof is not trivial.  

Drawings of connected graphs on a torus may have nonsimple regions.
For example, in the first drawing of $K_4$ below,
\begin{figure}[ht!]
\begin{minipage}{.45\textwidth}
\centering
\includegraphics{mppics/pic-110}
\end{minipage}
\hfill
\begin{minipage}{.45\textwidth}
\centering
\includegraphics{mppics/pic-111}
\end{minipage}

\medskip

\begin{minipage}{.45\textwidth}
\centering
\caption*{$p=4$, $q=6$, $r=3$; the regions $B$, $C$ are simple, and $A$ is not.}
\end{minipage}\hfill
\begin{minipage}{.45\textwidth}
\centering
\caption*{$p=4$, $q=6$, $r=2$; both regions $D$ and $E$ are simple.}
\end{minipage}
\vskip-4mm
\end{figure} 
the regions $B$ and $C$ are simple, and $A$ is not (it contains the meridian of the torus).
The second drawing of $K_4$ has only two regions, $D$ and $E$, and both of them are simple.

Further, we will use the following claim, which should be intuitively obvious.
We do not present its proof, but it is not hard;
a reader familiar with topology may consider it as an exercise.

{

\begin{wrapfigure}{r}{15 mm}
\vskip-4mm
\centering
\includegraphics{mppics/pic-121}
\medskip
\includegraphics{mppics/pic-120}
\end{wrapfigure}

\begin{thm}{Claim}\label{clm:cut}
Let $G$ be a drawing of a pseudograph on a torus or sphere and $R$ its simple region.
Suppose that a drawing $G'$ obtained from $G$ by adding a new edge $e$ in the region $R$.
\begin{enumerate}[(a)]
\item If both ends of $e$ are in $G$, then $e$ divides $R$ into two simple regions.
\item If only one the ends of $e$ is in $G$, then $e$ does not divide $R$ and the corresponding region of $G'$ is simple.
\end{enumerate}
\end{thm}

}

\begin{thm}{Advanced exercise}\label{ex:nonplanar-toroidal}
Suppose that  $G$ is a drawing of a connected nonlanar graph on the torus.
Show that each region of $G$ is simple.
\end{thm}

\section{Euler's formula}

\begin{thm}{Theorem}\label{thm:euler>=}
For any drawing $G$ of a connected pseudograph on the torus we have
\[p-q+r\ge 0,\]
where $p$, $q$, and $r$ denote the number of vertices, edges and regions in~$G$, respectively.

Moreover, equality holds if all regions of $G$ are simple.
\end{thm}

Note that the inequality might be strict.
For example, for the first drawing of $K_4$, above we have
\[p-q+r=4-6+3=1>0.\]
It happens since the region $A$ is not simple.
For the second drawing of $K_4$ we have the equality
\[p-q+r=4-6+2=0,\]
as it is supposed to be by the second part of the theorem.

In the proof we will use the following lemma which is a simple corollary of Claim~\ref{clm:cut}.
Given a drawing of pseudograph $G$ with $p$ vertices, $q$ edges and $r$ regions,
set 
\[\Sigma_G=p-q+r.\]

\begin{thm}{Lemma}\label{lem:euler}
Let $G$ be a drawing of a pseudograph on a torus or sphere.
Suppose that another drawing $G'$ of a connected pseudograph that is obtained from $G$ by adding a new edge $e$.
Then 
\[\Sigma_{G'}\le\Sigma_G.\eqlbl{G'=<G}\]
Moreover, if all regions of $G$ are simple, then 
\begin{enumerate}[(a)]
\item we have equality in \ref{G'=<G} and  
\item\label{lem:euler:simple} all regions  of $G'$ are simple as well.
\end{enumerate}

\end{thm}

\parit{Proof.}
Since $G'$ is connected, one of the ends of $e$ belongs to $G$.
Denote by $R$ the region of $G$ that contains $e$.

If the other end of $e$ is not in $G$, then $e$ does not divide its region.
In this case $G'$ has an extra vertex and an extra edge, and the number of regions did not change; that is,
$p'=p+1$, $q'=q+1$, and $r'=r$.
Therefore, 
\[\Sigma_{G'}=p'-q'+r'=p-q+r=\Sigma_G.\]

If the other end of $e$ is in $G$, then $e$ may divide $R$ into two regions or may not divide it.
According to Claim~\ref{clm:cut}, the latter may happen only if $R$ is not simple.
Note that after adding $e$, the number of vertices did not change.
That is,
$p'=p$, $q'=q+1$, $r'\le r+1$ and the equality holds if $R$ is simple.
Therefore, 
\[\Sigma_{G'}=p'-q'+r'\le p-q+r=\Sigma_G\]
and the equality holds if $R$ is simple.

According to Claim~\ref{clm:cut}, if $R$ is simple, then the region(s) of $G'$ that correspond to $R$ are simple as well.
The rest of the regions did not change.
Hence \ref{lem:euler:simple} follows.
\qeds

\parit{Proof of \ref{thm:euler>=}.}
We need to show that 
\[\Sigma_G\le 0\eqlbl{G=<0}\] 
for any drawing of connected pseudograph $G$ on the torus $T$.

\begin{wrapfigure}{r}{40 mm}
\vskip-0mm
\centering
\includegraphics{mppics/pic-115}
\vskip2mm
\end{wrapfigure}

Let $H$ be a drawing of the pseudograph formed by one meridian and one parallel as shown.
Note that it has 1 vertex, 2 edges, and 1 region;
therefore 
\[\Sigma_H=1-2+1=0.\eqlbl{H=0}\]

Without loss of generality, we may assume that the edges of $G$ and $H$ intersect and have only a finite number of points of intersection.
The latter can be achieved by perturbing the drawing of $G$.

Let us subdivide the graphs $G$ and $H$ by adding a new vertex at every crossing point of the graphs.
The obtained graphs, say $\bar G$ and $\bar H$, are subgraphs of a bigger graph, say $W$, formed by all edges and vertices of $\bar G$ and $\bar H$.
Note that adding a vertex on an edge increases the number of vertices and edges by $1$ and the number of regions stays the same. 
Since the subdivision is obtained by adding a finite number of extra vertices, we get that
\[\Sigma_G=\Sigma_{\bar G}\quad\text{and}\quad\Sigma_H=\Sigma_{\bar H}.
\eqlbl{G=G,H=H}\]

By construction, $W$ is connected.
If $W\ne G$, then there is an edge $e$ of $W$ that is not in $G$ but has one of its ends in $G$.
By adding $e$ to $G$ and applying the procedure recursively, we will obtain $W$ in a finite number of steps.
Applying Lemma~\ref{lem:euler} at each step, we get
\[\Sigma_W\le \Sigma_{\bar G}.
\eqlbl{W<G'}\]

Analogously, $W$ can be obtained from $H$ in a finite number of steps by adding one edge at a time.
Since the only region of $H$ is simple, the obtained drawings will have only simple regions.
Therefore, applying Lemma~\ref{lem:euler} at each step, we get
\[\Sigma_W= \Sigma_{\bar H}.
\eqlbl{W=G'}\]

Finally, \ref{H=0}, \ref{G=G,H=H}, \ref{W<G'}, and \ref{W=G'} imply \ref{G=<0};
indeed
\[\Sigma_G=\Sigma_{\bar G}\ge \Sigma_W=\Sigma_{\bar H}=\Sigma_H=0.\]
\qedsf

\begin{thm}{Exercise}\label{ex:toroidal-girth}
Suppose that $G$ is a toroidal graph with girth~$\ge4$.
Show that 
\[q\le 2\cdot p,\]
where $p$ and $q$ denote number of vertices and edges in $G$.
\end{thm}

\begin{thm}{Exercise}\label{ex:K5-torus}
Is there a drawing of $K_5$ on the torus with only square regions?
If ``yes'', then draw a square diagram; if ``no'', explain why.
\end{thm}

\section{The seven-color theorem}

\begin{thm}{Theorem}\label{thm:7-colors}
The chromatic number of any toroidal graph cannot exceed~$7$.
\end{thm}

Recall that $K_7$ is toroidal; see diagram on page \pageref{K5-toroidal}.
Therefore, the theorem gives an optimal bound.
In the proof, we will use the following lemma, which is a simple corollary of Euler's formula (\ref{thm:euler>=}).

\begin{thm}{Lemma}\label{cor:q=<3p}
Let $G$ be a toroidal graph with $p$ vertices and $q$ edges.
Then 
\[q\le 3\cdot p.\]

\end{thm}

\parit{Proof.}
Choose a drawing of $G$ on a torus.
By Euler's inequality we have
\[p-q+r\ge 0,\]
where $r$ denotes the number of regions in the drawing.

Note that each region in the drawing has at least $3$ sides
and each edge of $G$ appears twice as a side of a region.
Therefore, 
\[3\cdot r\le 2\cdot q.\]
These two inequalities imply the corollary.
\qeds

\parit{Proof of \ref{thm:7-colors}.}
Suppose that there is a toroidal graph $G$ that requires $8$ colors.
Choose a critical subgraph $H$ in $G$ with chromatic number 8.
Denote by $p$ and $q$ the number of vertices and edges in $H$.

By \cite[Theorem 2.1.3]{hartsfield-ringel}, each vertex in $H$ has a degree of at least $7$.
By the handshake lemma \cite[Theorem 1.1.1]{hartsfield-ringel}, we have 
\[7\cdot p\le 2\cdot q.\]

On the other hand, by Lemma \ref{cor:q=<3p}, we have
\[q\le 3\cdot p.\]
These two inequalities contradict each other.
 \qeds

\section{A remark about forbidden minors}

It is straightforward to see that if $G$ is toroidal, then any graph obtained from $G$ by the deletion of a vertex or an edge or by contracting an edge is also toroidal.

A graph that can be obtained from a given graph $G$ by applying a sequence of such operations is called a \index{minor}\emph{minor} of $G$.
The pseudograph $G$ is considered to be a minor of itself.
The other minors require at least one deletion or contraction; they are called \index{proper minor}\emph{proper minors}.

Note that the statement above implies that any minor of a toroidal graph is toroidal, or in other words, \textit{toridality is inherited by minors}.

The following deep result was proved by Neil Robertson and Paul Seymour \cite{robertson-seymour}.

\begin{thm}{Theorem}
Any property of pseudographs that is inherited by minors can be 
described by a finite set of \index{forbidden minior}\emph{forbidden minors};
that is,  a pseudograph meets the property if and only if none of its minors are forbidden.
\end{thm}

For example, note that deletion and contraction do not create cycles;
that is, any minor of a forest is a forest.
Forests can be described by one forbidden minor --- a pseudograph formed by one loop.
Indeed, if a pseudograph has a cycle, then by a sequence of deletions and contractions, one can get a single loop from it.

\begin{thm}{Exercise}\label{ex:forbidden-minors}
Describe the following classes of graphs by a single forbidden minor.
\begin{enumerate}[(a)]
 \item\label{ex:forbidden-minors:tree5} Graphs that do not contain trees with 5 end vertices.
 \item\label{ex:forbidden-minors:cycle} Graphs in which any two cycles have at most one common vertex.
\end{enumerate}
\end{thm}

A more complicated example is given by the Pontryagin--Kuratowski theorem;
it states that \textit{planar graphs are characterized by two forbidden minors: $K_5$ and $K_{3,3}$.}

Since toroidality is inherited by its minors, it can be described by a set of forbidden minors.
The complete list of forbidden minors for this problem has to be huge and is not yet known.

{

\begin{wrapfigure}{r}{30 mm}
\vskip-7mm
\centering
\includegraphics{mppics/pic-116}
\vskip0mm
\end{wrapfigure}

\begin{thm}{Advanced exercise}\label{ex:nontoroidal}
Show that the graph on the diagram is not toroidal,
but every proper minor of this graph is toroidal.
\end{thm}

}

\section{Other surfaces}

One may draw graphs on the so-called the \index{surface of genus $g$}\emph{surfaces of genus $g$}.
These surfaces can be obtained by attaching $g$ toruses to each other; 
\begin{figure}[ht!]%{r}{30 mm}
\vskip-0mm
\centering
\includegraphics{mppics/pic-122}
\vskip-0mm
\end{figure}
a surface of genus $g=4$ is shown in the picture.
Namely, to construct a surface $S'$ of genus $g+1$, start with a surface $S$ of genus $g$ and a torus $T$, drill a hole in each and reconnect them to each other as shown.%
\footnote{The described construction is called \index{connected sum}\emph{connected sum}; so we can say that \textit{connected sum of a surface of genus $g$ and a torus is a surface of genus $g+1$}.}
\begin{figure}[ht!]%{r}{30 mm}
\vskip-0mm
\centering
\includegraphics{mppics/pic-123}
\vskip-0mm
\end{figure}

It is natural to assume that the sphere has genus $0$ and the torus has genus $1$.
(In general, genus tells how many disjoint closed curves one could draw on the surface so that they do not cut the surface into pieces.)

The simple regions of a drawing on a surface of genus $g$ can be defined the same way as on the torus.
Euler's formula given in Theorem~\ref{thm:euler>=} admits the following straightforward generalization

\begin{thm}{Theorem}\label{thm:euler>=genus}
For any drawing $G$ of a pseudograph on the surface of genus $g$ we have
\[p-q+r\ge 1-2\cdot g,\]
where $p$, $q$, and $r$ denote the number of vertices, edges and regions of~$G$.

Moreover, equality holds if all regions of $G$ are simple.
\end{thm}

The seven-color theorem also admits the following straightforward generalization;
it was proved by Percy John Heawood \cite{heawood}.

\begin{thm}{Theorem}
If a graph $G$ admits a drawing on a surface of genus $g\ge 1$, 
then its chromatic number cannot exceed 
\[\frac{7+\sqrt{1+48\cdot g}}2.\]
\end{thm}

Note that for a sphere (that is, for $g=0$) the formula gives $4$, which is the right bound for the chromatic number for planar graphs. 
However, this is just a coincidence; the proof of Heawood works only for $g\ge 1$.

The estimate in the last theorem is sharp.
The latter was proved by constructing a drawing on the surface of genus $g$
of the complete graphs $K_n$ for any $n\le\frac{7+\sqrt{1+48\cdot g}}2$.
The final step in this construction was made by Gerhhard Ringel and Ted Youngs \cite{ringel-youngs}.
The solution uses the so-called \index{rotation}\emph{rotations of graphs}; which is a combinatoric way to encode a drawing of a graph on a surface.
This is the subject of \cite[Chapter 10]{hartsfield-ringel}.

\begin{wrapfigure}{r}{40 mm}
\vskip-0mm
\centering
\includegraphics{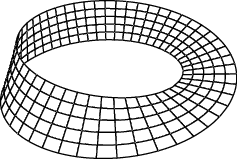}
\vskip-0mm
\end{wrapfigure}

One may consider nonoriented surfaces; for example, the so-called Möbius strip shown on the diagram.
It turns out that it is possible to generalize Euler's inequality and get an exact upper bound for the chromatic numbers for graphs that can be drawn on such surfaces.
(The problem is slightly harder for the so-called Klein bottle, but still, its difficulty is not comparable with the four-color theorem.)

\begin{thm}{Exercise}\label{ex:K6moebius}
Draw the complete graph $K_6$ on a Möbius strip (assume it is made from a transparent material).
\end{thm}

\chapter{Rado graph}

In this chapter, we consider one graph with many surprising properties.
Unlike most of the graphs we have considered so far, this graph has an infinite set of vertices.

\section{Definition}

Recall that a set is \index{countable set}\emph{countable} if it can be enumerated by natural numbers $1,2,\dots$;
it might be infinite or finite.

A \index{countable graph}\emph{countable graph} is a graph with a countable set of vertices;
the set of vertices can be infinite or finite, but it cannot be empty since we always assume that a graph has a nonempty set of vertices.

\begin{thm}{Definition}
A \index{Rado graph}\emph{Rado graph} is a countable graph satisfying the following property:

Given two finite disjoint sets of vertices $V$ and $W$, there exists a
vertex $v\notin V\cup W$ that is adjacent to any vertex in $V$ and nonadjacent to any vertex in $W$.
\end{thm}

The property in the definition will be called the \index{Rado property}\emph{Rado property}; so we can say that \textit{for the sets of vertices $V$ and $W$ in a graph, the Rado property holds} or \textit{does not hold}.

\begin{thm}{Exercise}\label{ex:rado-infty}
Show that any Rado graph has an infinite number of vertices. 
\end{thm}

\begin{thm}{Exercise}\label{ex:rado-diam}
Show that any Rado graph has diameter 2.
\end{thm}

\section{Stability}

The following exercises show that the Rado property is very stable --- 
small changes cannot destroy it.

\begin{thm}{Exercise}\label{ex:rado-partition}
Assume the set of vertices of a Rado graph is partitioned into two subsets.
Show that the subgraph induced by one of these subsets is Rado.
\end{thm}

\begin{thm}{Exercise}\label{ex:R-e-v-rev}
Let $R$ be a countable graph. 
\begin{enumerate}[(a)]
 \item Assume  $e$ is an edge in $R$. 
 Show that $R-e$ is a Rado graph if and only if so is $R$.
 
 \item Assume $v$ is a vertex in $R$. 
 Show that $R-v$ is a Rado graph if so is $R$.
 
 \item Assume $v$ is a vertex in $R$. 
 Consider the graph $R'$ obtained from $R$ by replacing each edge from $v$ by a non-edge, and each non-edge from $v$ by an edge (leaving the rest unchanged).
 Show that $R'$ is a Rado graph if and only if so is $R$. 
\end{enumerate}

\end{thm}

\begin{thm}{Exercise}\label{ex:rado-link}
Let $R$ be a Rado graph.
Assume that $Z$ is the set of all vertices in $R$ adjacent to a given vertex $z$.
Show that the subgraph induced by $Z$ is Rado.
\end{thm}

\section{Existence}

\begin{thm}{Theorem}
There exists a Rado graph.
\end{thm}

\parit{Proof.}
Let $G$ be a finite graph.
Denote by $G'$ the graph obtained from $G$ according to the following rule:
for each subset $V$ of vertices in $G$ add a vertex $v$ and connect it to all the vertices in $V$.

Note that if $G$ has $p$ vertices, then $G'$ has $p+2^p$ vertices --- it has $p$ vertices of $G$ and $2^p$ additional vertices --- one for each of $2^p$ subsets of the $p$-element set (including the empty set).

The original graph $G$ is an induced subgraph in $G'$.
Note also that $G'$ is finite --- it has $p+2^p$ vertices.

By construction, the Rado property holds in $G'$ for any two sets $V$ and $W$ of vertices in $G$ --- the required vertex $v$ is the vertex in $G'$ that corresponds to the subset $V$.

Let $G_1$ be a graph with one vertex.
\begin{figure}[ht!]%{r}{20 mm}
\vskip-0mm
\centering
\includegraphics{mppics/pic-71}
\vskip-0mm
\end{figure}
By repeating the construction, we get a sequence of graphs 
$G_1, G_2,G_3,\dots$,
such that $G_{n+1}=G_n'$ for any~$n$.
The graphs $G_1, G_2,G_3$ are shown on the diagram.%
\footnote{It would be hard to draw $G_4$ since it contains $1+2^1+2^{3}+2^{11}=2059$ vertices, and 
it is impossible to draw $G_5$ --- it has $2059+2^{2059}$ vertices which exceeds by many orders the number of particles in the observable universe.}

Since $G_n$ is a subgraph of $G_{n+1}$ for any $n$, we can consider the union of the graphs in the sequence $(G_n)$; denote it by $R$.
By construction, each graph $G_n$ is a subgraph of $R$ induced by finitely many vertices.
Moreover,  any vertex or edge of $R$ belongs to any $G_n$ with a sufficiently large~$n$.

It remains to show that $R$ is Rado.
Indeed, any two finite sets of vertices $V$ and $W$ belong to $G_n$ for some $n$.
From above, the Rado property holds for $V$ and $W$ in $G_{n+1}$, and therefore in $R$.
\qeds

\parbf{Another construction.}
One could also construct a Rado graph by directly specifying which vertices are adjacent.
Namely, consider the graph $R$ as on the diagram with vertices $r_0,r_1,\dots$
\begin{figure}[ht!]%{r}{45 mm}
\centering
\begin{lpic}[t(-0 mm),b(0 mm),r(0 mm),l(0 mm)]{pics/Rado_graph(.5)}
\lbl[]{5.5,39;$r_0$}
\lbl[]{24.5,39;$r_1$}
\lbl[]{43.5,39;$r_2$}
\lbl[]{62.5,39;$r_3$}
\lbl[]{81.6,39;$r_4$}
\lbl[]{100.7,39;$r_5$}
\lbl[]{119.8,39;$r_6$}
\lbl[]{138.9,39;$r_7$}
\lbl[]{158,39;$r_8$}
\end{lpic}
\end{figure}
such that $r_i$ is adjacent to $r_j$ for some $i<j$ if the $i$-th bit of the binary representation of $j$ is 1.

For instance,  vertex $r_0$ is adjacent to all $r_n$ with odd $n$, because the numbers whose 0-th bit is nonzero are exactly the odd numbers.
Vertex $r_1$ is adjacent to $r_0$ (since 1 is odd) and to all $r_n$ with $n\equiv 2$ or $3 \pmod 4$;
and so on.

\begin{thm}{Exercise}\label{ex:rado-costructive}
Show that the described graph is Rado.
\end{thm}

\section{Uniqueness}

In this section, we will prove that any two Rado graphs are isomorphic, so essentially, there is only one Rado graph.
First, let us prove a simpler statement.

\begin{thm}{Theorem}\label{thm:rado-subgraph}
Let $R$ be a Rado graph.
Then any countable graph $G$ (finite or infinite) is isomorphic to an induced subgraph of $R$. 
\end{thm}

\parit{Proof.}
Enumerate the vertices of $G$ as $v_1 , v_2 , \dots$ (the sequence might be finite or infinite).

It is sufficient to construct a sequence $r_1,r_2,\dots$ of vertices in $R$ such that $r_i$ is adjacent to $r_j$ if and only if $v_i$ is adjacent to $v_j$.
In this case, the graph $G$ is isomorphic to the subgraph of $R$ induced by $\{r_1,r_2\dots\}$.

We may choose any vertex of $R$ as $r_1$.
Suppose that the sequence $r_1,\dots,r_n$ is constructed.
If $G$ has $n$ vertices, then the required sequence is already constructed.
Otherwise, note that the Rado property implies that there is a vertex $r_{n+1}$ in $R$ that is adjacent to $r_i$ for $i\le n$ if and only if $v_{n+1}$ is adjacent to $v_i$.

Clearly, the new vertex $r_{n+1}$ meets all the required properties.
Repeating this procedure infinitely many times, or until the sequence $(v_n)$ terminates, produces the required sequence $(r_n)$.
\qeds

\begin{thm}{Exercise}\label{ex:rado-path10}
Show that any two vertices in a Rado graph can be connected by a path of length 10.
\end{thm}

\begin{thm}{Theorem}\label{thm:rado-isom}
Any two Rado graphs $R$ and $S$ are isomorphic.

Moreover any isomorphism $f_0\:S_0\to R_0$ between finite induced subgraphs in $R$ and $S$ can be extended to an isomorphism $f\:S\z\to R$.
\end{thm}

Note that Theorem~\ref{thm:rado-subgraph} implies that $R$ is isomorphic to an induced subgraph in $S$ and the other way around --- $S$ is isomorphic to an induced subgraph in $R$.
For finite graphs these two properties would imply that the graphs are isomorphic; see Exercise~\ref{ex:finite-subgraphs}.
As the following example shows, it does not hold for infinite graphs.
It is instructive to understand this example before going into the proof.

\begin{figure}[ht!]%{r}{20 mm}
\vskip-0mm
\centering
\includegraphics{mppics/pic-72}
\vskip-0mm
\end{figure}
 
The first graph $T$ on the diagram has an infinite number of vertices, non of which has degree 3.
The second graph $T'$ has exactly one vertex of degree 3. 
Therefore, these two graphs are not isomorphic.

Deleting the marked vertices from one graph produces the other one. 
Therefore, $T$ is isomorphic to a subgraph of $T'$ and the other way around.

The proof below uses the same construction as in the proof of Theorem~\ref{thm:rado-subgraph}, but it is applied \textit{back and forth} to ensure that the constructed subgraphs contain all the vertices of the original graph.

\parit{Proof.}
Once we have proved the second statement,
the first statement will follow if you apply it to single-vertex subgraphs $R_0$ and $S_0$.

Since the graphs are countable,
we can enumerate the vertices of $R$ and $S$, as $r_1 , r_2 , \dots$ and $s_1, s_2,\dots$, respectively. 
We will construct a sequence of induced subgraphs $R_n$ in $R$ and $S_n$ in $S$ with a sequence isomorphisms $f_n\:R_n\to S_n$.

Suppose that an isomorphism $f_n\:R_n\to S_n$ is constructed. 

If $n$ is even, set $m$ to be the smallest index such that $r_m$ not in $R_n$.
The Rado property guarantees that there is a vertex $s_k$ such that for any vertex $r_i$ in $R_n$, $s_k$ is adjacent to $f_n(r_i)$ if and only if $r_m$ is adjacent to $r_i$.
Set $R_{n+1}$ to be the graph induced by vertices of $R_n$ and $r_m$;
further set $S_{n+1}$ to be the graph induced by vertices of $S_n$ and $s_k$.
The isomorphism $f_n$ can be extended to the isomorphism $f_{n+1}\:R_{n+1}\to S_{n+1}$ by
setting $f_{n+1}(r_m)=s_k$.  

If $n$ is odd, we do the same, but backwards.
Let $m$ be the smallest index such that $s_m$ not in $S_n$.
The Rado property guarantees that there is a vertex $r_k$ which is adjacent to a vertex $r_i$ in $R_n$ if and only if $f_n(r_i)$ is adjacent to $s_m$.
Set $R_{n+1}$ to be the graph induced by vertices of $R_n$ and $r_k$;
further set $S_{n+1}$ to be the graph induced by vertices of $S_n$ and $s_m$.
The isomorphism $f_n$ can be extended to the isomorphism $f_{n+1}\:R_{n+1}\to S_{n+1}$ by
setting $f_{n+1}(r_k)=s_m$.

Note that if $f_n(r_i)=s_j$, then $f_m(r_i)=s_j$ for all $m\ge n$.
Therefore, we can define $f(r_i)=s_j$ if $f_n(r_i)=s_j$ for some $n$.

By construction we get that 
\begin{itemize}
\item $f_n(r_i)$ is defined for any $n>2\cdot i$. Therefore, $f$ is defined at any vertex of $R$.
\item $s_j$ lies in the range of $f_n$ for any $n>2\cdot j$.
Therefore, the range of $f$ contains all the vertices of $S$.
\item $r_i$ is adjacent to $r_j$ if and only if $f(r_i)$ is adjacent to $f(r_j)$.
\end{itemize}
Therefore, $f\:R\to S$ is an isomorphism.
\qeds

\begin{thm}{Exercise}\label{ex:rado-isom-generalization}
Explain how to modify the proof of \ref{thm:rado-isom} to prove the following theorem.
\end{thm}

\begin{thm}{Theorem}
Let $R$ be a Rado graph.
A countable graph $G$ is isomorphic to a spanning subgraph of $R$ if and only if, given any finite set $V$ of vertices of $G$, there is a vertex $w$ that is not adjacent to any vertex in $V$.
\end{thm}

\begin{thm}{Exercise}\label{ex:Rv>w}
Let $v$ and $w$ be two vertices in a Rado graph $R$.
Show that there is an isomorphism from $R$ to itself that sends $v$ to $w$.
\end{thm}

\begin{thm}{Exercise}\label{ex:finite-subgraphs}
Let $G$ and $H$ be two finite graphs.
Assume $G$ is isomorphic to a subgraph of $H$ 
and the other way around --- $H$ is isomorphic to a subgraph of $G$.
Show that $G$ is isomorphic to $H$.
\end{thm}

\section{The random graph}

The following theorem explains why a Rado graph is also named \textit{the random graph}.

\begin{thm}{Theorem}\label{thm:the-random-graph}
Suppose $V$ is an infinite countable set of vertices.
Let us connect each pair in $V$ by an edge independently and with probability $\tfrac12$.
Then, with probability 1, we obtain a Rado graph.
\end{thm}

\parit{Proof.}
It is sufficient to show that for two given finite sets of vertices $V$ and $W$, the Rado property fails with probability 0.

Assume $n=|V|+|W|$; that is, $n$ is the total number of vertices in $V$ and $W$.
The probability that a given vertex $v$ outside of $V$ and $W$ satisfies the Rado property for $V$ and $W$
is $\tfrac1{2^n}$.
Therefore, probability that a given vertex $v$ does \textit{not} satisfy this property
is $1-\tfrac1{2^n}$.

Note that events that a given vertex does not satisfy the property are independent.
Therefore, the probability that $N$ different vertices $v_1,\dots,v_N$ outside of $V$ and $W$ do not satisfy the Rado property for $V$ and $W$ is 
\[(1-\tfrac1{2^n})^N.\]
This value tends to 0 as $N \to \infty$; 
therefore the event that no vertex is correctly joined has probability 0.
\qeds

\begin{thm}{Exercise}\label{ex:rado-radnom}
Let $0<\alpha<1$.
Assume an infinite countable graph is chosen at random by selecting edges independently with probability $\alpha$ from the set of 2-element subsets of the vertex set.
Show that with probability 1, the resulting graph is a Rado graph.
\end{thm}

\section{Remarks}

The {}\emph{Rado graph} is also called the \index{Erd\H{o}s–R\'enyi graph}\emph{Erd\H{o}s–R\'enyi graph} or \index{random graph}\emph{random graph};
it was first discovered by Wilhelm Ackermann, rediscovered later by
Paul Erd\H{o}s and Alfr\'ed R\'enyi, and yet again by Richard Rado. 
An analogous object in metric geometry was discovered much earlier by Pavel Urysohn. 

Theorem \ref{thm:the-random-graph} was proved by Paul Erd\H{o}s and Alfr\'ed R\'enyi.
A good survey on the subject is written by Peter Cameron~\cite{cameron}.

\chapter{Rewriting systems}

Suppose a certain set of objects and a set of elementary transformations for these objects are given, and we are interested in the following problem.

\textit{Under what conditions do these transformations reduce each object to a single form?}

As we will see, this question is equivalent to an elegant statement about graphs (Lemma \ref{lem:diamond}).

For a precise formulation, we need to introduce rewriting systems.
They are understood in two senses: narrow and broad.
In the narrow sense, certain word transformations are considered.
These are the so-called \textit{string rewriting systems}, which will be defined in the next section.
In the broad sense, we talk about any objects and transformations; these are called \textit{abstract rewriting systems}.

\section{String rewriting systems}

Let us choose a finite set of symbols and call it an \index{alphabet}\emph{alphabet}.
Any finite strings (including the empty string) of symbols in the alphabet will be called a \index{word}\emph{word}.
If the alphabet has only two symbols $\{a,b\}$, then the words are 
$\emptyset$, $a$, $b$, $aa$, $ab$, $ba$, $bb$, $aaa$, and so on.

Suppose a list of pairs of words, $(s_1,t_1),\dots (s_n,t_n)$, is given;
we are allowed to change any occurrence of $s_i$ to $t_i$.
More precisely, a word $x=ls_ir$ can be exchanged for the word $y=lt_ir$;
here $l$ and $r$ are arbitrary words.
In this case we write $x\to y$.

For instance, we might be allowed to change $s_1=ba$ to $t_1=ab$ and $s_2=aba$ to $t_2=\emptyset$ anywhere in the word.
All possible words that can be reached from $bbaab$ are shown on the diagram.

{

\begin{wrapfigure}{o}{40 mm}
\vskip-0mm
\centering
\begin{tikzpicture}[node distance=1.5cm, auto]
\node (bbaab) {$bbaab$};
\node (babab) [below of=bbaab, node distance=.9cm] {$babab$};
\node (abbab) [below of=babab, node distance=.9cm] {$abbab$};
\node (baabb) [left of=babab] {$baabb$};
\node (ababb) [left of=abbab] {$ababb$};
\node (aabbb) [below of=ababb, node distance=.9cm] {$aabbb$};
\node (b) [right of=aabbb] {};
\node (bb) [right of=b] {$bb$};
\draw[->] (bbaab) to (babab);
\draw[->] (babab) to (abbab);
\draw[->] (babab) to (baabb);
\draw[->] (babab) to (bb);
\draw[->] (ababb) to (bb);
\draw[->] (baabb) to (ababb);
\draw[->] (abbab) to (ababb);
\draw[->] (ababb) to (aabbb);
\end{tikzpicture}
\vskip-2mm
\end{wrapfigure}

We have just described the so-called \index{string rewriting system}\emph{string rewriting systems}.
To describe the given example, we will use the notation
\[\langle\, a,b \mid ba\to ab,\, aba\to\emptyset\,\rangle.\]
So, after $\langle$, we list the symbols of the alphabet, and after $\mid$, we list all allowed replacements, closing it with $\rangle$.

}

\section{Abstract rewriting systems}

An \index{abstract rewriting system}\emph{abstract rewriting system} is a \index{digraph}\emph{digraph};
that is, a pseudograph (typically infinite) on each edge of which one of two directions is chosen.
The vertices of the digraph are often referred to as \index{object}\emph{objects}, and directed edges as \index{transformation}\emph{transformations}.
For two objects $x$ and $y$, we write $x \to y$ if there is a directed edge that starts at vertex $x$ and ends at vertex $y$; it defines a relation on the set of objects.

We will be interested in the relation $\to$ on the set of vertices.
Of course, if we know the digraph, then we know the relation $\to$.

A string rewriting system can be interpreted as an abstract rewriting system with the set of vertices formed by
all words in the alphabet.

Let us introduce two new binary relations on the set of vertices of the digraph that will be denoted by $\rightsquigarrow$ and $\sim$.

The first relation $x\rightsquigarrow y$ means that following the arrows in the digraph,
one can get from $x$ to $y$; it includes the case $x=y$.
More precisely, there is a sequence of objects $x=x_0\to x_1\to\dots\to x_n=y$ for some integer $n\ge 0$.

The expression $x \sim y$ says that one can walk along edges in the digraph from $x$ to $y$;
it is allowed to go in both directions of the edges.
In other words, $x \sim y$ means that $x$ and $y$ lie in the same connected component of the undirected pseudograph.

\section{Terminating systems}

A system is called \index{terminating system}\emph{terminating} if there are no infinite sequences $x_0 \to x_1 \to \dots $
In other words, one cannot travel around the digraph indefinitely, following the directions of the arrows.

An object $w$ (a vertex) is said to be \index{irreducible object}\emph{irreducible} if there is no object $x$ such that $w \to x$.
In other words, there is no way out of $w$.

It is clear that in a terminating rewriting system, starting at vertex $x$ and traveling along the arrows, we will arrive at some irreducible vertex, say~$w$.
In the latter case, we will say that $w$ is a \index{normal form}\emph{normal form} of $x$.

Any object in the terminating system has at least one normal form, and there can be many of them.
If $x$ has a {}\emph{unique normal form}, then it will be denoted by \index{$[x]$ (normal form)}$[x]$.

\begin{thm}{Exercise}\label{ex:examples}
Give an example of 
\begin{enumerate}[(a)]
\item a finite nonterminating system.
\item a terminating system such that a vertex $v_0$ has an arbitrary long sequence $v_0 \to v_1 \to \dots \to v_n$.
\item an infinite nonterminating  system such that any object has a unique normal form.
\end{enumerate}

\end{thm}

\section{Confluent systems}

\begin{wrapfigure}{r}{20 mm}
\vskip-0mm
\centering
\begin{tikzpicture}[node distance=1.1cm, auto]
\node (x) {$x$};
\node (zz) [below of=x] {$z$};
\node (yy) [right of=x] {$y$};
\node (v) [below of=yy] {$w$};
\draw[->,decorate,decoration=zigzag] (x) to (zz);
\draw[->,decorate,decoration=zigzag] (x) to (yy);
\draw[->,decorate,decoration=zigzag] (yy) to (v);
\draw[->,decorate,decoration=zigzag] (zz) to (v);
\end{tikzpicture}
\vskip-0mm
\end{wrapfigure}

A system is called \index{confluent system}\emph{confluent} if the conditions $x\rightsquigarrow  y$ and $x\rightsquigarrow  z$
for any three objects $x$, $y$, $z$ imply that
there is an object $w$ such that $y\rightsquigarrow  w$ and $z\rightsquigarrow  w$.
Informally, it means that if two people start at $x$ and follow the arrows for some time,
then, going further along the arrows, they can always converge at some vertex $w$.

\begin{thm}{Lemma}\label{lem:x->y}
Suppose two objects $x$ and $y$ in a confluent system have normal forms 
$v$ and $w$, respectively.
If $x\rightsquigarrow y$, then $v=w$.
\end{thm}

\parit{Proof.}
Since $x\rightsquigarrow y$ and $x\rightsquigarrow v$, we get $v\rightsquigarrow v'$ and $y\rightsquigarrow v'$ for some object $v'$.
Since $v$ is irreducible, we get $v=v'$;
that is $y\rightsquigarrow v$.

Since $y\rightsquigarrow v$ and $y\rightsquigarrow w$, we get $v\rightsquigarrow w'$ and $w\rightsquigarrow w'$ for some object $w'$.
Since $v$ and $w$ are irreducible, we get $v=w'=w$.
\qeds

A terminating confluent system is called \index{convergent}\emph{convergent}.

\begin{thm}{Theorem}
In a convergent system, any object has a unique normal form.
Moreover, for any two objects $x$ and $y$ in the system we have
\[x\sim y\qquad\Longleftrightarrow\qquad [x]=[y].\]

\end{thm}

\parit{Proof.}
Since the system is terminating, any object $x$ has a normal form.

Let $y$ and $z$ be normal forms of $x$;
in particular $x\rightsquigarrow y$ and $x\rightsquigarrow z$.
Since the system is confluent, $y\rightsquigarrow w$ and $z\rightsquigarrow w$ for some object $w$.
However, since both objects $x$ and $y$ are irreducible, we have $y= w$ and $z= w$;
that is, any object $x$ has a unique normal form, which proves the main statement.

Suppose $w=[x]=[y]$.
Then $x\rightsquigarrow w$ and $y\rightsquigarrow w$.
Therefore, $x\sim y$;
we proved the if part of the second statement.

Now suppose $x\sim y$;
so there is a sequence $x=x_0,x_1,\dots,x_n=y$ such that $x_i\to x_{i+1}$ or $x_{i+1}\to x_i$ for each $i$.
By \ref{lem:x->y},
\[[x]=[x_0]=[x_1]=\dots=[x_n]=[y].\]
Hence the only-if part follows.
\qeds

\section{Proving termination}

So, how do we prove that a system is terminating?
Typically, we construct a relation ``$\succ$'' on the set of vertices such that 
$x \to y$ implies $x\succ y$.
If any decreasing sequence with respect to $\succ$ is terminating,
then our original system is terminating as well. 

More precisely, ``$\succ$'' is a strict partial order;
that is, ``$\succ$'' is a relation that satisfies the following conditions for all objects $x$, $y$, and $z$:
\begin{itemize}
\item $x\nsucc x$.
\item If $x\succ y$, then  $y\nsucc x$.
\item If $x\succ y$ and $y\succ z$ then  $x\succ z$.
\end{itemize}
If there is no infinite sequence $x_0\succ x_1\succ\dots$,
and $x\to y$ implies $x\succ y$, then our system is terminating.

\begin{thm}{Exercise}\label{ex:x+1}
Let $A$ be an abstract rewriting system
with a set of objects formed by positive integers 
and the transformations $x\to \tfrac x2$ if $x$ is even and $x\to x+1$ if $x>1$ is odd.
Show that $A$ is terminating.
\end{thm}

For string rewriting systems, the words can be ordered first by length, and if the lengths are equal, lexicographically.
This is the so-called \index{shortlex order}\emph{shortlex order}.
Note that for this order, there is no infinite decreasing sequence.
In this case, to ensure that the system is terminating, it is sufficient that when moving from one word to another along the arrow, either the length of the word decreases, or it does not change, but the word decreases relative to the order of words in the dictionary.

\begin{thm}{Exercise}\label{ex:balls}
Show that the following process always terminates.
There is a box with a finite number of black and white balls.
Each step consists of removing an arbitrary ball from the box.
If it happens to be a black ball, one also adds an arbitrary (but finite) number of white balls to the box.
\end{thm}

\begin{thm}{Exercise}\label{ex:ab>ba,aba>}
Show that the following string rewriting system 
\[\langle\, a,b \mid ba\to ab,\, aba\to\emptyset\,\rangle\]
is terminating.
\end{thm}

\section{Proving confluence}

Checking confluence is a more difficult problem.
One needs to keep track of all pairs of paths leaving the same vertex.
It turns out that this is not necessary.

\begin{wrapfigure}{r}{20 mm}
\vskip-0mm
\centering
\begin{tikzpicture}[node distance=1.1cm, auto]
\node (x) {$x$};
\node (zz) [below of=x] {$z$};
\node (yy) [right of=x] {$y$};
\node (v) [below of=yy] {$w$};
\draw[->] (x) to (zz);
\draw[->] (x) to (yy);
\draw[->,decorate,decoration=zigzag] (yy) to (v);
\draw[->,decorate,decoration=zigzag] (zz) to (v);
\end{tikzpicture}
\vskip-0mm
\end{wrapfigure}

An abstract rewriting system is called \index{locally confluent}\emph{locally confluent} if for any objects $x$, $y$, $z$, from the conditions $x \to y$, $x\to z$, follows the existence of an object $w$ such that $y\rightsquigarrow  w$, $z\rightsquigarrow  w$.
This means that if two people started from the same vertex and each walked along one directed edge, then they can converge somewhere, following the arrows.

It is clear that every confluent system is locally confluent.

\begin{thm}{Exercise}\label{ex:not-confluent}
Construct a system that is locally confluent but not confluent.
\end{thm}

\begin{thm}{Diamond lemma}\label{lem:diamond}
A terminating rewriting system is confluent if and only if it is locally confluent.
\end{thm}

Let us describe an induction-type argument that will be used in the proof.

Suppose we want to prove a statement $P(x)$ for every object $x$ of the terminating system.
Then it is sufficient to prove the induction step:
\textit{$P(x)$ holds under the assumption that $P(y)$ holds for all downstream objects $y$}.
In other words, we need to show that \textit{if $P(y)$ holds for any object $y\ne x$ such that $x\rightsquigarrow y$, then $P(x)$ holds as well}.

Let us explain why it works.
Suppose we proved the induction step for some statement $P$ in a terminating system, but $P(x_0)$ does not hold for an object $x_0$.
Then we may choose another object $x_1\ne x_0$ such that $P(x_1)$ does not hold and $x_0\rightsquigarrow x_1$.
Continuing this way, we get an infinite descending sequence $x_0\rightsquigarrow x_1\rightsquigarrow x_2\dots$ of so that $x_0\ne x_1$, $x_1\ne x_2$, and so on.
The latter contradicts that our system is terminating.

Note that applying the step of induction to the empty set of objects,
we get that $P(x)$ holds for any irreducible object, and we can continue further.
It explains why this version of induction requires no base.

\parit{Proof.}
The only-if part is trivial;
let us prove the if part.

Given an object $x$, we need to show that if $x\rightsquigarrow y$ and $x\rightsquigarrow  z$, then there is $w$ such that $y\rightsquigarrow w$ and $z\rightsquigarrow  w$.
Denote the latter statement by $P(x)$.

If $x=y$ or $x=z$, then there is nothing to prove.
Therefore, we can assume that $x\ne y$ and $x\ne z$.
Consider the first edges $x\to y'$ and $x \to z'$ of the paths defining $x\rightsquigarrow y$ and $x\rightsquigarrow  z$; 
so we have $x\to y'$, $y' \rightsquigarrow  y$, $x \to z'$, and $z' \rightsquigarrow  z$.

\begin{wrapfigure}[9]{r}{40 mm}
\vskip-6mm
\centering
\begin{tikzpicture}[node distance=1.5cm, auto]
\node (x) {$x$};
\node (zz) [below of=x] {$z'$};
\node (z) [below of=zz] {$z$};
\node (zzz) [right of=z] {$z''$};
\node (yy) [right of=x] {$y'$};
\node (y) [right of=yy] {$y$};
\node (yyy) [below of=y] {$y''$};
\node (v) [below of=yy] {$v$};
\node (w) [below of=yyy] {$w$};
\draw[->] (x) to (zz);
\draw[->] (x) to (yy);
\draw[->,decorate,decoration=zigzag] (yy) to (y);
\draw[->,decorate,decoration=zigzag] (yy) to (v);
\draw[->,decorate,decoration=zigzag] (zz) to (v);
\draw[->,decorate,decoration=zigzag] (zz) to (z);
\draw[->,decorate,decoration=zigzag] (z) to (zzz);
\draw[->,decorate,decoration=zigzag] (v) to (zzz);
\draw[->,decorate,decoration=zigzag] (v) to (yyy);
\draw[->,decorate,decoration=zigzag] (y) to (yyy);
\draw[->,decorate,decoration=zigzag] (zzz) to (w);
\draw[->,decorate,decoration=zigzag] (yyy) to (w);
\end{tikzpicture}
\vskip-0mm
\end{wrapfigure}

Applying the local confluence condition to a pair of edges $x \to y'$, $x \to z'$, we find a vertex $v$ such that $y' \rightsquigarrow  v$ and $z' \rightsquigarrow  v$.

Next, the vertices $y'$ and $z'$ are descending with respect to $x$, and the same applies to vertex $v$. Therefore, by the induction hypothesis, the statements $P(y')$, $P(z')$, and $P(v)$ hold true.

Since $y' \rightsquigarrow  y$ and $y' \rightsquigarrow  v$ we can get $y''$ such that $y\rightsquigarrow  y''$ and $v\rightsquigarrow y''$.
Likewise, we get $z''$ such that $z\rightsquigarrow  z''$ and $v\rightsquigarrow  z''$.

Finally, we use $P(v)$.
Since $v\rightsquigarrow  y''$ and $v \rightsquigarrow z''$, we can find $w$ such that $y'' \rightsquigarrow  w$ and $z'' \rightsquigarrow  w$.

Finally, $y \rightsquigarrow y'' \rightsquigarrow  w$ implies $y \rightsquigarrow  w$.
Likewise, we get $z \rightsquigarrow  w$.
\qeds

\section{Applications}

Suppose $x$ is a word in a string rewriting system
from which two arrows emanate $x\to y$ and $x\to z$.
So, $y$ and $z$ are obtained from $x$ by two elementary transformations, replacing some of its occurrences $s_i$ and $s_j$ with $t_i$ and $t_j$, respectively.
If occurrences do not overlap, then we can apply these two transformations in any order and get the same result, say $w$.
In this case, local confluence holds.

If the occurrences overlap, then two cases are possible:
they can either be contained in one another (a rarer case) or overlap by some non-empty word.
Without loss of generality, we can pass to the part of our word that is the union of occurrences $s_i$ and $s_j$.
The following two cases are subject to verification.

\parit{Case A.} Suppose $s_{i}$ and $s_j$ overlap by non-empty word $q$;
we may assume that  $s_{i}=lq$ and  $s_{j}=qr$, for some words $l$ and $r$.
That is, we have a word $lqr$, from which we can go along the edge to both the word $t_{i}r$ and the word $lt_{j}$.
By making transformations in these words, we want to reduce them to one word; that is, we need to find a word $w$ such that $t_{i}r\rightsquigarrow  w$ and $lt_{j}\rightsquigarrow w$.

\parit{Case B.} Suppose  $s_{j}$ is contained in $s_{i}$;
in other words $s_i=ls_{j}r$ for some words $l$ and $r$.
In this case, from the vertex $s_{i}$ we can move along the edge to both $t_{i}$ and $lt_{j}r$.
The goal is to find a word $w$ such that $t_{i}\rightsquigarrow w$, $lt_{j}r\rightsquigarrow w$.

\begin{thm}{Example}\label{em:baaba>a}
\[\langle a,b \mid baaba\to a\rangle.\]
\end{thm}

Applying the shortlex order, we find that the system is terminating.

To check confluence, we need to examine cases A and B.
There is only one word here, but it overlaps with itself in $ba$.
By gluing together two such words that overlap in $ba$, we come to checking the word $baabaaba$.
Applying the transformation at the left and right ends, we get two irreducible words
$aaba$ and $baaa$.
It follows that the system is not confluent and, therefore, not convergent.

\medskip

In this example, the words $aaba$ and $baaa$ are equivalent to the original word $baabaaba$.
Therefore, if we add $baaa\to aaba$ to the existing rules, then the word equivalence relation will remain unchanged.
It leads to the following example.

\begin{thm}{Example}\label{em:baaba>a,baaa>aaba}
\[\langle a, b \mid baaba\to a, baaa\to aaba\rangle.\]
\end{thm}

Again, applying the shortlex order, we find that the system is terminating.

Now we need to go thru all the options for overlapping words.
There are only two such options.
One of them is old and does not need to be checked.
Now the words $aaba$ and $baaa$ from the previous example are automatically converted to $aaba$ in the new system.
But we have a new relationship, and therefore a new overlap has arisen: $baaba$ overlaps with $baaa$ by $ba$. (Here you must always be careful, since a new word can give several new overlaps;
two words can overlap not in one, but in many ways --- all this must be checked.)

Gluing the overlapping words together leads to $baabaaa$
and we can get $(baaba)aa\to aaa$
and $baa(baaa)\to baaaaba$.
Since $(baaa)aba \z\to  aa(baaba) \to aaa$, we get that the system is confluent,
and therefore convergent.

\begin{thm}{Exercise}\label{ex:convergent}
Check if the following systems are terminating and conluent.

\begin{enumerate}[(a)]
 \item\label{ex:convergent:a} \[\langle a, b \mid baaba\to a, aaba \to  baaa\rangle;\]
 \item\label{ex:convergent:b} \[\langle a,b,c \mid ba\to ab, ca\to ac, cb\to bc \rangle;\]
 \item\label{ex:convergent:c} \[\left\langle
\begin{matrix}
x,&y,&z,
\\
X,&Y,&Z
\end{matrix}
\,
\middle| 
\,
\begin{matrix}
xX\to \emptyset,& yY\to \emptyset,& zZ\to \emptyset,
\\
Xx\to \emptyset,& Yy\to \emptyset,& Zz\to \emptyset
\end{matrix}
\,
\right\rangle.\]
\end{enumerate}

\end{thm}

\section{Comments}

This chapter is inspired by a blog post by Victor Guba \cite{guba}.
Several examples and exercises are taken from the book by Franz Baader and Tobias Nipkow \cite{baader-nipkow},
which is very good for further study.

The diamond lemma (\ref{lem:diamond}) is also known as \index{Newman's lemma}\emph{Newman's lemma};
it was proved by Max Newman \cite[Theorem 3]{newman}.

The replenishment of the system in \ref{em:baaba>a} to the system in \ref{em:baaba>a,baaa>aaba} is called  the \index{Knuth--Bendix procedure}\emph{Knuth--Bendix procedure}.
It may produce a convergent system after several iterations, but this process may also run forever;
in our example, just one step was needed.

\appendix

%\input{extra.tex}
%+Shannon switching game/Hex game???
%+De Bruijn--Erdős theorem???

\newgeometry{top=0.9in, bottom=0.9in,inner=0.55in, outer=0.45in}
%\newgeometry{top=1.025in, bottom=1.025in,inner=0.5in, outer=0.625in, paperwidth=6.125in, paperheight=9.25in}%with bleed
{\footnotesize

\backmatter

\chapter{Hints}

\raggedcolumns\setlength{\multicolsep}{10mm}
\spell{\begin{multicols}{2}}{}

\refstepcounter{chapter}
\setcounter{eqtn}{0}

\parbf{\ref{ex:cannibals}.}
Here is a part of the needed graph.

\begin{center}
\begin{tikzpicture}[scale=1.8,
  thick,main node/.style={circle,draw,font=\sffamily\bfseries,minimum size=1mm}]
  \node[main node] (1) at (0,0) {{\small${4^*_4}{\parallel}{0_0}$}};
  \node[main node] (2) at (1,-1){{\small${2_4}{\parallel}{2^*_0}$}};
  \node[main node] (3) at (0,-1){{\small${3_3}{\parallel}{1^*_1}$}};
  \node[main node] (4) at (1,-2){{\small${3^*_4}{\parallel}{1_0}$}};
    \node[main node] (5) at (-1,-1){{\small${3_4}{\parallel}{1^*_0}$}};
  \path[every node/.style={font=\sffamily\small}]
   (1) edge node[auto]{$2_0$}(2)
   (3) edge node[auto]{$1_1$}(1)
   (5) edge node[auto]{$1_0$}(1)
   (4) edge node[auto]{$0_1$}(3)
   (2) edge node[auto]{$1_0$}(4);
\end{tikzpicture}
\end{center}

\refstepcounter{chapter}
\setcounter{eqtn}{0}

\parbf{\ref{ex:r(2,n)}.}
Use the following observation:
\textit{if there is no blue $K_2$, then all edges are red.}

\parbf{\ref{ex:K8+K17}.} 
Assuming such a subgraph exists.
Fix a vertex $v$.
Note that we can assume that the subgraph contains $v$; otherwise rotate it.
In each case, draw the subgraph induced by the vertices connected to $v$.
(If uncertain, see definition of \textit{induced subgraph}.)

\refstepcounter{chapter}
\setcounter{eqtn}{0}

\parbf{\ref{ex:number(ham-cycles)}.}
Show that $K_{100}$ has $\tfrac{99!}2$ Hamiltonian cycles.
Find the probability $P$ that that a given Hamiltonian cycle is monochromatic.
Show that $N\cdot P$ is the expected number of Hamiltonian cycles.
Estimate $N\cdot P$; you can use
Stirling's inequality
\[n!>\sqrt{2\cdot \pi\cdot n}\cdot \left(\frac{n}{e}\right)^n.\]

\parbf{\ref{ex:Qn-dist};} \textit{\ref{Pn}}.
Use that the distance between two vertices is the number of different digits in their sequences.

\parit{\ref{kPn}}. Use that expected value of sum is the sum of expected values. 

\parit{\ref{ex:Qn-dist:end}}. Show and use that $1.05\cdot0.95<1$.

\parbf{\ref{ex:lin-Qn};} \textit{\ref{ex:lin-Qn:n+1}}.
Identify the vertices of $Q_n$ by $n$-dimensional vectors with $\pm1$ components.
Observe that two vertices $v$ and $w$ are at a distance $>\tfrac n2$ if and only if $\langle v,w\rangle <0$;
here $\langle\ ,\ \rangle$ denotes the scalar product (also known as dot product).

Suppose $v_1,\dots,v_k$ are $n$-dimensional vectors such that $\langle v_i,v_j\rangle <0$ if $i< j$.
Denote by $w_i$ the projection of $v_i$ to the hyperplane perpendicular to $v_k$.
Show that $\langle w_i,w_j\rangle <0$ if $i<j<k$.
Apply it with induction by $k$.

\parit{\ref{ex:lin-Qn:2n}}
Choose a hyperplane $H$ thru the center of the cube that does not pass thru any vertices of the cube.
It subdivides the vertices in our set into two collections.

Project one of these collection to $H$ and argue as in \textit{\ref{ex:lin-Qn:n+1}} to show that it has at most $n$ vertices.

\parit{Remark.}
The equality in \textit{\ref{ex:lin-Qn:2n}} holds only if $n$ is a power of 2.
Try to prove it.

\parbf{\ref{ex:prob(isom)}.} Identify vertices of $H$ and $G_n$.
Calculate probability that $H=G_n$; it gives the first inequality.
Note that there are $n!$ ways to identify vertices of $H$ and $G_n$;
it implies the second inequality.

The first inequality becomes equality if permuting vertices of $H$ gives the same graph.
There only two graphs of that type; try to find both.

The second inequality becomes equality if $H$ has no symmetries;
that is, there is no nontrivial permutation of vertices that induces an isomorphism from $H$ to itself.
Try to construct such a graph with 10 vertices.

\parbf{\ref{ex:diam=2}.}
Let $P_n$ be the probability that two vertices $v$ and $w$ in $G_n$ cannot be connected by a path of length $2$.

Find $P_n$.
Show that $\binom n2\cdot P_n$ is the expected number of such pairs.
Apply Markov's inequality.
You should get that a typical graph has diameter at most 2.

Finally, find probability that $G_n$ has diameter $<2$.
Make a conclusion.

\parbf{\ref{ex:typ(K100)}.}
Denote by $P$ the probability that the subgraph induced by $100$ vertices in $G_n$ is $K_{100}$;
note that $P>0$.

Choose disjoint 100-element subsets $S_1,\dots,S_k$ of vertices in $G_n$.
We may assume that $k>\tfrac n{100}-1$; in particular, $k\to \infty$ as $n\to \infty$.

Show that probability that at least one of $S_i$ induces a subgraph isomorphic to $K_{100}$ is at least 
\[1-(1-P)^k.\]
Make a conclusion.

\refstepcounter{chapter}
\setcounter{eqtn}{0}

\parbf{\ref{ex:bridge}.}
Show and use that any spanning tree of $G$ contains the bridge.

\parbf{\ref{ex:zig-zag}.}
Use induction on $n$ and/or mimic the proof of \ref{thm:fans} with the following analogs of $Z_n$:

\begin{Figure}
\centering
\includegraphics{mppics/pic-371}
\end{Figure}

\begin{wrapfigure}{r}{20 mm}
\vskip-4mm
\centering
\includegraphics{mppics/pic-39}
\end{wrapfigure}

\parbf{\ref{ex:ladder}.}
To construct the recursive relation, in addition to the ladders $L_n$, you will need two of its analogs --- $L_n'$ and $L_n''$, shown on the diagram.

\parbf{\ref{ex:wheel}.} Start with a spoke of the wheel.
The diagram will contain wheels, and several analogs of fans.

\refstepcounter{chapter}
\setcounter{eqtn}{0}

\parbf{\ref{ex:n(walks)}.}
Apply induction on $n$.

\parbf{\ref{ex:Kirchhoff-row}.}
Use property \ref{3} on page \pageref{3}.

\parbf{\ref{ex:minor>graph}.}
Revert the steps of the construction of the Kirchhoff minor.
(Since the Kirchhoff minor is $5{\times}5$, 
the graph should have 6 vertices.)

\parbf{\ref{ex:sum-kirchhoff}.}
Apply \ref{ex:Kirchhoff-row}.

\parbf{\ref{ex:K33W6Q3}.}
Apply the construction of Kirchhoff minor.

\parbf{\ref{ex:det}.}
Apply property \ref{3} twice and property \ref{2} (see page \pageref{3}).

\parbf{\ref{ex:s(Kp-e)}.}
We can assume that $e$ connects the last two vertices of $K_p$.
In this case
\[
M=\left(
\begin{matrix}
p{-}1&-1&\cdots&-1
\\
-1&\ddots&\ddots&\vdots
\\
\vdots&\ddots&p{-}1&-1
\\
-1&\cdots&-1&p{-}2
\end{matrix}
\right)
\]
is the Kirchhoff minor of $K_p-e$.
It remains to find $\det M$.

Follow the calculations in the proof of the Cayley formula.

\parbf{\ref{ex:s(Kmn)}.}
Note that
\[
M=
\left(
\begin{matrix}
3&0&0&0&-1&-1
\\
0&3&0&0&-1&-1
\\
0&0&3&0&-1&-1
\\
0&0&0&3&-1&-1
\\
-1&-1&-1&-1&4&0
\\
-1&-1&-1&-1&0&4
\end{matrix}
\right)\]
is a Kirchhoff minor of $K_{4,3}$.

First, show that $\det M=3^3\cdot 4^2$,
and then generalize the argument to $K_{m,n}$.

\refstepcounter{chapter}
\setcounter{eqtn}{0}

\parbf{\ref{ex:PG=PHPH}.}
Use that the subgrahs $H_1$ and $H_2$ can be colored independently.

\parbf{\ref{ex:PWn}.} Let $v$ be the center of $W_n$.
Suppose we have $x+1$ choices for color of $v$;
once this choice is made, the rest of $W_n-v$ has to be colored in the remaining $x$ colors.
Show and use that each such coloring corresponds to a coloring of $C_n$ in $x$ colors. 

\parbf{\ref{ex:PGpqn}.} 
By \ref{thm:chromatic-polynomial}, $P_G$ is a monic polynomial of degree $p$.
Show that $P_H(0)=0$ for any connected graph $H$.
Use \ref{ex:PG=PHPH} to show that for any $k<n$ the coefficient of $P_G$ in front of $x^k$ vanish.

To show that $a_{p-1}=q$,
use the deletion-minus-contraction formula in an induction on $q$.

\parbf{\ref{ex:P(tree)}.} The only-if part follows from \ref{ex:PTCpFnLn}\ref{ex:PTCpFnLn:tree}.
For the if part, apply \ref{ex:PGpqn} to show that $G$ is connected, has $p$ vertices and $p-1$ edges.
Make a conclusion.

\parbf{\ref{ex:chrom(K_p)}.} Apply the induction on $p$.

\parbf{\ref{ex:P=nonisom}.} Use \ref{ex:PTCpFnLn}\ref{ex:PTCpFnLn:tree}.

\parbf{\ref{ex:MG=MHMH}.} Use that the matchings in $H_1$ and $H_2$ can be chosen independently.

\parbf{\ref{ex:matchings}.} Use the definition of $M_G$.

\parbf{\ref{ex:deletion-deletion-total};}
\textit{\ref{ex:deletion-deletion}.}
Show and use that 
$$m_n(G)=m_n(G-e)+m_{n-1}(G-[e]).$$

\parit{\ref{ex:deletion-deletion-K}.}
Apply it $p$ times to all edges at one vertex of $K_p$.

\parbf{\ref{ex:SG}.} Use the definition of $S_G$.

\parbf{\ref{ex:S(Kp)}.}
Calculate $N=\tfrac{\partial^{k-1} }{\partial x_1^{k-1}}S_{K_n}(0,1,\dots,1)$ using the expression for $S_{K_n}$ given in Theorem~\ref{thm:spanning-tree-polynomial} and determine how much a tree with degree $d$ at the first vertex contributes to the value $N$.

\parbf{\ref{ex:S(Kmn)}.}
Modify the proof of Theorem~\ref{thm:spanning-tree-polynomial}.

\refstepcounter{chapter}
\setcounter{eqtn}{0}

\parbf{\ref{ex:B=xA}+\ref{ex:B=xA'}.} Apply the definition of exponential generating function.

\parbf{\ref{ex:exp(Fn)}.}
For \textit{\ref{ex:exp(Fn):F''}},
apply the definitions of Fibonacci numbers and exponential generating function.

For \textit{\ref{ex:exp(Fn):F(x)}}, solve the differential equation in \textit{\ref{ex:exp(Fn):F''}} and use that $f_0=f_1=1$.

For \textit{\ref{ex:exp(Fn):Binet}}, use the Taylor expansion 
\[e^x=1+\tfrac x{1!}+\tfrac {x^2}{2!}+\tfrac {x^3}{3!}+\dots\]

\parbf{\ref{ex:perfect-matching}.} 
Apply \ref{ex:deletion-deletion-total}\textit{\ref{ex:deletion-deletion-K}}.

\parbf{\ref{ex:an+nan-1}.}
Apply \ref{ex:deletion-deletion-total}\textit{\ref{ex:deletion-deletion-K}} for $x=1$.

\parbf{\ref{ex:ex:2-factor};} \textit{\ref{ex:2-factor:cn}.}
Enumerate the vertices of $K_n$.
Choose a cycle and list the numbers of the vertices in the order they appear on the cycle after $n$.
Count all possible orders.
Show and use that we have exactly two different orders for one cycle.

\refstepcounter{chapter}
\setcounter{eqtn}{0}

\parbf{\ref{ex:neighbor-trees}.}
Suppose $G$ is a minimal connected graph with two cycles;
that is, removing an edge or a vertex from $G$ makes it disconnected or reduces the number of cycles.

Show that $G$ contains two spanning trees that are not neighbors.
Use it to prove the general case.

\parbf{\ref{ex:w>2w}.} Denote the two weights by $\weight_1$ and $\weight_2$.
Suppose $T$ is a spanning tree and $T'$ is its neighbor. 
Show that 
\begin{align*}
\weight_1(T)&\le \weight_1(T')
\\
&\Updownarrow
\\
\weight_2(T)&\le \weight_2(T').
\end{align*}
Apply \ref{thm:mst-iff}.

\parbf{\ref{ex:PB}.} Apply \ref{thm:mst-iff}.

\parbf{\ref{ex:KPB}.}
Compare the number of edges that have to be checked in each algorithm.
Note that the organization of data might give an essential difference.
Think which steps could be done parallelly (say, on different computers).

\parbf{\ref{ex:deleting-algorithm}.}
Part \textit{\ref{ex:deleting-algorithm:a}} should be evident.
For part \textit{\ref{ex:deleting-algorithm:b}}, use \ref{thm:mst-iff}.

\parit{\ref{ex:deleting-algorithm:c}.}
Start to walk from an arbitrary vertex without coming back.
Stop at the end vertex or at the first vertex that was visited twice.
Decide what to do in each case and start over.

\refstepcounter{chapter}
\setcounter{eqtn}{0}

\parbf{\ref{ex:bigraph-matching}.} Show and use that $P$ has odd length.

\parbf{\ref{ex:1-factor};} \textit{\ref{ex:1-factor:a}.}
Note that a 1-factor is a matching.
Apply \ref{thm:marriage}.

\parit{\ref{ex:1-factor:b}.}
Remove the 1-factor provided by \textit{\ref{ex:1-factor:a}}. Apply \textit{\ref{ex:1-factor:a}} again and again.

\parit{Remark.}
If $r=2^n$ for an integer $n\ge 1$, then $G$ has an Euler's circuit. 
Note that the total number of edges in $G$ is even, so we can delete all odd edges from the circuit.
The obtained graph $G'$ is regular with degree $2^{n-1}$.
Repeating the described procedure recursively $n$ times, 
we will end up with a 1-factor of $G$.

There is a tricky way to make this idea work for arbitrary $r$, not necessarily a power of $2$; 
it was discovered by Noga Alon [see \ncite{alon} and also \ncite{kalai}]. 

\parbf{\ref{ex:no-1-factor}.}
Assume the graph, say $G$, has a matching.
Then the central vertex $v$ is matched with one of its three neighbors, say $w$.
Show and use that $G-v-w$ has even number of vertices in each connected component.

\parbf{\ref{ex:kids}.}
Consider the bigraph with vertices labeled by rows and countries.
Connect a row-vertex to a country-vertex if a kid from the corresponding county stands in the corresponding row.
Try to apply \ref{thm:marriage}.

\parbf{\ref{ex:sons(king)}.}
Consider the bigraph with vertices labeled by 23 old parts and 24 new parts.
Connect an old-part vertex to a new-part vertex if the corresponding parts overlap.
Try to apply \ref{thm:marriage}.

\parbf{\ref{ex:nxn-table}.}
Consider the bigraph with vertices labeled by rows and columns.
Connect a row vertex to a column vertex if a positive number stands in the common cell.
Try to apply \ref{thm:marriage}.

\parbf{\ref{ex:camel17}.}
Add $s$ boys that everyone likes.
Apply \ref{thm:marriage} and remove these $s$ boys.

\parbf{\ref{ex:two-paths}.}
Assume that two $M$-alternating paths starting from $l$ and $r$ can have a common vertex~$v$.
Show and use that there is an $M$-alternating paths starting from $l$ to $r$.

\parbf{\ref{ex:rooks}.} Consider the bigraph with a vertex for each column and row.
Connect a row vertex to a column vertex if their common cell is marked.
Apply \ref{thm:vertex-cover}.

\parbf{\ref{ex:min-cut-marriage}.} 
Let $m$ be the minimal number of edges one can remove from $\hat G$ 
so that there will be no directed path from $s$ to $t$.
Show that the assumptions in \ref{thm:marriage} imply that $m=|L|$.
Construct a bijection between directed paths from $s$ to $t$ in $\hat G$ 
corresponds to an edge in to an edge in $G$.
Apply \ref{thm:mincut}.

\refstepcounter{chapter}
\setcounter{eqtn}{0}

\parbf{\ref{ex:crossing1}.}
Attach a small handle to the sphere so that one edge can go thru it over the crossing.

\parbf{\ref{ex:toroidal-graphs}.}
\textit{(a)}
Modify the diagram for $K_7$.

\parit{(b)}+\textit{(c)} Start with the following square diagrams and draw the remaining edges.
\vskip2mm
\includegraphics{mppics/pic-119}
\hskip5mm
\includegraphics{mppics/pic-118}

\parbf{\ref{ex:nonplanar-toroidal}.}
Suppose that $\Delta$ is a nonsimple region in the drawing.
Cut $\Delta$ from the sphere and glue instead a disc or two to get a sphere.

\parbf{\ref{ex:toroidal-girth}.}
Show that $q\ge 2\cdot r$ and apply
\ref{thm:euler>=}.

\parbf{\ref{ex:K5-torus}.}
Assume it is possible.
Show that the drawing has 5 squares.
Try to glue 5 squares to each other side-to-side to obtain a torus. 

\parbf{\ref{ex:forbidden-minors}.}
Show that the answers are $K_{1,5}$ and $K_4-e$.

\parbf{\ref{ex:nontoroidal}.}
Formally, one needs to show that deletion and contraction of every edge makes the graph toroidal.
(Indeed, deleting a vertex produces a minor in a graph after deleting an edge.)
But the graph has many symmetries.
It reduces the number of cases to four:
two deletions and two constructions.

\parbf{\ref{ex:K6moebius}.} 
Make a Möbius strip; better use transparent paper.
Redraw the diagram on the strip and try to draw the remaining 6 edges.

\begin{Figure}
\centering
\includegraphics{mppics/pic-125}
\end{Figure}

\refstepcounter{chapter}
\setcounter{eqtn}{0}

\parbf{\ref{ex:rado-infty}.}
Assume it is finite,
then apply the definition including all the vertices in $V\cap W$.

\parbf{\ref{ex:rado-diam}.}
Apply the definition for $V=\emptyset$ and a single vertex in $V$.
Conclude that the diameter is at least 2.

Choose two vertices $x$ and $y$.
Apply the definition for $V=\{x,y\}$ and $W=\emptyset$.
Conclude that there is a path of length 2 from $x$ to $y$,
and therefore the diameter is at most 2.

\parbf{\ref{ex:rado-partition}.}
Let $P$ and $Q$ be the induced subgraphs in the Rado graph $R$.
Assume $P$ is not Rado; that is, there is a pair of finite vertex sets $V$ and $W$ in $P$, such that any vertex $v$ in $R$ that meet the Rado property for $V$ and $W$ does not lie in $P$ (and therefore it lies in~$Q$).
Use $V$ and $W$ to show that $Q$ is Rado.

\parbf{\ref{ex:R-e-v-rev}+\ref{ex:rado-link}+\ref{ex:rado-costructive}.}
Check the definition in each case.
For \ref{ex:rado-costructive}, use the binary numeral system.

\parbf{\ref{ex:rado-path10}.}
Argue as in \ref{thm:rado-subgraph}.

\parbf{\ref{ex:rado-isom-generalization}.}
Find how to weaken the conditions at the odd steps of the construction.

\parbf{\ref{ex:Rv>w}.}
Apply \ref{thm:rado-isom}.

\parbf{\ref{ex:finite-subgraphs}.}
Show and use that the graphs have the same number of vertices and edges.

\parbf{\ref{ex:rado-radnom}.}
Choose two finite subsets of vertices $V$ and~$W$, and a vertex $v\notin V\cup W$.
Show that $v$ satisfies the definition with a fixed positive probability $\beta$.
Conclude that the definition holds with probability 1.

Finally, use that there is only countably many of choices for $V$ and $W$.

\parbf{\ref{ex:x+1}.}
Use the order 
\[1\prec 2\prec 4\prec 3\prec 6\prec 5\prec 8\prec 7\prec\dots\]

\parit{Remark.}
Proving that a given rewriting system is terminating might be difficult.
For example,
it is unknown if the following system is terminating:
the vertex set formed by all positive integers and the rules $x\to \tfrac x2$ if $x$ is even and $x\to 3\cdot x+1$ if $x>1$ is add.
This is the so-called $(3\cdot x+1)$-problem.

\parbf{\ref{ex:balls}.}
Consider the following order on pairs of integers:
$(m,n)\succ(m',n')$ in two cases: if $m>m'$, or if $m=m'$ and $n>n'$.
(In fact it is the lexicographic or on pairs.)
Apply this order assuming that the pair $(m,n)$ describes a box with $m$ white balls and $n$ black balls.

\parbf{\ref{ex:ab>ba,aba>}.}
Use the shortlex order.

\parbf{\ref{ex:not-confluent}.}
Try to construct such system with four objects,
two irreducible and two reducible.

\parbf{\ref{ex:convergent}.}
\textit{\ref{ex:convergent:a}} It is terminating but not confluent.

\parit{Remark.} Note that the systems in \ref{em:baaba>a}, \ref{em:baaba>a,baaa>aaba}, and \textit{\ref{ex:convergent:a}} have the same equivalence relation $\sim$.
The Knuth--Bendix procedure that led us from \ref{em:baaba>a} to \ref{em:baaba>a,baaa>aaba} lasts indefinitely, but all new relations that arise have the form
\[bb(ab)^{m}aaa\to a(ba)^m.\]
All of them, taken at $m=0,1,2,\dots$ together with the relation $aaba\to baaa$, form a convergent system.

\begin{wrapfigure}{r}{20 mm}
\vskip-2mm
{
\centering
\begin{tikzpicture}[node distance=1.3cm, auto]
\node (x) {$cba$};
\node (y) [right of=x] {$bca$};
\node (zz) [below of=x, node distance=.9cm] {$cab$};
\node (v) [below of=y, node distance=.9cm] {$abc$};
\draw[->] (x) to (zz);
\draw[->] (x) to (y);
\draw[->,decorate,decoration=zigzag] (y) to (v);
\draw[->,decorate,decoration=zigzag] (zz) to (v);
\end{tikzpicture}
}
\end{wrapfigure}

\parit{\ref{ex:convergent:b}}
Use lexicographic order to show that the system is terminating.
For confluence, check the following diagram and show that this is the only thing to check.

\parit{\ref{ex:convergent:c}}
Note that each rule shortens the word.
Use it to show that the system is terminating.

There are 6 possible overlaps here, but they are all of the same type.
It is sufficient to consider one word $xXx$.

\parit{Remark.}
А discerning reader could notice that this example describes a free group of rank 3;
capital letters correspond to inverses. 

\spell{\end{multicols}}{}
\newpage

}

\refstepcounter{chapter}
\setcounter{eqtn}{0}
\chapter*{Corrections and additions}
\addcontentsline{toc}{chapter}{Corrections and additions}

Here we include corrections and additions to \cite{hartsfield-ringel}.

\section{Correction to 3.2.2}

The proof of Theorem 3.2.2 about decomposition of a cubic graph with a bridge into 1-factors
does not explain why ``each bank has an odd number of vertices''.

This is true since the 1-factor containing the bridge breaks all the vertices of each bank into pairs, except for the end vertex of the bridge.

\section{Correction to 8.4.1}

There is an inaccuracy in the proof of Theorem 8.4.1 about stretchable planar graphs.
Namely, in the planar drawing of $G-h$, the region $R$ might be unbounded.

To fix this inaccuracy, one needs to prove a slightly stronger statement.
Namely that \textit{any planar drawing of the maximal planar graph $G$ can be stretched}.
In other words, given a planar drawing of $G$, there is a stretched drawing of $G$ 
and a bijection between the bounded (necessarily triangular) regions such that corresponding triangles have the same edges of $G$ as the sides.

The rest of the proof remains unchanged.

\newpage
\phantomsection
\sloppy{\scriptsize\input{graph-theory.ind}}

\def\emph{\textit}
\printbibliography[heading=bibintoc]
\fussy

\end{document}